\input amstex

%

\def\next{AMS-SEKR}\ifx\styname\next \endinput\fi
\catcode`\@=11
\def\styname{AMS-SEKR}
\def\styversion{2.0}
{\W@{}\W@{\styname.STY - Version \styversion}\W@{}}
\hyphenation{acad-e-my acad-e-mies af-ter-thought anom-aly anom-alies
an-ti-deriv-a-tive an-tin-o-my an-tin-o-mies apoth-e-o-ses apoth-e-o-sis
ap-pen-dix ar-che-typ-al as-sign-a-ble as-sist-ant-ship as-ymp-tot-ic
asyn-chro-nous at-trib-uted at-trib-ut-able bank-rupt bank-rupt-cy
bi-dif-fer-en-tial blue-print busier busiest cat-a-stroph-ic
cat-a-stroph-i-cally con-gress cross-hatched data-base de-fin-i-tive
de-riv-a-tive dis-trib-ute dri-ver dri-vers eco-nom-ics econ-o-mist
elit-ist equi-vari-ant ex-quis-ite ex-tra-or-di-nary flow-chart
for-mi-da-ble forth-right friv-o-lous ge-o-des-ic ge-o-det-ic geo-met-ric
griev-ance griev-ous griev-ous-ly hexa-dec-i-mal ho-lo-no-my ho-mo-thetic
ideals idio-syn-crasy in-fin-ite-ly in-fin-i-tes-i-mal ir-rev-o-ca-ble
key-stroke lam-en-ta-ble light-weight mal-a-prop-ism man-u-script
mar-gin-al meta-bol-ic me-tab-o-lism meta-lan-guage me-trop-o-lis
met-ro-pol-i-tan mi-nut-est mol-e-cule mono-chrome mono-pole mo-nop-oly
mono-spline mo-not-o-nous mul-ti-fac-eted mul-ti-plic-able non-euclid-ean
non-iso-mor-phic non-smooth par-a-digm par-a-bol-ic pa-rab-o-loid
pa-ram-e-trize para-mount pen-ta-gon phe-nom-e-non post-script pre-am-ble
pro-ce-dur-al pro-hib-i-tive pro-hib-i-tive-ly pseu-do-dif-fer-en-tial
pseu-do-fi-nite pseu-do-nym qua-drat-ics quad-ra-ture qua-si-smooth
qua-si-sta-tion-ary qua-si-tri-an-gu-lar quin-tes-sence quin-tes-sen-tial
re-arrange-ment rec-tan-gle ret-ri-bu-tion retro-fit retro-fit-ted
right-eous right-eous-ness ro-bot ro-bot-ics sched-ul-ing se-mes-ter
semi-def-i-nite semi-ho-mo-thet-ic set-up se-vere-ly side-step sov-er-eign
spe-cious spher-oid spher-oid-al star-tling star-tling-ly
sta-tis-tics sto-chas-tic straight-est strange-ness strat-a-gem strong-hold
sum-ma-ble symp-to-matic syn-chro-nous topo-graph-i-cal tra-vers-a-ble
tra-ver-sal tra-ver-sals treach-ery turn-around un-at-tached un-err-ing-ly
white-space wide-spread wing-spread wretch-ed wretch-ed-ly Brown-ian
Eng-lish Euler-ian Feb-ru-ary Gauss-ian Grothen-dieck Hamil-ton-ian
Her-mit-ian Jan-u-ary Japan-ese Kor-te-weg Le-gendre Lip-schitz
Lip-schitz-ian Mar-kov-ian Noe-ther-ian No-vem-ber Rie-mann-ian
Schwarz-schild Sep-tem-ber
form per-iods Uni-ver-si-ty cri-ti-sism for-ma-lism}
\Invalid@\nofrills
\Invalid@\usualspace
\newif\ifnofrills@
\def\nofrills@#1#2{\relaxnext@
  \DN@{\ifx\next\nofrills
    \nofrills@true\let#2\relax\DN@\nofrills{\nextii@}%
  \else
    \nofrills@false\def#2{#1}\let\next@\nextii@\fi
\next@}}
\def\usualspace@#1{\ifnofrills@\def\usualspace{#1}\fi}
\def\addto#1#2{\csname \expandafter\eat@\string#1@\endcsname
  \expandafter{\the\csname \expandafter\eat@\string#1@\endcsname#2}}
\newdimen\bigsize@
\def\big@#1#2{{\hbox{$\left#2\vcenter to#1\bigsize@{}%
  \right.\nulldelimiterspace\z@\m@th$}}}
\def\big{\big@\@ne}
\def\Big{\big@{1.5}}
\def\bigg{\big@\tw@}
\def\Bigg{\big@{2.5}}
\def\raggedcenter@{\leftskip\z@ plus.4\hsize \rightskip\leftskip
 \parfillskip\z@ \parindent\z@ \spaceskip.3333em \xspaceskip.5em
 \pretolerance9999\tolerance9999 \exhyphenpenalty\@M
 \hyphenpenalty\@M \let\\\linebreak}
\def\upperspecialchars{\def\ss{SS}\let\i=I\let\j=J\let\ae\AE\let\oe\OE
  \let\o\O\let\aa\AA\let\l\L}
\def\uppercasetext@#1{%
  {\spaceskip1.2\fontdimen2\the\font plus1.2\fontdimen3\the\font
   \upperspecialchars\uctext@#1$\m@th\aftergroup\eat@$}}
\def\uctext@#1$#2${\endash@#1-\endash@$#2$\uctext@}
\def\endash@#1-#2\endash@{\uppercase{#1}\if\notempty{#2}--\endash@#2\endash@\fi}
\def\runaway@#1{\DN@{#1}\ifx\envir@\next@
  \Err@{You seem to have a missing or misspelled \string\end#1 ...}%
  \let\envir@\empty\fi}
\newif\iftemp@
\def\notempty#1{TT\fi\def\test@{#1}\ifx\test@\empty\temp@false
  \else\temp@true\fi \iftemp@}
\font@\tensmc=cmcsc10
\font@\sevenex=cmex7
\font@\sevenit=cmti7
\font@\eightrm=cmr8 
\font@\sixrm=cmr6 
\font@\eighti=cmmi8     \skewchar\eighti='177 
\font@\sixi=cmmi6       \skewchar\sixi='177   
\font@\eightsy=cmsy8    \skewchar\eightsy='60 
\font@\sixsy=cmsy6      \skewchar\sixsy='60   
\font@\eightex=cmex8
\font@\eightbf=cmbx8 
\font@\sixbf=cmbx6   
\font@\eightit=cmti8 
\font@\eightsl=cmsl8 
\font@\eightsmc=cmcsc8
\font@\eighttt=cmtt8 


\loadmsam
\loadmsbm
\loadeufm
\UseAMSsymbols
\newtoks\tenpoint@
\def\tenpoint{\normalbaselineskip12\p@
 \abovedisplayskip12\p@ plus3\p@ minus9\p@
 \belowdisplayskip\abovedisplayskip
 \abovedisplayshortskip\z@ plus3\p@
 \belowdisplayshortskip7\p@ plus3\p@ minus4\p@
 \textonlyfont@\rm\tenrm \textonlyfont@\it\tenit
 \textonlyfont@\sl\tensl \textonlyfont@\bf\tenbf
 \textonlyfont@\smc\tensmc \textonlyfont@\tt\tentt
 \textonlyfont@\bsmc\tenbsmc
 \ifsyntax@ \def\big##1{{\hbox{$\left##1\right.$}}}%
  \let\Big\big \let\bigg\big \let\Bigg\big
 \else
  \textfont\z@=\tenrm  \scriptfont\z@=\sevenrm  \scriptscriptfont\z@=\fiverm
  \textfont\@ne=\teni  \scriptfont\@ne=\seveni  \scriptscriptfont\@ne=\fivei
  \textfont\tw@=\tensy \scriptfont\tw@=\sevensy \scriptscriptfont\tw@=\fivesy
  \textfont\thr@@=\tenex \scriptfont\thr@@=\sevenex
        \scriptscriptfont\thr@@=\sevenex
  \textfont\itfam=\tenit \scriptfont\itfam=\sevenit
        \scriptscriptfont\itfam=\sevenit
  \textfont\bffam=\tenbf \scriptfont\bffam=\sevenbf
        \scriptscriptfont\bffam=\fivebf
  \setbox\strutbox\hbox{\vrule height8.5\p@ depth3.5\p@ width\z@}%
  \setbox\strutbox@\hbox{\lower.5\normallineskiplimit\vbox{%
        \kern-\normallineskiplimit\copy\strutbox}}%
 \setbox\z@\vbox{\hbox{$($}\kern\z@}\bigsize@=1.2\ht\z@
 \fi
 \normalbaselines\rm\ex@.2326ex\jot3\ex@\the\tenpoint@}
\newtoks\eightpoint@
\def\eightpoint{\normalbaselineskip10\p@
 \abovedisplayskip10\p@ plus2.4\p@ minus7.2\p@
 \belowdisplayskip\abovedisplayskip
 \abovedisplayshortskip\z@ plus2.4\p@
 \belowdisplayshortskip5.6\p@ plus2.4\p@ minus3.2\p@
 \textonlyfont@\rm\eightrm \textonlyfont@\it\eightit
 \textonlyfont@\sl\eightsl \textonlyfont@\bf\eightbf
 \textonlyfont@\smc\eightsmc \textonlyfont@\tt\eighttt
 \textonlyfont@\bsmc\eightbsmc
 \ifsyntax@\def\big##1{{\hbox{$\left##1\right.$}}}%
  \let\Big\big \let\bigg\big \let\Bigg\big
 \else
  \textfont\z@=\eightrm \scriptfont\z@=\sixrm \scriptscriptfont\z@=\fiverm
  \textfont\@ne=\eighti \scriptfont\@ne=\sixi \scriptscriptfont\@ne=\fivei
  \textfont\tw@=\eightsy \scriptfont\tw@=\sixsy \scriptscriptfont\tw@=\fivesy
  \textfont\thr@@=\eightex \scriptfont\thr@@=\sevenex
   \scriptscriptfont\thr@@=\sevenex
  \textfont\itfam=\eightit \scriptfont\itfam=\sevenit
   \scriptscriptfont\itfam=\sevenit
  \textfont\bffam=\eightbf \scriptfont\bffam=\sixbf
   \scriptscriptfont\bffam=\fivebf
 \setbox\strutbox\hbox{\vrule height7\p@ depth3\p@ width\z@}%
 \setbox\strutbox@\hbox{\raise.5\normallineskiplimit\vbox{%
   \kern-\normallineskiplimit\copy\strutbox}}%
 \setbox\z@\vbox{\hbox{$($}\kern\z@}\bigsize@=1.2\ht\z@
 \fi
 \normalbaselines\eightrm\ex@.2326ex\jot3\ex@\the\eightpoint@}

\font@\twelverm=cmr10 scaled\magstep1
\font@\twelveit=cmti10 scaled\magstep1
\font@\twelvesl=cmsl10 scaled\magstep1
\font@\twelvesmc=cmcsc10 scaled\magstep1
\font@\twelvett=cmtt10 scaled\magstep1
\font@\twelvebf=cmbx10 scaled\magstep1
\font@\twelvei=cmmi10 scaled\magstep1
\font@\twelvesy=cmsy10 scaled\magstep1
\font@\twelveex=cmex10 scaled\magstep1
\font@\twelvemsa=msam10 scaled\magstep1
\font@\twelveeufm=eufm10 scaled\magstep1
\font@\twelvemsb=msbm10 scaled\magstep1
\newtoks\twelvepoint@
\def\twelvepoint{\normalbaselineskip15\p@
 \abovedisplayskip15\p@ plus3.6\p@ minus10.8\p@
 \belowdisplayskip\abovedisplayskip
 \abovedisplayshortskip\z@ plus3.6\p@
 \belowdisplayshortskip8.4\p@ plus3.6\p@ minus4.8\p@
 \textonlyfont@\rm\twelverm \textonlyfont@\it\twelveit
 \textonlyfont@\sl\twelvesl \textonlyfont@\bf\twelvebf
 \textonlyfont@\smc\twelvesmc \textonlyfont@\tt\twelvett
 \textonlyfont@\bsmc\twelvebsmc
 \ifsyntax@ \def\big##1{{\hbox{$\left##1\right.$}}}%
  \let\Big\big \let\bigg\big \let\Bigg\big
 \else
  \textfont\z@=\twelverm  \scriptfont\z@=\tenrm  \scriptscriptfont\z@=\sevenrm
  \textfont\@ne=\twelvei  \scriptfont\@ne=\teni  \scriptscriptfont\@ne=\seveni
  \textfont\tw@=\twelvesy \scriptfont\tw@=\tensy \scriptscriptfont\tw@=\sevensy
  \textfont\thr@@=\twelveex \scriptfont\thr@@=\tenex
        \scriptscriptfont\thr@@=\tenex
  \textfont\itfam=\twelveit \scriptfont\itfam=\tenit
        \scriptscriptfont\itfam=\tenit
  \textfont\bffam=\twelvebf \scriptfont\bffam=\tenbf
        \scriptscriptfont\bffam=\sevenbf
  \setbox\strutbox\hbox{\vrule height10.2\p@ depth4.2\p@ width\z@}%
  \setbox\strutbox@\hbox{\lower.6\normallineskiplimit\vbox{%
        \kern-\normallineskiplimit\copy\strutbox}}%
 \setbox\z@\vbox{\hbox{$($}\kern\z@}\bigsize@=1.4\ht\z@
 \fi
 \normalbaselines\rm\ex@.2326ex\jot3.6\ex@\the\twelvepoint@}

\def\headfonts{\twelvepoint\bf}

\font@\fourteenrm=cmr10 scaled\magstep2
\font@\fourteenit=cmti10 scaled\magstep2
\font@\fourteensl=cmsl10 scaled\magstep2
\font@\fourteensmc=cmcsc10 scaled\magstep2
\font@\fourteentt=cmtt10 scaled\magstep2
\font@\fourteenbf=cmbx10 scaled\magstep2
\font@\fourteeni=cmmi10 scaled\magstep2
\font@\fourteensy=cmsy10 scaled\magstep2
\font@\fourteenex=cmex10 scaled\magstep2
\font@\fourteenmsa=msam10 scaled\magstep2
\font@\fourteeneufm=eufm10 scaled\magstep2
\font@\fourteenmsb=msbm10 scaled\magstep2
\newtoks\fourteenpoint@
\def\fourteenpoint{\normalbaselineskip15\p@
 \abovedisplayskip18\p@ plus4.3\p@ minus12.9\p@
 \belowdisplayskip\abovedisplayskip
 \abovedisplayshortskip\z@ plus4.3\p@
 \belowdisplayshortskip10.1\p@ plus4.3\p@ minus5.8\p@
 \textonlyfont@\rm\fourteenrm \textonlyfont@\it\fourteenit
 \textonlyfont@\sl\fourteensl \textonlyfont@\bf\fourteenbf
 \textonlyfont@\smc\fourteensmc \textonlyfont@\tt\fourteentt
 \textonlyfont@\bsmc\fourteenbsmc
 \ifsyntax@ \def\big##1{{\hbox{$\left##1\right.$}}}%
  \let\Big\big \let\bigg\big \let\Bigg\big
 \else
  \textfont\z@=\fourteenrm  \scriptfont\z@=\twelverm  \scriptscriptfont\z@=\tenrm
  \textfont\@ne=\fourteeni  \scriptfont\@ne=\twelvei  \scriptscriptfont\@ne=\teni
  \textfont\tw@=\fourteensy \scriptfont\tw@=\twelvesy \scriptscriptfont\tw@=\tensy
  \textfont\thr@@=\fourteenex \scriptfont\thr@@=\twelveex
        \scriptscriptfont\thr@@=\twelveex
  \textfont\itfam=\fourteenit \scriptfont\itfam=\twelveit
        \scriptscriptfont\itfam=\twelveit
  \textfont\bffam=\fourteenbf \scriptfont\bffam=\twelvebf
        \scriptscriptfont\bffam=\tenbf
  \setbox\strutbox\hbox{\vrule height12.2\p@ depth5\p@ width\z@}%
  \setbox\strutbox@\hbox{\lower.72\normallineskiplimit\vbox{%
        \kern-\normallineskiplimit\copy\strutbox}}%
 \setbox\z@\vbox{\hbox{$($}\kern\z@}\bigsize@=1.7\ht\z@
 \fi
 \normalbaselines\rm\ex@.2326ex\jot4.3\ex@\the\fourteenpoint@}

\def\chapheadfonts{\fourteenpoint\bf}

\font@\seventeenrm=cmr10 scaled\magstep3
\font@\seventeenit=cmti10 scaled\magstep3
\font@\seventeensl=cmsl10 scaled\magstep3
\font@\seventeensmc=cmcsc10 scaled\magstep3
\font@\seventeentt=cmtt10 scaled\magstep3
\font@\seventeenbf=cmbx10 scaled\magstep3
\font@\seventeeni=cmmi10 scaled\magstep3
\font@\seventeensy=cmsy10 scaled\magstep3
\font@\seventeenex=cmex10 scaled\magstep3
\font@\seventeenmsa=msam10 scaled\magstep3
\font@\seventeeneufm=eufm10 scaled\magstep3
\font@\seventeenmsb=msbm10 scaled\magstep3
\newtoks\seventeenpoint@
\def\seventeenpoint{\normalbaselineskip18\p@
 \abovedisplayskip21.6\p@ plus5.2\p@ minus15.4\p@
 \belowdisplayskip\abovedisplayskip
 \abovedisplayshortskip\z@ plus5.2\p@
 \belowdisplayshortskip12.1\p@ plus5.2\p@ minus7\p@
 \textonlyfont@\rm\seventeenrm \textonlyfont@\it\seventeenit
 \textonlyfont@\sl\seventeensl \textonlyfont@\bf\seventeenbf
 \textonlyfont@\smc\seventeensmc \textonlyfont@\tt\seventeentt
 \textonlyfont@\bsmc\seventeenbsmc
 \ifsyntax@ \def\big##1{{\hbox{$\left##1\right.$}}}%
  \let\Big\big \let\bigg\big \let\Bigg\big
 \else
  \textfont\z@=\seventeenrm  \scriptfont\z@=\fourteenrm  \scriptscriptfont\z@=\twelverm
  \textfont\@ne=\seventeeni  \scriptfont\@ne=\fourteeni  \scriptscriptfont\@ne=\twelvei
  \textfont\tw@=\seventeensy \scriptfont\tw@=\fourteensy \scriptscriptfont\tw@=\twelvesy
  \textfont\thr@@=\seventeenex \scriptfont\thr@@=\fourteenex
        \scriptscriptfont\thr@@=\fourteenex
  \textfont\itfam=\seventeenit \scriptfont\itfam=\fourteenit
        \scriptscriptfont\itfam=\fourteenit
  \textfont\bffam=\seventeenbf \scriptfont\bffam=\fourteenbf
        \scriptscriptfont\bffam=\twelvebf
  \setbox\strutbox\hbox{\vrule height14.6\p@ depth6\p@ width\z@}%
  \setbox\strutbox@\hbox{\lower.86\normallineskiplimit\vbox{%
        \kern-\normallineskiplimit\copy\strutbox}}%
 \setbox\z@\vbox{\hbox{$($}\kern\z@}\bigsize@=2\ht\z@
 \fi
 \normalbaselines\rm\ex@.2326ex\jot5.2\ex@\the\seventeenpoint@}

\font@\rrrrrm=cmr10 scaled\magstep4
\font@\bigtitlefont=cmbx10 scaled\magstep4

\parindent1pc
\normallineskiplimit\p@
\newdimen\indenti \indenti=2pc
\def\pageheight#1{\vsize#1}
\def\pagewidth#1{\hsize#1%
   \captionwidth@\hsize \advance\captionwidth@-2\indenti}
\pagewidth{30pc} \pageheight{47pc}
\def\topmatter{%
 \ifx\undefined\msafam
 \else\font@\eightmsa=msam8 \font@\sixmsa=msam6
   \ifsyntax@\else \addto\tenpoint{\textfont\msafam=\tenmsa
              \scriptfont\msafam=\sevenmsa \scriptscriptfont\msafam=\fivemsa}%
     \addto\eightpoint{\textfont\msafam=\eightmsa \scriptfont\msafam=\sixmsa
              \scriptscriptfont\msafam=\fivemsa}%
   \fi
 \fi
 \ifx\undefined\msbfam
 \else\font@\eightmsb=msbm8 \font@\sixmsb=msbm6
   \ifsyntax@\else \addto\tenpoint{\textfont\msbfam=\tenmsb
         \scriptfont\msbfam=\sevenmsb \scriptscriptfont\msbfam=\fivemsb}%
     \addto\eightpoint{\textfont\msbfam=\eightmsb \scriptfont\msbfam=\sixmsb
         \scriptscriptfont\msbfam=\fivemsb}%
   \fi
 \fi
 \ifx\undefined\eufmfam
 \else \font@\eighteufm=eufm8 \font@\sixeufm=eufm6
   \ifsyntax@\else \addto\tenpoint{\textfont\eufmfam=\teneufm
       \scriptfont\eufmfam=\seveneufm \scriptscriptfont\eufmfam=\fiveeufm}%
     \addto\eightpoint{\textfont\eufmfam=\eighteufm
       \scriptfont\eufmfam=\sixeufm \scriptscriptfont\eufmfam=\fiveeufm}%
   \fi
 \fi
 \ifx\undefined\eufbfam
 \else \font@\eighteufb=eufb8 \font@\sixeufb=eufb6
   \ifsyntax@\else \addto\tenpoint{\textfont\eufbfam=\teneufb
      \scriptfont\eufbfam=\seveneufb \scriptscriptfont\eufbfam=\fiveeufb}%
    \addto\eightpoint{\textfont\eufbfam=\eighteufb
      \scriptfont\eufbfam=\sixeufb \scriptscriptfont\eufbfam=\fiveeufb}%
   \fi
 \fi
 \ifx\undefined\eusmfam
 \else \font@\eighteusm=eusm8 \font@\sixeusm=eusm6
   \ifsyntax@\else \addto\tenpoint{\textfont\eusmfam=\teneusm
       \scriptfont\eusmfam=\seveneusm \scriptscriptfont\eusmfam=\fiveeusm}%
     \addto\eightpoint{\textfont\eusmfam=\eighteusm
       \scriptfont\eusmfam=\sixeusm \scriptscriptfont\eusmfam=\fiveeusm}%
   \fi
 \fi
 \ifx\undefined\eusbfam
 \else \font@\eighteusb=eusb8 \font@\sixeusb=eusb6
   \ifsyntax@\else \addto\tenpoint{\textfont\eusbfam=\teneusb
       \scriptfont\eusbfam=\seveneusb \scriptscriptfont\eusbfam=\fiveeusb}%
     \addto\eightpoint{\textfont\eusbfam=\eighteusb
       \scriptfont\eusbfam=\sixeusb \scriptscriptfont\eusbfam=\fiveeusb}%
   \fi
 \fi
 \ifx\undefined\eurmfam
 \else \font@\eighteurm=eurm8 \font@\sixeurm=eurm6
   \ifsyntax@\else \addto\tenpoint{\textfont\eurmfam=\teneurm
       \scriptfont\eurmfam=\seveneurm \scriptscriptfont\eurmfam=\fiveeurm}%
     \addto\eightpoint{\textfont\eurmfam=\eighteurm
       \scriptfont\eurmfam=\sixeurm \scriptscriptfont\eurmfam=\fiveeurm}%
   \fi
 \fi
 \ifx\undefined\eurbfam
 \else \font@\eighteurb=eurb8 \font@\sixeurb=eurb6
   \ifsyntax@\else \addto\tenpoint{\textfont\eurbfam=\teneurb
       \scriptfont\eurbfam=\seveneurb \scriptscriptfont\eurbfam=\fiveeurb}%
    \addto\eightpoint{\textfont\eurbfam=\eighteurb
       \scriptfont\eurbfam=\sixeurb \scriptscriptfont\eurbfam=\fiveeurb}%
   \fi
 \fi
 \ifx\undefined\cmmibfam
 \else \font@\eightcmmib=cmmib8 \font@\sixcmmib=cmmib6
   \ifsyntax@\else \addto\tenpoint{\textfont\cmmibfam=\tencmmib
       \scriptfont\cmmibfam=\sevencmmib \scriptscriptfont\cmmibfam=\fivecmmib}%
    \addto\eightpoint{\textfont\cmmibfam=\eightcmmib
       \scriptfont\cmmibfam=\sixcmmib \scriptscriptfont\cmmibfam=\fivecmmib}%
   \fi
 \fi
 \ifx\undefined\cmbsyfam
 \else \font@\eightcmbsy=cmbsy8 \font@\sixcmbsy=cmbsy6
   \ifsyntax@\else \addto\tenpoint{\textfont\cmbsyfam=\tencmbsy
      \scriptfont\cmbsyfam=\sevencmbsy \scriptscriptfont\cmbsyfam=\fivecmbsy}%
    \addto\eightpoint{\textfont\cmbsyfam=\eightcmbsy
      \scriptfont\cmbsyfam=\sixcmbsy \scriptscriptfont\cmbsyfam=\fivecmbsy}%
   \fi
 \fi
 \let\topmatter\relax}
\def\chapterno@{\uppercase\expandafter{\romannumeral\chaptercount@}}
\newcount\chaptercount@
\def\chapter{\nofrills@{\afterassignment\chapterno@
                        CHAPTER \global\chaptercount@=}\chapter@
 \DNii@##1{\leavevmode\hskip-\leftskip
   \rlap{\vbox to\z@{\vss\centerline{\eightpoint
   \chapter@##1\unskip}\baselineskip2pc\null}}\hskip\leftskip
   \nofrills@false}%
 \FN@\next@}
\newbox\titlebox@

\def\title{\nofrills@{\relax}\title@%
 \DNii@##1\endtitle{\global\setbox\titlebox@\vtop{\tenpoint\bf
 \raggedcenter@\ignorespaces
 \baselineskip1.3\baselineskip\title@{##1}\endgraf}%
 \ifmonograph@ \edef\next{\the\leftheadtoks}\ifx\next\empty
    \leftheadtext{##1}\fi
 \fi
 \edef\next{\the\rightheadtoks}\ifx\next\empty \rightheadtext{##1}\fi
 }\FN@\next@}
\newbox\authorbox@
\def\author#1\endauthor{\global\setbox\authorbox@
 \vbox{\tenpoint\smc\raggedcenter@\ignorespaces
 #1\endgraf}\relaxnext@ \edef\next{\the\leftheadtoks}%
 \ifx\next\empty\leftheadtext{#1}\fi}
\newbox\affilbox@
\def\affil#1\endaffil{\global\setbox\affilbox@
 \vbox{\tenpoint\raggedcenter@\ignorespaces#1\endgraf}}
\newcount\addresscount@
\addresscount@\z@
\def\address#1\endaddress{\global\advance\addresscount@\@ne
  \expandafter\gdef\csname address\number\addresscount@\endcsname
  {\vskip12\p@ minus6\p@\noindent\eightpoint\smc\ignorespaces#1\par}}
\def\email{\nofrills@{\eightpoint{\it E-mail\/}:\enspace}\email@
  \DNii@##1\endemail{%
  \expandafter\gdef\csname email\number\addresscount@\endcsname
  {\def\usualspace{{\it\enspace}}\smallskip\noindent\eightpoint\email@
  \ignorespaces##1\par}}%
 \FN@\next@}
\def\thedate@{}
\def\date#1\enddate{\gdef\thedate@{\tenpoint\ignorespaces#1\unskip}}
\def\thethanks@{}
\def\thanks#1\endthanks{\gdef\thethanks@{\eightpoint\ignorespaces#1.\unskip}}
\def\thekeywords@{}
\def\keywords{\nofrills@{{\it Key words and phrases.\enspace}}\keywords@
 \DNii@##1\endkeywords{\def\thekeywords@{\def\usualspace{{\it\enspace}}%
 \eightpoint\keywords@\ignorespaces##1\unskip.}}%
 \FN@\next@}
\def\thesubjclass@{}
\def\subjclass{\nofrills@{{\rm2000 {\it Mathematics Subject
   Classification\/}.\enspace}}\subjclass@
 \DNii@##1\endsubjclass{\def\thesubjclass@{\def\usualspace
  {{\rm\enspace}}\eightpoint\subjclass@\ignorespaces##1\unskip.}}%
 \FN@\next@}
\newbox\abstractbox@
\def\abstract{\nofrills@{{\smc Abstract.\enspace}}\abstract@
 \DNii@{\setbox\abstractbox@\vbox\bgroup\noindent$$\vbox\bgroup
  \def\envir@{abstract}\advance\hsize-2\indenti
  \usualspace@{{\enspace}}\eightpoint \noindent\abstract@\ignorespaces}%
 \FN@\next@}
\def\endabstract{\par\unskip\egroup$$\egroup}
\def\widestnumber#1#2{\begingroup\let\head\null\let\subhead\empty
   \let\subsubhead\subhead
   \ifx#1\head\global\setbox\tocheadbox@\hbox{#2.\enspace}%
   \else\ifx#1\subhead\global\setbox\tocsubheadbox@\hbox{#2.\enspace}%
   \else\ifx#1\key\bgroup\let\endrefitem@\egroup
        \key#2\endrefitem@\global\refindentwd\wd\keybox@
   \else\ifx#1\no\bgroup\let\endrefitem@\egroup
        \no#2\endrefitem@\global\refindentwd\wd\nobox@
   \else\ifx#1\page\global\setbox\pagesbox@\hbox{\quad\bf#2}%
   \else\ifx#1\item\setboxz@h{#2}\global\rosteritemwd\wdz@
        \global\advance\rosteritemwd by.5\parindent
   \else\message{\string\widestnumber is not defined for this option
   (\string#1)}%
\fi\fi\fi\fi\fi\fi\endgroup}
\newif\ifmonograph@
\def\Monograph{\monograph@true \let\headmark\rightheadtext
  \let\varindent@\indent \def\headfont@{\bf}\def\proclaimheadfont@{\smc}%
  \def\demofont@{\smc}}
\let\varindent@\indent

\newbox\tocheadbox@    \newbox\tocsubheadbox@
\newbox\tocbox@
\def\toc{\toc@{Contents}}
\def\newtocdefs{%
   \def \title##1\endtitle
       {\penaltyandskip@\z@\smallskipamount
        \hangindent\wd\tocheadbox@\noindent{\bf##1}}%
   \def \chapter##1{%
        Chapter \uppercase\expandafter{\romannumeral##1.\unskip}\enspace}%
   \def \specialhead##1\endspecialhead
       {\par\hangindent\wd\tocheadbox@ \noindent##1\par}%
   \def \head##1 ##2\endhead
       {\par\hangindent\wd\tocheadbox@ \noindent
        \if\notempty{##1}\hbox to\wd\tocheadbox@{\hfil##1\enspace}\fi
        ##2\par}%
   \def \subhead##1 ##2\endsubhead
       {\par\vskip-\parskip {\normalbaselines
        \advance\leftskip\wd\tocheadbox@
        \hangindent\wd\tocsubheadbox@ \noindent
        \if\notempty{##1}\hbox to\wd\tocsubheadbox@{##1\unskip\hfil}\fi
         ##2\par}}%
   \def \subsubhead##1 ##2\endsubsubhead
       {\par\vskip-\parskip {\normalbaselines
        \advance\leftskip\wd\tocheadbox@
        \hangindent\wd\tocsubheadbox@ \noindent
        \if\notempty{##1}\hbox to\wd\tocsubheadbox@{##1\unskip\hfil}\fi
        ##2\par}}}
\def\toc@#1{\relaxnext@
   \def\page##1%
       {\unskip\penalty0\null\hfil
        \rlap{\hbox to\wd\pagesbox@{\quad\hfil##1}}\hfilneg\penalty\@M}%
 \DN@{\ifx\next\nofrills\DN@\nofrills{\nextii@}%
      \else\DN@{\nextii@{{#1}}}\fi
      \next@}%
 \DNii@##1{%
\ifmonograph@\bgroup\else\setbox\tocbox@\vbox\bgroup
   \centerline{\headfont@\ignorespaces##1\unskip}\nobreak
   \vskip\belowheadskip \fi
   \setbox\tocheadbox@\hbox{0.\enspace}%
   \setbox\tocsubheadbox@\hbox{0.0.\enspace}%
   \leftskip\indenti \rightskip\leftskip
   \setbox\pagesbox@\hbox{\bf\quad000}\advance\rightskip\wd\pagesbox@
   \newtocdefs
 }%
 \FN@\next@}
\def\endtoc{\par\egroup}
\let\pretitle\relax
\let\preauthor\relax
\let\preaffil\relax
\let\predate\relax
\let\preabstract\relax
\let\prepaper\relax
\def\dedicatory #1\enddedicatory{\def\preabstract{{\medskip
  \eightpoint\it \raggedcenter@#1\endgraf}}}
\def\thetranslator@{}
\def\translator#1\endtranslator{\def\thetranslator@{\nobreak\medskip
 \line{\eightpoint\hfil Translated by \uppercase{#1}\qquad\qquad}\nobreak}}
\outer\def\endtopmatter{\runaway@{abstract}%
 \edef\next{\the\leftheadtoks}\ifx\next\empty
  \expandafter\leftheadtext\expandafter{\the\rightheadtoks}\fi
 \ifmonograph@\else
   \ifx\thesubjclass@\empty\else \makefootnote@{}{\thesubjclass@}\fi
   \ifx\thekeywords@\empty\else \makefootnote@{}{\thekeywords@}\fi
   \ifx\thethanks@\empty\else \makefootnote@{}{\thethanks@}\fi
 \fi
  \pretitle
  \ifmonograph@ \topskip7pc \else \topskip4pc \fi
  \box\titlebox@
  \topskip10pt
  \preauthor
  \ifvoid\authorbox@\else \vskip2.5pc plus1pc \unvbox\authorbox@\fi
  \preaffil
  \ifvoid\affilbox@\else \vskip1pc plus.5pc \unvbox\affilbox@\fi
  \predate
  \ifx\thedate@\empty\else \vskip1pc plus.5pc \line{\hfil\thedate@\hfil}\fi
  \preabstract
  \ifvoid\abstractbox@\else \vskip1.5pc plus.5pc \unvbox\abstractbox@ \fi
  \ifvoid\tocbox@\else\vskip1.5pc plus.5pc \unvbox\tocbox@\fi
  \prepaper
  \vskip2pc plus1pc
}
\def\document{\let\fontlist@\relax\let\alloclist@\relax
  \tenpoint}

\newskip\aboveheadskip       \aboveheadskip1.8\bigskipamount
\newdimen\belowheadskip      \belowheadskip1.8\medskipamount

\def\headfont@{\smc}
\def\penaltyandskip@#1#2{\relax\ifdim\lastskip<#2\relax\removelastskip
      \ifnum#1=\z@\else\penalty@#1\relax\fi\vskip#2%
  \else\ifnum#1=\z@\else\penalty@#1\relax\fi\fi}
\def\nobreak{\penalty\@M
  \ifvmode\def\penalty@{\let\penalty@\penalty\count@@@}%
  \everypar{\let\penalty@\penalty\everypar{}}\fi}
\let\penalty@\penalty
\def\heading#1\endheading{\head#1\endhead}

\def\specialheadfont@{\bf}
\outer\def\specialhead{\par\penaltyandskip@{-200}\aboveheadskip
  \begingroup\interlinepenalty\@M\rightskip\z@ plus\hsize \let\\\linebreak
  \specialheadfont@\noindent\ignorespaces}
\def\endspecialhead{\par\endgroup\nobreak\vskip\belowheadskip}
\let\headmark\eat@
\newskip\subheadskip       \subheadskip\medskipamount
\def\subheadfont@{\bf}
\outer\def\subhead{\nofrills@{.\enspace}\subhead@
 \DNii@##1\endsubhead{\par\penaltyandskip@{-100}\subheadskip
  \varindent@{\usualspace@{{\subheadfont@\enspace}}%
 \subheadfont@\ignorespaces##1\unskip\subhead@}\ignorespaces}%
 \FN@\next@}
\outer\def\subsubhead{\nofrills@{.\enspace}\subsubhead@
 \DNii@##1\endsubsubhead{\par\penaltyandskip@{-50}\medskipamount
      {\usualspace@{{\it\enspace}}%
  \it\ignorespaces##1\unskip\subsubhead@}\ignorespaces}%
 \FN@\next@}
\def\proclaimheadfont@{\bf}
\outer\def\proclaim{\runaway@{proclaim}\def\envir@{proclaim}%
  \nofrills@{.\enspace}\proclaim@
 \DNii@##1{\penaltyandskip@{-100}\medskipamount\varindent@
   \usualspace@{{\proclaimheadfont@\enspace}}\proclaimheadfont@
   \ignorespaces##1\unskip\proclaim@
  \sl\ignorespaces}%
 \FN@\next@}
\outer\def\endproclaim{\let\envir@\relax\par\rm
  \penaltyandskip@{55}\medskipamount}
\def\demoheadfont@{\it}
\def\demo{\runaway@{proclaim}\nofrills@{.\enspace}\demo@
     \DNii@##1{\par\penaltyandskip@\z@\medskipamount
  {\usualspace@{{\demoheadfont@\enspace}}%
  \varindent@\demoheadfont@\ignorespaces##1\unskip\demo@}\rm
  \ignorespaces}\FN@\next@}
\def\enddemo{\par\medskip}
\def\qed{\ifhmode\unskip\nobreak\fi\quad\ifmmode\square\else$\m@th\square$\fi}
\let\remark\demo
\let\endremark\enddemo
\def\definition{\runaway@{proclaim}%
  \nofrills@{.\demoheadfont@\enspace}\definition@
        \DNii@##1{\penaltyandskip@{-100}\medskipamount
        {\usualspace@{{\demoheadfont@\enspace}}%
        \varindent@\demoheadfont@\ignorespaces##1\unskip\definition@}%
        \rm \ignorespaces}\FN@\next@}

\let\example\demo
\let\endexample\enddemo

\newdimen\rosteritemwd
\newcount\rostercount@
\newif\iffirstitem@
\let\plainitem@\item
\newtoks\everypartoks@
\def\par@{\everypartoks@\expandafter{\the\everypar}\everypar{}}
\def\roster{\edef\leftskip@{\leftskip\the\leftskip}%
 \relaxnext@
 \rostercount@\z@  
 \def\item{\FN@\rosteritem@}%
 \DN@{\ifx\next\runinitem\let\next@\nextii@\else
  \let\next@\nextiii@\fi\next@}%
 \DNii@\runinitem  
  {\unskip  
   \DN@{\ifx\next[\let\next@\nextii@\else
    \ifx\next"\let\next@\nextiii@\else\let\next@\nextiv@\fi\fi\next@}%
   \DNii@[####1]{\rostercount@####1\relax
    \enspace{\rm(\number\rostercount@)}~\ignorespaces}%
   \def\nextiii@"####1"{\enspace{\rm####1}~\ignorespaces}%
   \def\nextiv@{\enspace{\rm(1)}\rostercount@\@ne~}%
   \par@\firstitem@false  
   \FN@\next@}%
 \def\nextiii@{\par\par@  
  \penalty\@m\smallskip\vskip-\parskip
  \firstitem@true}%
 \FN@\next@}
\def\rosteritem@{\iffirstitem@\firstitem@false\else\par\vskip-\parskip\fi
 \leftskip3\parindent\noindent  
 \DNii@[##1]{\rostercount@##1\relax
  \llap{\hbox to2.5\parindent{\hss\rm(\number\rostercount@)}%
   \hskip.5\parindent}\ignorespaces}%
 \def\nextiii@"##1"{%
  \llap{\hbox to2.5\parindent{\hss\rm##1}\hskip.5\parindent}\ignorespaces}%
 \def\nextiv@{\advance\rostercount@\@ne
  \llap{\hbox to2.5\parindent{\hss\rm(\number\rostercount@)}%
   \hskip.5\parindent}}%
 \ifx\next[\let\next@\nextii@\else\ifx\next"\let\next@\nextiii@\else
  \let\next@\nextiv@\fi\fi\next@}

\newif\ifnextRunin@
\def\endroster{\relaxnext@
 \par\leftskip@  
 \penalty-50 \vskip-\parskip\smallskip  
 \DN@{\ifx\next\Runinitem\let\next@\relax
  \else\nextRunin@false\let\item\plainitem@  
   \ifx\next\par 
    \DN@\par{\everypar\expandafter{\the\everypartoks@}}%
   \else  
    \DN@{\noindent\everypar\expandafter{\the\everypartoks@}}%
  \fi\fi\next@}%
 \FN@\next@}
\newcount\rosterhangafter@
\def\Runinitem#1\roster\runinitem{\relaxnext@
 \rostercount@\z@ 
 \def\item{\FN@\rosteritem@}%
 \def\runinitem@{#1}%
 \DN@{\ifx\next[\let\next\nextii@\else\ifx\next"\let\next\nextiii@
  \else\let\next\nextiv@\fi\fi\next}%
 \DNii@[##1]{\rostercount@##1\relax
  \def\item@{{\rm(\number\rostercount@)}}\nextv@}%
 \def\nextiii@"##1"{\def\item@{{\rm##1}}\nextv@}%
 \def\nextiv@{\advance\rostercount@\@ne
  \def\item@{{\rm(\number\rostercount@)}}\nextv@}%
 \def\nextv@{\setbox\z@\vbox  
  {\ifnextRunin@\noindent\fi  
  \runinitem@\unskip\enspace\item@~\par  
  \global\rosterhangafter@\prevgraf}%
  \firstitem@false  
  \ifnextRunin@\else\par\fi  
  \hangafter\rosterhangafter@\hangindent3\parindent
  \ifnextRunin@\noindent\fi  
  \runinitem@\unskip\enspace 
  \item@~\ifnextRunin@\else\par@\fi  
  \nextRunin@true\ignorespaces}%
 \FN@\next@}
\def\footmarkform@#1{$\m@th^{#1}$}
\let\thefootnotemark\footmarkform@
\def\makefootnote@#1#2{\insert\footins
 {\interlinepenalty\interfootnotelinepenalty
 \eightpoint\splittopskip\ht\strutbox\splitmaxdepth\dp\strutbox
 \floatingpenalty\@MM\leftskip\z@\rightskip\z@\spaceskip\z@\xspaceskip\z@
 \leavevmode{#1}\footstrut\ignorespaces#2\unskip\lower\dp\strutbox
 \vbox to\dp\strutbox{}}}
\newcount\footmarkcount@
\footmarkcount@\z@
\def\footnotemark{\let\@sf\empty\relaxnext@
 \ifhmode\edef\@sf{\spacefactor\the\spacefactor}\/\fi
 \DN@{\ifx[\next\let\next@\nextii@\else
  \ifx"\next\let\next@\nextiii@\else
  \let\next@\nextiv@\fi\fi\next@}%
 \DNii@[##1]{\footmarkform@{##1}\@sf}%
 \def\nextiii@"##1"{{##1}\@sf}%
 \def\nextiv@{\iffirstchoice@\global\advance\footmarkcount@\@ne\fi
  \footmarkform@{\number\footmarkcount@}\@sf}%
 \FN@\next@}
\def\footnotetext{\relaxnext@
 \DN@{\ifx[\next\let\next@\nextii@\else
  \ifx"\next\let\next@\nextiii@\else
  \let\next@\nextiv@\fi\fi\next@}%
 \DNii@[##1]##2{\makefootnote@{\footmarkform@{##1}}{##2}}%
 \def\nextiii@"##1"##2{\makefootnote@{##1}{##2}}%
 \def\nextiv@##1{\makefootnote@{\footmarkform@{\number\footmarkcount@}}{##1}}%
 \FN@\next@}
\def\footnote{\let\@sf\empty\relaxnext@
 \ifhmode\edef\@sf{\spacefactor\the\spacefactor}\/\fi
 \DN@{\ifx[\next\let\next@\nextii@\else
  \ifx"\next\let\next@\nextiii@\else
  \let\next@\nextiv@\fi\fi\next@}%
 \DNii@[##1]##2{\footnotemark[##1]\footnotetext[##1]{##2}}%
 \def\nextiii@"##1"##2{\footnotemark"##1"\footnotetext"##1"{##2}}%
 \def\nextiv@##1{\footnotemark\footnotetext{##1}}%
 \FN@\next@}
\def\adjustfootnotemark#1{\advance\footmarkcount@#1\relax}
\def\footnoterule{\kern-3\p@
  \hrule width 5pc\kern 2.6\p@} 
\def\captionfont@{\smc}
\def\topcaption#1#2\endcaption{%
  {\dimen@\hsize \advance\dimen@-\captionwidth@
   \rm\raggedcenter@ \advance\leftskip.5\dimen@ \rightskip\leftskip
  {\captionfont@#1}%
  \if\notempty{#2}.\enspace\ignorespaces#2\fi
  \endgraf}\nobreak\bigskip}
\def\botcaption#1#2\endcaption{%
  \nobreak\bigskip
  \setboxz@h{\captionfont@#1\if\notempty{#2}.\enspace\rm#2\fi}%
  {\dimen@\hsize \advance\dimen@-\captionwidth@
   \leftskip.5\dimen@ \rightskip\leftskip
   \noindent \ifdim\wdz@>\captionwidth@ 
   \else\hfil\fi 
  {\captionfont@#1}\if\notempty{#2}.\enspace\rm#2\fi\endgraf}}
\def\@ins{\par\begingroup\def\vspace##1{\vskip##1\relax}%
  \def\captionwidth##1{\captionwidth@##1\relax}%
  \setbox\z@\vbox\bgroup} 
\def\block{\RIfMIfI@\nondmatherr@\block\fi
       \else\ifvmode\vskip\abovedisplayskip\noindent\fi
        $$\def\endblock{\par\egroup$$}\fi
  \vbox\bgroup\advance\hsize-2\indenti\noindent}
\def\endblock{\par\egroup}
\def\cite#1{{\rm[{\citefont@\m@th#1}]}}
\def\citefont@{\rm}
\def\refsfont@{\eightpoint}
\outer\def\Refs{\runaway@{proclaim}%
 \relaxnext@ \DN@{\ifx\next\nofrills\DN@\nofrills{\nextii@}\else
  \DN@{\nextii@{References}}\fi\next@}%
 \DNii@##1{\penaltyandskip@{-200}\aboveheadskip
  \line{\hfil\headfont@\ignorespaces##1\unskip\hfil}\nobreak
  \vskip\belowheadskip
  \begingroup\refsfont@\sfcode`.=\@m}%
 \FN@\next@}
\def\endRefs{\par\endgroup}
\newbox\nobox@            \newbox\keybox@           \newbox\bybox@
\newbox\paperbox@         \newbox\paperinfobox@     \newbox\jourbox@
\newbox\volbox@           \newbox\issuebox@         \newbox\yrbox@
\newbox\pagesbox@         \newbox\bookbox@          \newbox\bookinfobox@
\newbox\publbox@          \newbox\publaddrbox@      \newbox\finalinfobox@
\newbox\edsbox@           \newbox\langbox@
\newif\iffirstref@        \newif\iflastref@
\newif\ifprevjour@        \newif\ifbook@            \newif\ifprevinbook@
\newif\ifquotes@          \newif\ifbookquotes@      \newif\ifpaperquotes@
\newdimen\bysamerulewd@
\setboxz@h{\refsfont@\kern3em}
\bysamerulewd@\wdz@
\newdimen\refindentwd
\setboxz@h{\refsfont@ 00. }
\refindentwd\wdz@
\outer\def\ref{\begingroup \noindent\hangindent\refindentwd
 \firstref@true \def\nofrills{\def\refkern@{\kern3sp}}%
 \ref@}
\def\ref@{\book@false \bgroup\let\endrefitem@\egroup \ignorespaces}
\def\moreref{\endrefitem@\endref@\firstref@false\ref@}%
\def\transl{\endrefitem@\endref@\firstref@false
  \book@false
  \prepunct@
  \setboxz@h\bgroup \aftergroup\unhbox\aftergroup\z@
    \def\endrefitem@{\unskip\refkern@\egroup}\ignorespaces}%
\def\emptyifempty@{\dimen@\wd\currbox@
  \advance\dimen@-\wd\z@ \advance\dimen@-.1\p@
  \ifdim\dimen@<\z@ \setbox\currbox@\copy\voidb@x \fi}
\let\refkern@\relax
\def\endrefitem@{\unskip\refkern@\egroup
  \setboxz@h{\refkern@}\emptyifempty@}\ignorespaces
\def\refdef@#1#2#3{\edef\next@{\noexpand\endrefitem@
  \let\noexpand\currbox@\csname\expandafter\eat@\string#1box@\endcsname
    \noexpand\setbox\noexpand\currbox@\hbox\bgroup}%
  \toks@\expandafter{\next@}%
  \if\notempty{#2#3}\toks@\expandafter{\the\toks@
  \def\endrefitem@{\unskip#3\refkern@\egroup
  \setboxz@h{#2#3\refkern@}\emptyifempty@}#2}\fi
  \toks@\expandafter{\the\toks@\ignorespaces}%
  \edef#1{\the\toks@}}
\refdef@\no{}{. }
\refdef@\key{[\m@th}{] }
\refdef@\by{}{}
\def\bysame{\by\hbox to\bysamerulewd@{\hrulefill}\thinspace
   \kern0sp}
\def\manyby{\message{\string\manyby is no longer necessary; \string\by
  can be used instead, starting with version 2.0 of \styname.STY}\by}
\refdef@\paper{\ifpaperquotes@``\fi\it}{}
\refdef@\paperinfo{}{}
\def\jour{\endrefitem@\let\currbox@\jourbox@
  \setbox\currbox@\hbox\bgroup
  \def\endrefitem@{\unskip\refkern@\egroup
    \setboxz@h{\refkern@}\emptyifempty@
    \ifvoid\jourbox@\else\prevjour@true\fi}%
\ignorespaces}
\refdef@\vol{\ifbook@\else\bf\fi}{}
\refdef@\issue{no. }{}
\refdef@\yr{}{}
\refdef@\pages{}{}
\def\page{\endrefitem@\def\pp@{\def\pp@{pp.~}p.~}\let\currbox@\pagesbox@
  \setbox\currbox@\hbox\bgroup\ignorespaces}
\def\pp@{pp.~}
\def\book{\endrefitem@ \let\currbox@\bookbox@
 \setbox\currbox@\hbox\bgroup\def\endrefitem@{\unskip\refkern@\egroup
  \setboxz@h{\ifbookquotes@``\fi}\emptyifempty@
  \ifvoid\bookbox@\else\book@true\fi}%
  \ifbookquotes@``\fi\it\ignorespaces}
\def\inbook{\endrefitem@
  \let\currbox@\bookbox@\setbox\currbox@\hbox\bgroup
  \def\endrefitem@{\unskip\refkern@\egroup
  \setboxz@h{\ifbookquotes@``\fi}\emptyifempty@
  \ifvoid\bookbox@\else\book@true\previnbook@true\fi}%
  \ifbookquotes@``\fi\ignorespaces}
\refdef@\eds{(}{, eds.)}
\def\ed{\endrefitem@\let\currbox@\edsbox@
 \setbox\currbox@\hbox\bgroup
 \def\endrefitem@{\unskip, ed.)\refkern@\egroup
  \setboxz@h{(, ed.)}\emptyifempty@}(\ignorespaces}
\refdef@\bookinfo{}{}
\refdef@\publ{}{}
\refdef@\publaddr{}{}
\refdef@\finalinfo{}{}
\refdef@\lang{(}{)}

\let\refdef@\relax 
\def\ppunbox@#1{\ifvoid#1\else\prepunct@\unhbox#1\fi}
\def\nocomma@#1{\ifvoid#1\else\changepunct@3\prepunct@\unhbox#1\fi}
\def\changepunct@#1{\ifnum\lastkern<3 \unkern\kern#1sp\fi}
\def\prepunct@{\count@\lastkern\unkern
  \ifnum\lastpenalty=0
    \let\penalty@\relax
  \else
    \edef\penalty@{\penalty\the\lastpenalty\relax}%
  \fi
  \unpenalty
  \let\refspace@\ \ifcase\count@,
\or;\or.\or 
  \or\let\refspace@\relax
  \else,\fi
  \ifquotes@''\quotes@false\fi \penalty@ \refspace@
}
\def\transferpenalty@#1{\dimen@\lastkern\unkern
  \ifnum\lastpenalty=0\unpenalty\let\penalty@\relax
  \else\edef\penalty@{\penalty\the\lastpenalty\relax}\unpenalty\fi
  #1\penalty@\kern\dimen@}
\def\endref{\endrefitem@\lastref@true\endref@
  \par\endgroup \prevjour@false \previnbook@false }
\def\endref@{%
\iffirstref@
  \ifvoid\nobox@\ifvoid\keybox@\indent\fi
  \else\hbox to\refindentwd{\hss\unhbox\nobox@}\fi
  \ifvoid\keybox@
  \else\ifdim\wd\keybox@>\refindentwd
         \box\keybox@
       \else\hbox to\refindentwd{\unhbox\keybox@\hfil}\fi\fi
  \kern4sp\ppunbox@\bybox@
\fi 
  \ifvoid\paperbox@
  \else\prepunct@\unhbox\paperbox@
    \ifpaperquotes@\quotes@true\fi\fi
  \ppunbox@\paperinfobox@
  \ifvoid\jourbox@
    \ifprevjour@ \nocomma@\volbox@
      \nocomma@\issuebox@
      \ifvoid\yrbox@\else\changepunct@3\prepunct@(\unhbox\yrbox@
        \transferpenalty@)\fi
      \ppunbox@\pagesbox@
    \fi 
  \else \prepunct@\unhbox\jourbox@
    \nocomma@\volbox@
    \nocomma@\issuebox@
    \ifvoid\yrbox@\else\changepunct@3\prepunct@(\unhbox\yrbox@
      \transferpenalty@)\fi
    \ppunbox@\pagesbox@
  \fi 
  \ifbook@\prepunct@\unhbox\bookbox@ \ifbookquotes@\quotes@true\fi \fi
  \nocomma@\edsbox@
  \ppunbox@\bookinfobox@
  \ifbook@\ifvoid\volbox@\else\prepunct@ vol.~\unhbox\volbox@
  \fi\fi
  \ppunbox@\publbox@ \ppunbox@\publaddrbox@
  \ifbook@ \ppunbox@\yrbox@
    \ifvoid\pagesbox@
    \else\prepunct@\pp@\unhbox\pagesbox@\fi
  \else
    \ifprevinbook@ \ppunbox@\yrbox@
      \ifvoid\pagesbox@\else\prepunct@\pp@\unhbox\pagesbox@\fi
    \fi \fi
  \ppunbox@\finalinfobox@
  \iflastref@
    \ifvoid\langbox@.\ifquotes@''\fi
    \else\changepunct@2\prepunct@\unhbox\langbox@\fi
  \else
    \ifvoid\langbox@\changepunct@1%
    \else\changepunct@3\prepunct@\unhbox\langbox@
      \changepunct@1\fi
  \fi
}
\outer\def\enddocument{%
 \runaway@{proclaim}%
\ifmonograph@ 
\else
 \nobreak
 \thetranslator@
 \count@\z@ \loop\ifnum\count@<\addresscount@\advance\count@\@ne
 \csname address\number\count@\endcsname
 \csname email\number\count@\endcsname
 \repeat
\fi
 \vfill\supereject\end}

\def\headfont@{\headfonts}
\def\proclaimheadfont@{\bf}
\def\specialheadfont@{\bf}
\def\subheadfont@{\bf}
\def\demoheadfont@{\smc}

\newif\ifThisToToc \ThisToTocfalse
\newif\iftocloaded \tocloadedfalse

\def\C@L{\noexpand\Cal}\def\B@B{\noexpand\Bbb}\def\fR@K{\noexpand\frak}
\def\S@{\noexpand\S}\def\P@P{\noexpand\"}
\def\xpar{\\}

\def\writetoc#1{\iftocloaded\ifThisToToc\begingroup\def\totoc{}
  \def\Cal{\noexpand\C@L}\def\Bbb{\noexpand\B@B}
  \def\frak{\noexpand\fR@K}\def\goth{\frak}\def\S{\noexpand\S@}
  \def\"{\noexpand\P@P}
  \def\xpar{\par\penalty100000 }\def\idx##1{##1}\def\\{\xpar}
  \edef\next@{\write\toc{\noindent#1\leaderfill\noexpand\folio\par}}%
  \next@\endgroup\global\ThisToTocfalse\fi\fi}
\def\leaderfill{\leaders\hbox to 1em{\hss.\hss}\hfill}

\newif\ifindexloaded \indexloadedfalse
\def\idx#1{\ifindexloaded\begingroup\def\ign{}\def\it{}\def\/{}%
 \def\smc{}\def\bf{}\def\tt{}%
 \def\Cal{\noexpand\C@L}\def\Bbb{\noexpand\B@B}%
 \def\frak{\noexpand\fR@K}\def\goth{\frak}\def\S{\noexpand\S@}%
  \def\"{\noexpand\P@P}%
 {\edef\next@{\write\index{#1, \noexpand\folio}}\next@}%
 \endgroup\fi{#1}}
\def\ign#1{}

\def\totoc{\global\ThisToToctrue}

\outer\def\head#1\endhead{\par\penaltyandskip@{-200}\aboveheadskip
 {\headfont@\raggedcenter@\interlinepenalty\@M
 \ignorespaces#1\endgraf}\nobreak
 \vskip\belowheadskip
 \headmark{#1}\writetoc{#1}}

\outer\def\chaphead#1\endchaphead{\par\penaltyandskip@{-200}\aboveheadskip
 {\chapheadfonts\raggedcenter@\interlinepenalty\@M
 \ignorespaces#1\endgraf}\nobreak
 \vskip3\belowheadskip
 \headmark{#1}\writetoc{#1}}

\def\folio{{\foliofont@\ifnum\pageno<\z@ \romannumeral-\pageno
 \else\number\pageno \fi}}
\newtoks\leftheadtoks
\newtoks\rightheadtoks

\def\leftheadtext{\nofrills@{\relax}\lht@
  \DNii@##1{\leftheadtoks\expandafter{\lht@{##1}}%
    \mark{\the\leftheadtoks\noexpand\else\the\rightheadtoks}
    \ifsyntax@\setboxz@h{\def\\{\unskip\space\ignorespaces}%
        \headlinefont@##1}\fi}%
  \FN@\next@}
\def\rightheadtext{\nofrills@{\relax}\rht@
  \DNii@##1{\rightheadtoks\expandafter{\rht@{##1}}%
    \mark{\the\leftheadtoks\noexpand\else\the\rightheadtoks}%
    \ifsyntax@\setboxz@h{\def\\{\unskip\space\ignorespaces}%
        \headlinefont@##1}\fi}%
  \FN@\next@}
\def\NoRunningHeads{\global\runheads@false\global\let\headmark\eat@}

\newif\iffirstpage@     \firstpage@true
\newif\ifrunheads@      \runheads@true

\newdimen\fullhsize \fullhsize=\hsize
\newdimen\fullvsize \fullvsize=\vsize
\def\fullline{\hbox to\fullhsize}

\def\pagenumbers{\gdef\folio{\folio@}}

\let\norunningheads\NoRunningHeads
\def\userunningheads{\global\runheads@true}
\norunningheads

\headline={\def\chapter#1{\chapterno@. }%
  \def\\{\unskip\space\ignorespaces}\ifrunheads@\headlinefont@
    \ifodd\pageno\rightheadline \else\leftheadline\fi
   \else\hfil\fi\ifNoRunHeadline\global\NoRunHeadlinefalse\fi}
\let\folio@\folio
\def\foliofont@{\foliofont}
\def\foliofont{\eightrm}
\def\headlinefont@{\headlinefont}
\def\headlinefont{\eightpoint\smc}
\def\leftheadline{\rlap{\folio}\hfill
   \ifNoRunHeadline\else\iftrue\topmark\fi\fi \hfill}
\def\rightheadline{\hfill\ifNoRunHeadline
   \else \expandafter\fi
  \hfill \llap{\folio}}
\footline={{\eightpoint\bottremark}%
   \ifrunheads@\else\hfil{\let\foliofont\tenrm\folio}\fi\hfil}
\def\bottremark{}
 
\newif\ifNoRunHeadline      
\def\norunninghead{\global\NoRunHeadlinetrue}
\norunninghead

\output={\output@}
%
\newif\ifoffset\offsetfalse
\output={\output@}
\def\output@{%
 \ifoffset 
  \ifodd\count0\advance\hoffset by0.5truecm
   \else\advance\hoffset by-0.5truecm\fi\fi
 \shipout\vbox{%
  \makeheadline \pagebody \makefootline }%
 \advancepageno \ifnum\outputpenalty>-\@MM\else\dosupereject\fi}

\def\indexoutput#1{%
  \ifoffset 
   \ifodd\count0\advance\hoffset by0.5truecm
    \else\advance\hoffset by-0.5truecm\fi\fi
  \shipout\vbox{\makeheadline
  \vbox to\fullvsize{\boxmaxdepth\maxdepth%
  \ifvoid\topins\else\unvbox\topins\fi%
  #1 %
  \ifvoid\footins\else 
    \vskip\skip\footins
    \footnoterule
    \unvbox\footins\fi
  \ifr@ggedbottom \kern-\dimen@ \vfil \fi}%
  \baselineskip2pc
  \makefootline}%
 \global\advance\pageno\@ne
 \ifnum\outputpenalty>-\@MM\else\dosupereject\fi}
 
 \newbox\partialpage \newdimen\halfsize \halfsize=0.5\fullhsize
 \advance\halfsize by-0.5em

 \def\begindoublecolumns{\output={\indexoutput{\unvbox255}}%
   \begingroup \def\line{\fullline}
   \output={\global\setbox\partialpage=\vbox{\unvbox255\bigskip}}\eject
   \output={\doublecolumnout}\hsize=\halfsize \vsize=2\fullvsize}
 \def\enddoublecolumns{\output={\balancecolumns}\eject
  \endgroup \pagegoal=\fullvsize%
  \output={\output@}}
\def\doublecolumnout{\splittopskip=\topskip \splitmaxdepth=\maxdepth
  \dimen@=\fullvsize \advance\dimen@ by-\ht\partialpage
  \setbox0=\vsplit255 to \dimen@ \setbox2=\vsplit255 to \dimen@
  \indexoutput{\pagesofar} \unvbox255 \penalty\outputpenalty}
\def\pagesofar{\unvbox\partialpage
  \wd0=\hsize \wd2=\hsize \hbox to\fullhsize{\box0\hfil\box2}}
\def\balancecolumns{\setbox0=\vbox{\unvbox255} \dimen@=\ht0
  \advance\dimen@ by\topskip \advance\dimen@ by-\baselineskip
  \divide\dimen@ by2 \splittopskip=\topskip
  {\vbadness=10000 \loop \global\setbox3=\copy0
    \global\setbox1=\vsplit3 to\dimen@
    \ifdim\ht3>\dimen@ \global\advance\dimen@ by1pt \repeat}
  \setbox0=\vbox to\dimen@{\unvbox1} \setbox2=\vbox to\dimen@{\unvbox3}
  \pagesofar}

\tenpoint
\catcode`\@=\active

\def\smallheadings{\let\chapheadfonts\tenpoint\let\headfonts\tenpoint}

\tenpoint
\catcode`\@=\active

\def\LL{\leavevmode\setbox0=\hbox{L}\hbox to\wd0{\hss\char'40L}}
\def\al{\alpha}

\def\ep{\varepsilon}

\def\si{\sigma}

\def\De{\Delta}

\def\Om{\Omega}


\def\today{\ifcase\month\or
 January\or February\or March\or April\or May\or June\or
 July\or August\or September\or October\or November\or December\fi
 \space\number\day, \number\year}

\def\({\left(}
\def\){\right)}
\def\[{\left[}
\def\]{\right]}

\def\3{\ss}
\catcode`\@=11
\def\dddot#1{\vbox{\ialign{##\crcr
      .\hskip-.5pt.\hskip-.5pt.\crcr\noalign{\kern1.5\p@\nointerlineskip}
      $\hfil\displaystyle{#1}\hfil$\crcr}}}

\newif\iftab@\tab@false
\newif\ifvtab@\vtab@false
\def\tab{\bgroup\tab@true\vtab@false\vst@bfalse\Strich@false%
   \def\\{\global\hline@@false%
     \ifhline@\global\hline@false\global\hline@@true\fi\cr}
   \edef\l@{\the\leftskip}\ialign\bgroup\hskip\l@##\hfil&&##\hfil\cr}
\def\endtab{\cr\egroup\egroup}
\def\vtab{\vtop\bgroup\vst@bfalse\vtab@true\tab@true\Strich@false%
   \bgroup\def\\{\cr}\ialign\bgroup&##\hfil\cr}
\def\endvtab{\cr\egroup\egroup\egroup}
\def\stab{\D@cke0.5pt\null 
 \bgroup\tab@true\vtab@false\vst@bfalse\Strich@true\Let@@\vspace@
 \normalbaselines\offinterlineskip
  \openup\spreadmlines@
 \edef\l@{\the\leftskip}\ialign
 \bgroup\hskip\l@##\hfil&&##\hfil\crcr}
\def\endstab{\crcr\egroup
 \egroup}
\newif\ifvst@b\vst@bfalse
\def\vstab{\D@cke0.5pt\null
 \vtop\bgroup\tab@true\vtab@false\vst@btrue\Strich@true\bgroup\Let@@\vspace@
 \normalbaselines\offinterlineskip
  \openup\spreadmlines@\bgroup}
\def\endvstab{\crcr\egroup\egroup
 \egroup\tab@false\Strich@false}

\newdimen\htstrut@
\htstrut@8.5\p@
\newdimen\htStrut@
\htStrut@12\p@
\newdimen\dpstrut@
\dpstrut@3.5\p@
\newdimen\dpStrut@
\dpStrut@3.5\p@
\def\openup{\afterassignment\@penup\dimen@=}
\def\@penup{\advance\lineskip\dimen@
  \advance\baselineskip\dimen@
  \advance\lineskiplimit\dimen@
  \divide\dimen@ by2
  \advance\htstrut@\dimen@
  \advance\htStrut@\dimen@
  \advance\dpstrut@\dimen@
  \advance\dpStrut@\dimen@}
\def\Let@@{\relax%
    \def\\{\global\hline@@false%
     \ifhline@\global\hline@false\global\hline@@true\fi\cr}%
    \iffalse}\fi}
\def\matrix{\null\,\vcenter\bgroup
 \tab@false\vtab@false\vst@bfalse\Strich@false\Let@@\vspace@
 \normalbaselines\openup\spreadmlines@\ialign
 \bgroup\hfil$\m@th##$\hfil&&\quad\hfil$\m@th##$\hfil\crcr
 \Mathstrut@\crcr\noalign{\kern-\baselineskip}}
\def\endmatrix{\crcr\Mathstrut@\crcr\noalign{\kern-\baselineskip}\egroup
 \egroup\,}
\def\smatrix{\D@cke0.5pt\null\,
 \vcenter\bgroup\tab@false\vtab@false\vst@bfalse\Strich@true\Let@@\vspace@
 \normalbaselines\offinterlineskip
  \openup\spreadmlines@\ialign
 \bgroup\hfil$\m@th##$\hfil&&\quad\hfil$\m@th##$\hfil\crcr}
\def\endsmatrix{\crcr\egroup
 \egroup\,\Strich@false}
\newdimen\D@cke
\def\Dicke#1{\global\D@cke#1}
\newtoks\tabs@\tabs@{&}
\newif\ifStrich@\Strich@false
\newif\iff@rst

\def\Stricherr@{\iftab@\ifvtab@\errmessage{\noexpand\s not allowed
     here. Use \noexpand\vstab!}%
  \else\errmessage{\noexpand\s not allowed here. Use \noexpand\stab!}%
  \fi\else\errmessage{\noexpand\s not allowed
     here. Use \noexpand\smatrix!}\fi}
\def\format{\ifvst@b\else\crcr\fi\egroup\iffalse{\fi\ifnum`}=0 \fi\format@}
\def\format@#1\\{\def\preamble@{#1}%
 \def\Str@chfehlt##1{\ifx##1\s\Stricherr@\fi\ifx##1\\\let\Next\relax%
   \else\let\Next\Str@chfehlt\fi\Next}%
 \def\c{\hfil\noexpand\ifhline@@\hbox{\vrule height\htStrut@%
   depth\dpstrut@ width\z@}\noexpand\fi%
   \ifStrich@\hbox{\vrule height\htstrut@ depth\dpstrut@ width\z@}%
   \fi\iftab@\else$\m@th\fi\the\hashtoks@\iftab@\else$\fi\hfil}%
 \def\r{\hfil\noexpand\ifhline@@\hbox{\vrule height\htStrut@%
   depth\dpstrut@ width\z@}\noexpand\fi%
   \ifStrich@\hbox{\vrule height\htstrut@ depth\dpstrut@ width\z@}%
   \fi\iftab@\else$\m@th\fi\the\hashtoks@\iftab@\else$\fi}%
 \def\l{\noexpand\ifhline@@\hbox{\vrule height\htStrut@%
   depth\dpstrut@ width\z@}\noexpand\fi%
   \ifStrich@\hbox{\vrule height\htstrut@ depth\dpstrut@ width\z@}%
   \fi\iftab@\else$\m@th\fi\the\hashtoks@\iftab@\else$\fi\hfil}%
 \def\s{\ifStrich@\ \the\tabs@\vrule width\D@cke\the\hashtoks@%
          \fi\the\tabs@\ }%
 \def\sa{\ifStrich@\vrule width\D@cke\the\hashtoks@%
            \the\tabs@\ %
            \fi}%
 \def\se{\ifStrich@\ \the\tabs@\vrule width\D@cke\the\hashtoks@\fi}%
 \def\cd{\hfil\noexpand\ifhline@@\hbox{\vrule height\htStrut@%
   depth\dpstrut@ width\z@}\noexpand\fi%
   \ifStrich@\hbox{\vrule height\htstrut@ depth\dpstrut@ width\z@}%
   \fi$\dsize\m@th\the\hashtoks@$\hfil}%
 \def\rd{\hfil\noexpand\ifhline@@\hbox{\vrule height\htStrut@%
   depth\dpstrut@ width\z@}\noexpand\fi%
   \ifStrich@\hbox{\vrule height\htstrut@ depth\dpstrut@ width\z@}%
   \fi$\dsize\m@th\the\hashtoks@$}%
 \def\ld{\noexpand\ifhline@@\hbox{\vrule height\htStrut@%
   depth\dpstrut@ width\z@}\noexpand\fi%
   \ifStrich@\hbox{\vrule height\htstrut@ depth\dpstrut@ width\z@}%
   \fi$\dsize\m@th\the\hashtoks@$\hfil}%
 \ifStrich@\else\Str@chfehlt#1\\\fi%
 \setbox\z@\hbox{\xdef\Preamble@{\preamble@}}\ifnum`{=0 \fi\iffalse}\fi
 \ialign\bgroup\span\Preamble@\crcr}
\newif\ifhline@\hline@false
\newif\ifhline@@\hline@@false
\def\hlinefor#1{\multispan@{\strip@#1 }\leaders\hrule height\D@cke\hfill%
    \global\hline@true\ignorespaces}
\def\Item "#1"{\par\noindent\hangindent2\parindent%
  \hangafter1\setbox0\hbox{\rm#1\enspace}\ifdim\wd0>2\parindent%
  \box0\else\hbox to 2\parindent{\rm#1\hfil}\fi\ignorespaces}
\def\ITEM #1"#2"{\par\noindent\hangafter1\hangindent#1%
  \setbox0\hbox{\rm#2\enspace}\ifdim\wd0>#1%
  \box0\else\hbox to 0pt{\rm#2\hss}\hskip#1\fi\ignorespaces}
\def\item"#1"{\par\noindent\hang%
  \setbox0=\hbox{\rm#1\enspace}\ifdim\wd0>\the\parindent%
  \box0\else\hbox to \parindent{\rm#1\hfil}\enspace\fi\ignorespaces}
\let\plainitem@\item
\catcode`\@=13

\catcode`\@=11
\font\tenln    = line10
\font\tenlnw   = linew10

\newskip\Einheit \Einheit=0.5cm
\newcount\xcoord \newcount\ycoord
\newdimen\xdim \newdimen\ydim \newdimen\PfadD@cke \newdimen\Pfadd@cke

\newcount\@tempcnta
\newcount\@tempcntb

\newdimen\@tempdima
\newdimen\@tempdimb

\newdimen\@wholewidth
\newdimen\@halfwidth

\newcount\@xarg
\newcount\@yarg
\newcount\@yyarg
\newbox\@linechar
\newbox\@tempboxa
\newdimen\@linelen
\newdimen\@clnwd
\newdimen\@clnht

\newif\if@negarg

\def\@whilenoop#1{}
\def\@whiledim#1\do #2{\ifdim #1\relax#2\@iwhiledim{#1\relax#2}\fi}
\def\@iwhiledim#1{\ifdim #1\let\@nextwhile=\@iwhiledim
        \else\let\@nextwhile=\@whilenoop\fi\@nextwhile{#1}}

\def\@whileswnoop#1\fi{}
\def\@whilesw#1\fi#2{#1#2\@iwhilesw{#1#2}\fi\fi}
\def\@iwhilesw#1\fi{#1\let\@nextwhile=\@iwhilesw
         \else\let\@nextwhile=\@whileswnoop\fi\@nextwhile{#1}\fi}

\def\thinlines{\let\@linefnt\tenln \let\@circlefnt\tencirc
  \@wholewidth\fontdimen8\tenln \@halfwidth .5\@wholewidth}
\def\thicklines{\let\@linefnt\tenlnw \let\@circlefnt\tencircw
  \@wholewidth\fontdimen8\tenlnw \@halfwidth .5\@wholewidth}
\thinlines

\PfadD@cke1pt \Pfadd@cke0.5pt
\def\PfadDicke#1{\PfadD@cke#1 \divide\PfadD@cke by2 \Pfadd@cke\PfadD@cke \multiply\PfadD@cke by2}
\long\def\LOOP#1\REPEAT{\def\BODY{#1}\ITERATE}
\def\ITERATE{\BODY \let\next\ITERATE \else\let\next\relax\fi \next}
\let\REPEAT=\fi
\def\Punkt{\hbox{\raise-2pt\hbox to0pt{\hss$\ssize\bullet$\hss}}}
\def\DuennPunkt(#1,#2){\unskip
  \raise#2 \Einheit\hbox to0pt{\hskip#1 \Einheit
          \raise-2.5pt\hbox to0pt{\hss$\bullet$\hss}\hss}}
\def\NormalPunkt(#1,#2){\unskip
  \raise#2 \Einheit\hbox to0pt{\hskip#1 \Einheit
          \raise-3pt\hbox to0pt{\hss\twelvepoint$\bullet$\hss}\hss}}
\def\DickPunkt(#1,#2){\unskip
  \raise#2 \Einheit\hbox to0pt{\hskip#1 \Einheit
          \raise-4pt\hbox to0pt{\hss\fourteenpoint$\bullet$\hss}\hss}}
\def\Kreis(#1,#2){\unskip
  \raise#2 \Einheit\hbox to0pt{\hskip#1 \Einheit
          \raise-4pt\hbox to0pt{\hss\fourteenpoint$\circ$\hss}\hss}}

\def\Line@(#1,#2)#3{\@xarg #1\relax \@yarg #2\relax
\@linelen=#3\Einheit
\ifnum\@xarg =0 \@vline
  \else \ifnum\@yarg =0 \@hline \else \@sline\fi
\fi}

\def\@sline{\ifnum\@xarg< 0 \@negargtrue \@xarg -\@xarg \@yyarg -\@yarg
  \else \@negargfalse \@yyarg \@yarg \fi
\ifnum \@yyarg >0 \@tempcnta\@yyarg \else \@tempcnta -\@yyarg \fi
\ifnum\@tempcnta>6 \@badlinearg\@tempcnta0 \fi
\ifnum\@xarg>6 \@badlinearg\@xarg 1 \fi
\setbox\@linechar\hbox{\@linefnt\@getlinechar(\@xarg,\@yyarg)}%
\ifnum \@yarg >0 \let\@upordown\raise \@clnht\z@
   \else\let\@upordown\lower \@clnht \ht\@linechar\fi
\@clnwd=\wd\@linechar
\if@negarg \hskip -\wd\@linechar \def\@tempa{\hskip -2\wd\@linechar}\else
     \let\@tempa\relax \fi
\@whiledim \@clnwd <\@linelen \do
  {\@upordown\@clnht\copy\@linechar
   \@tempa
   \advance\@clnht \ht\@linechar
   \advance\@clnwd \wd\@linechar}%
\advance\@clnht -\ht\@linechar
\advance\@clnwd -\wd\@linechar
\@tempdima\@linelen\advance\@tempdima -\@clnwd
\@tempdimb\@tempdima\advance\@tempdimb -\wd\@linechar
\if@negarg \hskip -\@tempdimb \else \hskip \@tempdimb \fi
\multiply\@tempdima \@m
\@tempcnta \@tempdima \@tempdima \wd\@linechar \divide\@tempcnta \@tempdima
\@tempdima \ht\@linechar \multiply\@tempdima \@tempcnta
\divide\@tempdima \@m
\advance\@clnht \@tempdima
\ifdim \@linelen <\wd\@linechar
   \hskip \wd\@linechar
  \else\@upordown\@clnht\copy\@linechar\fi}

\def\@hline{\ifnum \@xarg <0 \hskip -\@linelen \fi
\vrule height\Pfadd@cke width \@linelen depth\Pfadd@cke
\ifnum \@xarg <0 \hskip -\@linelen \fi}

\def\@getlinechar(#1,#2){\@tempcnta#1\relax\multiply\@tempcnta 8
\advance\@tempcnta -9 \ifnum #2>0 \advance\@tempcnta #2\relax\else
\advance\@tempcnta -#2\relax\advance\@tempcnta 64 \fi
\char\@tempcnta}

\def\Vektor(#1,#2)#3(#4,#5){\unskip\leavevmode
  \xcoord#4\relax \ycoord#5\relax
      \raise\ycoord \Einheit\hbox to0pt{\hskip\xcoord \Einheit
         \Vector@(#1,#2){#3}\hss}}

\def\Vector@(#1,#2)#3{\@xarg #1\relax \@yarg #2\relax
\@tempcnta \ifnum\@xarg<0 -\@xarg\else\@xarg\fi
\ifnum\@tempcnta<5\relax
\@linelen=#3\Einheit
\ifnum\@xarg =0 \@vvector
  \else \ifnum\@yarg =0 \@hvector \else \@svector\fi
\fi
\else\@badlinearg\fi}

\def\@hvector{\@hline\hbox to 0pt{\@linefnt
\ifnum \@xarg <0 \@getlarrow(1,0)\hss\else
    \hss\@getrarrow(1,0)\fi}}

\def\@vvector{\ifnum \@yarg <0 \@downvector \else \@upvector \fi}

\def\@svector{\@sline
\@tempcnta\@yarg \ifnum\@tempcnta <0 \@tempcnta=-\@tempcnta\fi
\ifnum\@tempcnta <5
  \hskip -\wd\@linechar
  \@upordown\@clnht \hbox{\@linefnt  \if@negarg
  \@getlarrow(\@xarg,\@yyarg) \else \@getrarrow(\@xarg,\@yyarg) \fi}%
\else\@badlinearg\fi}

\def\@upline{\hbox to \z@{\hskip -.5\Pfadd@cke \vrule width \Pfadd@cke
   height \@linelen depth \z@\hss}}

\def\@downline{\hbox to \z@{\hskip -.5\Pfadd@cke \vrule width \Pfadd@cke
   height \z@ depth \@linelen \hss}}

\def\@upvector{\@upline\setbox\@tempboxa\hbox{\@linefnt\char'66}\raise
     \@linelen \hbox to\z@{\lower \ht\@tempboxa\box\@tempboxa\hss}}

\def\@downvector{\@downline\lower \@linelen
      \hbox to \z@{\@linefnt\char'77\hss}}

\def\@getlarrow(#1,#2){\ifnum #2 =\z@ \@tempcnta='33\else
\@tempcnta=#1\relax\multiply\@tempcnta \sixt@@n \advance\@tempcnta
-9 \@tempcntb=#2\relax\multiply\@tempcntb \tw@
\ifnum \@tempcntb >0 \advance\@tempcnta \@tempcntb\relax
\else\advance\@tempcnta -\@tempcntb\advance\@tempcnta 64
\fi\fi\char\@tempcnta}

\def\@getrarrow(#1,#2){\@tempcntb=#2\relax
\ifnum\@tempcntb < 0 \@tempcntb=-\@tempcntb\relax\fi
\ifcase \@tempcntb\relax \@tempcnta='55 \or
\ifnum #1<3 \@tempcnta=#1\relax\multiply\@tempcnta
24 \advance\@tempcnta -6 \else \ifnum #1=3 \@tempcnta=49
\else\@tempcnta=58 \fi\fi\or
\ifnum #1<3 \@tempcnta=#1\relax\multiply\@tempcnta
24 \advance\@tempcnta -3 \else \@tempcnta=51\fi\or
\@tempcnta=#1\relax\multiply\@tempcnta
\sixt@@n \advance\@tempcnta -\tw@ \else
\@tempcnta=#1\relax\multiply\@tempcnta
\sixt@@n \advance\@tempcnta 7 \fi\ifnum #2<0 \advance\@tempcnta 64 \fi
\char\@tempcnta}

\def\Diagonale(#1,#2)#3{\unskip\leavevmode
  \xcoord#1\relax \ycoord#2\relax
      \raise\ycoord \Einheit\hbox to0pt{\hskip\xcoord \Einheit
         \Line@(1,1){#3}\hss}}
\def\AntiDiagonale(#1,#2)#3{\unskip\leavevmode
  \xcoord#1\relax \ycoord#2\relax 
      \raise\ycoord \Einheit\hbox to0pt{\hskip\xcoord \Einheit
         \Line@(1,-1){#3}\hss}}
\def\Pfad(#1,#2),#3\endPfad{\unskip\leavevmode
  \xcoord#1 \ycoord#2 \thicklines\ZeichnePfad#3\endPfad\thinlines}
\def\ZeichnePfad#1{\ifx#1\endPfad\let\next\relax
  \else\let\next\ZeichnePfad
    \ifnum#1=1
      \raise\ycoord \Einheit\hbox to0pt{\hskip\xcoord \Einheit
         \vrule height\Pfadd@cke width1 \Einheit depth\Pfadd@cke\hss}%
      \advance\xcoord by 1
    \else\ifnum#1=2
      \raise\ycoord \Einheit\hbox to0pt{\hskip\xcoord \Einheit
        \hbox{\hskip-\PfadD@cke\vrule height1 \Einheit width\PfadD@cke depth0pt}\hss}%
      \advance\ycoord by 1
    \else\ifnum#1=3
      \raise\ycoord \Einheit\hbox to0pt{\hskip\xcoord \Einheit
         \Line@(1,1){1}\hss}
      \advance\xcoord by 1
      \advance\ycoord by 1
    \else\ifnum#1=4
      \raise\ycoord \Einheit\hbox to0pt{\hskip\xcoord \Einheit
         \Line@(1,-1){1}\hss}
      \advance\xcoord by 1
      \advance\ycoord by -1
    \else\ifnum#1=5
      \advance\xcoord by -1
      \raise\ycoord \Einheit\hbox to0pt{\hskip\xcoord \Einheit
         \vrule height\Pfadd@cke width1 \Einheit depth\Pfadd@cke\hss}%
    \else\ifnum#1=6
      \advance\ycoord by -1
      \raise\ycoord \Einheit\hbox to0pt{\hskip\xcoord \Einheit
        \hbox{\hskip-\PfadD@cke\vrule height1 \Einheit width\PfadD@cke depth0pt}\hss}%
    \else\ifnum#1=7
      \advance\xcoord by -1
      \advance\ycoord by -1
      \raise\ycoord \Einheit\hbox to0pt{\hskip\xcoord \Einheit
         \Line@(1,1){1}\hss}
    \else\ifnum#1=8
      \advance\xcoord by -1
      \advance\ycoord by +1
      \raise\ycoord \Einheit\hbox to0pt{\hskip\xcoord \Einheit
         \Line@(1,-1){1}\hss}
    \fi\fi\fi\fi
    \fi\fi\fi\fi
  \fi\next}
\def\hSSchritt{\leavevmode\raise-.4pt\hbox to0pt{\hss.\hss}\hskip.2\Einheit
  \raise-.4pt\hbox to0pt{\hss.\hss}\hskip.2\Einheit
  \raise-.4pt\hbox to0pt{\hss.\hss}\hskip.2\Einheit
  \raise-.4pt\hbox to0pt{\hss.\hss}\hskip.2\Einheit
  \raise-.4pt\hbox to0pt{\hss.\hss}\hskip.2\Einheit}
\def\vSSchritt{\vbox{\baselineskip.2\Einheit\lineskiplimit0pt
\hbox{.}\hbox{.}\hbox{.}\hbox{.}\hbox{.}}}
\def\DSSchritt{\leavevmode\raise-.4pt\hbox to0pt{%
  \hbox to0pt{\hss.\hss}\hskip.2\Einheit
  \raise.2\Einheit\hbox to0pt{\hss.\hss}\hskip.2\Einheit
  \raise.4\Einheit\hbox to0pt{\hss.\hss}\hskip.2\Einheit
  \raise.6\Einheit\hbox to0pt{\hss.\hss}\hskip.2\Einheit
  \raise.8\Einheit\hbox to0pt{\hss.\hss}\hss}}
\def\dSSchritt{\leavevmode\raise-.4pt\hbox to0pt{%
  \hbox to0pt{\hss.\hss}\hskip.2\Einheit
  \raise-.2\Einheit\hbox to0pt{\hss.\hss}\hskip.2\Einheit
  \raise-.4\Einheit\hbox to0pt{\hss.\hss}\hskip.2\Einheit
  \raise-.6\Einheit\hbox to0pt{\hss.\hss}\hskip.2\Einheit
  \raise-.8\Einheit\hbox to0pt{\hss.\hss}\hss}}
\def\SPfad(#1,#2),#3\endSPfad{\unskip\leavevmode
  \xcoord#1 \ycoord#2 \ZeichneSPfad#3\endSPfad}
\def\ZeichneSPfad#1{\ifx#1\endSPfad\let\next\relax
  \else\let\next\ZeichneSPfad
    \ifnum#1=1
      \raise\ycoord \Einheit\hbox to0pt{\hskip\xcoord \Einheit
         \hSSchritt\hss}%
      \advance\xcoord by 1
    \else\ifnum#1=2
      \raise\ycoord \Einheit\hbox to0pt{\hskip\xcoord \Einheit
        \hbox{\hskip-2pt \vSSchritt}\hss}%
      \advance\ycoord by 1
    \else\ifnum#1=3
      \raise\ycoord \Einheit\hbox to0pt{\hskip\xcoord \Einheit
         \DSSchritt\hss}
      \advance\xcoord by 1
      \advance\ycoord by 1
    \else\ifnum#1=4
      \raise\ycoord \Einheit\hbox to0pt{\hskip\xcoord \Einheit
         \dSSchritt\hss}
      \advance\xcoord by 1
      \advance\ycoord by -1
    \else\ifnum#1=5
      \advance\xcoord by -1
      \raise\ycoord \Einheit\hbox to0pt{\hskip\xcoord \Einheit
         \hSSchritt\hss}%
    \else\ifnum#1=6
      \advance\ycoord by -1
      \raise\ycoord \Einheit\hbox to0pt{\hskip\xcoord \Einheit
        \hbox{\hskip-2pt \vSSchritt}\hss}%
    \else\ifnum#1=7
      \advance\xcoord by -1
      \advance\ycoord by -1
      \raise\ycoord \Einheit\hbox to0pt{\hskip\xcoord \Einheit
         \DSSchritt\hss}
    \else\ifnum#1=8
      \advance\xcoord by -1
      \advance\ycoord by 1
      \raise\ycoord \Einheit\hbox to0pt{\hskip\xcoord \Einheit
         \dSSchritt\hss}
    \fi\fi\fi\fi
    \fi\fi\fi\fi
  \fi\next}
\def\Koordinatenachsen(#1,#2){\unskip
 \hbox to0pt{\hskip-.5pt\vrule height#2 \Einheit width.5pt depth1 \Einheit}%
 \hbox to0pt{\hskip-1 \Einheit \xcoord#1 \advance\xcoord by1
    \vrule height0.25pt width\xcoord \Einheit depth0.25pt\hss}}
\def\Koordinatenachsen(#1,#2)(#3,#4){\unskip
 \hbox to0pt{\hskip-.5pt \ycoord-#4 \advance\ycoord by1
    \vrule height#2 \Einheit width.5pt depth\ycoord \Einheit}%
 \hbox to0pt{\hskip-1 \Einheit \hskip#3\Einheit 
    \xcoord#1 \advance\xcoord by1 \advance\xcoord by-#3 
    \vrule height0.25pt width\xcoord \Einheit depth0.25pt\hss}}
\def\Gitter(#1,#2){\unskip \xcoord0 \ycoord0 \leavevmode
  \LOOP\ifnum\ycoord<#2
    \loop\ifnum\xcoord<#1
      \raise\ycoord \Einheit\hbox to0pt{\hskip\xcoord \Einheit\Punkt\hss}%
      \advance\xcoord by1
    \repeat
    \xcoord0
    \advance\ycoord by1
  \REPEAT}
\def\Gitter(#1,#2)(#3,#4){\unskip \xcoord#3 \ycoord#4 \leavevmode
  \LOOP\ifnum\ycoord<#2
    \loop\ifnum\xcoord<#1
      \raise\ycoord \Einheit\hbox to0pt{\hskip\xcoord \Einheit\Punkt\hss}%
      \advance\xcoord by1
    \repeat
    \xcoord#3
    \advance\ycoord by1
  \REPEAT}
\def\Label#1#2(#3,#4){\unskip \xdim#3 \Einheit \ydim#4 \Einheit
  \def\lo{\advance\xdim by-.5 \Einheit \advance\ydim by.5 \Einheit}%
  \def\llo{\advance\xdim by-.25cm \advance\ydim by.5 \Einheit}%
  \def\loo{\advance\xdim by-.5 \Einheit \advance\ydim by.25cm}%
  \def\o{\advance\ydim by.25cm}%
  \def\ro{\advance\xdim by.5 \Einheit \advance\ydim by.5 \Einheit}%
  \def\rro{\advance\xdim by.25cm \advance\ydim by.5 \Einheit}%
  \def\roo{\advance\xdim by.5 \Einheit \advance\ydim by.25cm}%
  \def\l{\advance\xdim by-.30cm}%
  \def\r{\advance\xdim by.30cm}%
  \def\lu{\advance\xdim by-.5 \Einheit \advance\ydim by-.6 \Einheit}%
  \def\llu{\advance\xdim by-.25cm \advance\ydim by-.6 \Einheit}%
  \def\luu{\advance\xdim by-.5 \Einheit \advance\ydim by-.30cm}%
  \def\u{\advance\ydim by-.30cm}%
  \def\ru{\advance\xdim by.5 \Einheit \advance\ydim by-.6 \Einheit}%
  \def\rru{\advance\xdim by.25cm \advance\ydim by-.6 \Einheit}%
  \def\ruu{\advance\xdim by.5 \Einheit \advance\ydim by-.30cm}%
  #1\raise\ydim\hbox to0pt{\hskip\xdim
     \vbox to0pt{\vss\hbox to0pt{\hss$#2$\hss}\vss}\hss}%
}
\catcode`\@=13

\hsize13cm
\vsize19cm
\newdimen\fullhsize
\newdimen\fullvsize
\newdimen\halfsize
\fullhsize13cm
\fullvsize19cm
\halfsize=0.5\fullhsize
\advance\halfsize by-0.5em

\magnification1200

\TagsOnRight

\def\AnBBAA{1}
\def\ArmDAA{2}
\def\AthaAI{3}
\def\AtReAA{4}
\def\AtTzAA{5}
\def\BesDAA{6}
\def\BjBrAB{7}
\def\BRWaAA{8}
\def\BRWaAB{9}
\def\CartAA{10}
\def\ChaFAA{11}
\def\EdelAA{12}
\def\FoReAA{13}
\def\FoReAB{14}
\def\GoJaAS{15}
\def\HumpAC{16}
\def\KratCB{17}
\def\KrewAC{18}
\def\PropAI{19}
\def\ReivAG{20}
\def\SimiAD{21}
\def\StanBZ{22}
\def\StanAP{23}
\def\StemAZ{24}

\def\AA{2.1}
\def\AB{2.2}
\def\AC{3.1}
\def\AD{4.1}
\def\AE{4.2}
\def\AF{4.3}
\def\AFa{4.4}
\def\AFb{4.5}
\def\AFc{4.6}
\def\AG{4.7}
\def\AGa{4.8}
\def\AGb{4.9}
\def\AGc{4.10}
\def\AGd{4.11}
\def\AR{4.12}
\def\SIa{4.13}
\def\SIb{4.14}
\def\SIc{4.15}
\def\SId{4.16}
\def\AHa{4.17}
\def\AH{4.18}
\def\AHb{4.19}
\def\AHc{4.20a}
\def\AHd{4.20b}
\def\AI{5.1}
\def\AJ{5.2}
\def\AKa{6.1}
\def\AKb{6.2}
\def\AKc{6.3}
\def\AKe{6.4}
\def\AKf{6.5}
\def\AKg{6.6}
\def\AKi{6.7}
\def\AKj{6.8}
\def\AK{6.9}
\def\ALa{7.1}
\def\ALb{7.2}
\def\ALc{7.3}
\def\ALd{7.4}
\def\ALe{7.5}
\def\ALf{7.6}
\def\ALg{7.7}
\def\ALh{7.8}
\def\ALgg{7.9}
\def\ALhh{7.10}
\def\ALi{7.11}
\def\ALj{7.12}
\def\ALk{7.13}
\def\ALl{7.14}
\def\ALm{7.15}
\def\ALmm{7.16}
\def\ALmmm{7.17}
\def\ALn{7.18}
\def\ALo{7.19}
\def\ALp{7.20}
\def\ALq{7.21}
\def\ALpp{7.22}
\def\ALr{7.23}
\def\AL{7.24}
\def\AM{8.1}
\def\AMa{8.2}
\def\AN{8.3}
\def\AO{9.1}
\def\AOa{9.2}
\def\AOb{9.3}
\def\AOc{9.4}
\def\AOd{9.6}
\def\AOe{9.5}
\def\AT{A.1}
\def\AU{A.2}
\def\AV{A.3}
\def\AW{A.4}
\def\AX{A.5}
\def\AP{A.6}
\def\AS{A.7}

\def\PA{1}
\def\PB{2}
\def\PE{3}
\def\PF{4}
\def\LBa{5}
\def\PC{6}
\def\PD{7}
\def\TA{8}

\def\al{\alpha}

\def\rk{\operatorname{rk}}

\topmatter 
\title The $M$-triangle of generalised non-crossing partitions
for the types $E_7$ and $E_8$
\endtitle 
\author C.~Krattenthaler$^{\dagger}$
\endauthor 
\affil 
Fakult\"at f\"ur Mathematik, Universit\"at Wien,\\
Nordbergstra{\ss}e~15, A-1090 Vienna, Austria.\\
WWW: \tt http://www.mat.univie.ac.at/\~{}kratt
\endaffil
\address Fakult\"at f\"ur Mathematik, Universit\"at Wien,
Nordbergstra{\ss}e~15, A-1090 Vienna, Austria.
WWW: \tt http://www.mat.univie.ac.at/\~{}kratt
\endaddress
\dedicatory
D\'edi\'e \`a Xavier Viennot
\enddedicatory
\thanks $^\dagger$Research partially supported 
by the Austrian Science Foundation FWF, grant S9607-N13,
in the framework of the National Research Network
``Analytic Combinatorics and Probabilistic Number Theory,"
and by
the ``Algebraic Combinatorics" Programme during Spring 2005 
of the Institut Mittag--Leffler of the Royal Swedish Academy of
Sciences%
\endthanks
\subjclass Primary 05E15;
 Secondary 05A05 05A15 05A19 06A07\linebreak 20F55 
\endsubjclass
\keywords root systems, 
generalised non-crossing partitions, M\"obius function,
$M$-triangle,
generalised cluster complex, face numbers, $F$-triangle
\endkeywords
\abstract 
The $M$-triangle of a ranked locally finite
poset $P$ is the generating function $\sum
_{u,w\in P} ^{}\mu(u,w)\,x^{\rk u}y^{\rk w}$, where $\mu(.,.)$ is the
M\"obius function of $P$. We compute the $M$-triangle of
Armstrong's poset of $m$-divisible non-crossing partitions for the 
root systems of type $E_7$ and
$E_8$. For the other types except $D_n$ this had been
accomplished in the earlier paper
``The $F$-triangle of the generalised cluster complex"
[in: ``Topics in Discrete Mathematics,'' M.~Klazar, J.~Kratochvil, M.~Loebl,
J.~Matou\v sek, R.~Thomas and P.~Valtr, eds., Springer--Verlag, Berlin,
New York, 2006, pp.~93--126].
Altogether, this almost settles Armstrong's
$F=M$ Conjecture, predicting a surprising relation between
the $M$-triangle of the $m$-divisible partitions poset and the
$F$-triangle (a certain refined face count) of the generalised
cluster complex of Fomin and Reading,
the only gap remaining in type $D_n$. Moreover, we prove
a reciprocity result for this $M$-triangle, 
again with the possible exception of type $D_n$.
Our results are based on the calculation of
certain decomposition numbers for the
reflection groups of types $E_7$ and $E_8$, which carry in fact finer
information than the $M$-triangle does. The decomposition numbers for the
other exceptional reflection groups had been computed in the earlier
paper. As an aside,
we show that there is a closed form product formula for the type $A_n$
decomposition numbers, leaving the problem of computing the type $B_n$
and type $D_n$ decomposition numbers open.
\endabstract
\endtopmatter
\document

\subhead 1. Introduction\endsubhead
The lattice of non-crossing partitions of Kreweras \cite{\KrewAC} is
a now classical object of study in combinatorics with many fascinating
properties (see \cite{\SimiAD}). 
In \cite{\EdelAA}, Edelman generalised non-crossing partitions 
to $m$-divisible non-crossing
partitions and showed that they, too, have many beautiful properties.
Recently, 
Bessis \cite{\BesDAA} and Brady and Watt \cite{\BRWaAA} gave a 
uniform definition of non-crossing partitions associated to
root systems, which includes Kreweras non-crossing partitions as
``type $A_n$" non-crossing partitions, as well as Reiner's
\cite{\ReivAG} ``type $B_n$" non-crossing partitions.
(As it turned out, Reiner did not have the ``right" definition of
``type $D_n$" non-crossing partitions). 
The question whether Edelman's \cite{\EdelAA}
$m$-divisible non-crossing partitions also allow for an extension to
root systems was answered by Armstrong \cite{\ArmDAA}, who introduced
$m$-divisible non-crossing partitions for all root systems uniformly.
(He calls them also ``generalised" non-crossing partitions for root
systems.)

Extending an earlier conjecture of Chapoton \cite{\ChaFAA}
(which, in the meantime, has become a theorem due to
Athanasiadis \cite{\AthaAI}), Armstrong predicted a relation between
the ``$M$-triangle" of his $m$-divisible non-crossing
partitions poset and the ``$F$-triangle" of the generalised cluster
complex of Fomin and Reading \cite{\FoReAA}. (Roughly speaking, 
the $M$-triangle is a bi-variate generating function involving the
M\"obius function, while the $F$-triangle is a refined face count.
The reader is referred to Section~2 for the definition of the
$M$-triangle, and to Section~8 for the definition of the
$F$-triangle.) In the sequel, we shall refer to Armstrong's
conjecture as the ``$F=M$ Conjecture."

In \cite{\KratCB}, the $F$-triangle of the generalised cluster
complex was computed for all types, and the $M$-triangle of the
$m$-divisible non-crossing partitions poset was computed for all types
except for the types $D_n$, $E_7$, and $E_8$. Thus, the $F=M$
Conjecture could be verified for all types except the afore-mentioned
three. The purpose of this paper is to compute the $M$-triangle of
the $m$-divisible non-crossing partitions poset for the exceptional
types $E_7$ and $E_8$. Thus, the $F=M$ Conjecture remains open only
for $D_n$. It should be noted, however, that a conjectural expression
for the $M$-triangle of the $m$-divisible non-crossing partitions
poset appears in \cite{\KratCB, Sec.~11, Prop.~D}. 
Moreover, since, in the case of type $D_n$, Section~11 of \cite{\KratCB}
contains a proof for $m=1$ (independent from
the one in \cite{\AthaAI}), the results from \cite{\KratCB} together
with the results of the present paper provide a case-by-case proof
of the $m=1$ case of the $F=M$ Conjecture, that is, of
Chapoton's (ex-)conjecture.

In order to compute the $M$-triangle of the $m$-divisible partitions
poset for $E_7$ and $E_8$, we compute certain decomposition numbers
in the corresponding reflection groups (see Section~3 for the
definition). These decomposition numbers carry, in fact, finer
information than the $M$-triangle. They
are therefore of intrinsic interest.
For the other reflection groups of
exceptional type, they had been already computed in \cite{\KratCB}, 
but not for the types $A_n$,
$B_n$, and $D_n$. In Section~10, we point out that, in fact, a result of
Goulden and Jackson \cite{\GoJaAS, Theorem~3.2} on the minimal factorization
of a long cycle implies a closed form product formula for
the type $A_n$ decomposition numbers. 
We plan to address the problem of computing the type $B_n$ and
type $D_n$ decomposition numbers in a future publication.

Our paper is organised as follows.
The next section contains the definition and basics of non-crossing
partitions for root systems and the $m$-divisible non-crossing
partitions. There, we define also the main object of this paper,
the $M$-triangle of the $m$-divisible non-crossing partitions poset.
Section~3 recalls the formula from \cite{\KratCB} which
expresses the $M$-triangle in terms of the afore-mentioned
decomposition numbers and the characteristic polynomials of 
non-crossing partitions posets of lower rank. 
Section~4 explains our strategy
of how to compute these decomposition numbers for the types $E_7$ and
$E_8$. (This strategy is actually applicable for any specific
reflection group, except that the computational difficulties increase
significantly with the rank of the reflection group.)
The intermediate Section~5 recalls from \cite{\KratCB} how these
decomposition numbers can be used to compute the characteristic polynomial
of the non-crossing partitions poset corresponding to the reflection
group. The programme from Section~4 is then implemented in Section~6
to compute the $M$-triangle of the $m$-divisible non-crossing
partitions poset of type $E_7$ and in Section~7 to compute the
$M$-triangle of the $m$-divisible non-crossing
partitions poset of type $E_8$. The purpose of Section~8 is to
briefly explain the $F=M$ Conjecture, and why the results from
Sections~6 and 7 together with results from \cite{\KratCB} prove it
for the types $E_7$ and $E_8$. Section~9 presents a curious
observation which results from our explicit expressions in this paper
and in \cite{\KratCB} for the
$M$-triangle of the $m$-divisible non-crossing partitions poset:
a reciprocity relation which sets the $M$-triangle 
with parameter $m$ in relation with the $M$-triangle with parameter
$-m$. It would be interesting to find an intrinsic explanation of
this phenomenon. We conclude the paper by addressing
the type $A_n$ decomposition numbers. 
Theorem~9 in Section~10 states the afore-mentioned result by Goulden and
Jackson in the language of the present paper, 
and Theorem~10 gives the implied result on arbitrary type
$A_n$ decomposition numbers. An Appendix lists
the decomposition numbers for the types 
$A_1,A_2,A_3,A_4,A_5,A_6,A_7,D_4,D_5,D_6,D_7,E_6$, which are needed
in our calculations of Sections~6 and 7.

\subhead 2. Generalised non-crossing partitions\endsubhead
In this section we recall the definition of Armstrong's \cite{\ArmDAA}
$m$-divisible non-crossing partitions poset, and we define its
$M$-triangle. 

Let $\Phi$ be a finite root system of rank $n$. (We refer the reader to
\cite{\HumpAC} for all root system terminology.)
For an element $\al\in\Phi$, let $t_\al$
denote the reflection in the central hyperplane perpendicular to
$\al$. Let $W=W(\Phi)$ be the group generated by these reflections. By
definition, any element $w$ of $W$ can be represented as a product 
$w=t_1t_2\cdots t_\ell$, where the $t_i$'s are reflections. We call
the minimal number of reflections which is needed for such a
product representation the {\it absolute length\/} of $w$, and we
denote it by $\ell_T(w)$. We then define the {\it absolute order} on
$W$, denoted by $\le_T$, by 
$$u\le_T w\quad \text{if and only if}\quad
\ell_T(w)=\ell_T(u)+\ell_T(u^{-1}w).$$
It can be shown that this is equivalent to the statement that any
shortest product representation of $u$ by reflections
occurs as an initial segment in some shortest product representation
of $w$ by reflections. 

We can now define the {\it non-crossing partition lattice
$NC(\Phi)$}. Let $c$ be a {\it Coxeter element\/} in $W$, that is, the
product of all reflections corresponding to the simple roots. Then
$NC(\Phi)$ is defined to be
the restriction of the partial order $\le_T$ to the set
of all elements which are less than or equal to $c$ in absolute order. 
This definition makes sense because, regardless
of the chosen Coxeter element $c$, the resulting poset
is always the same up to isomorphism.
It can be shown that $NC(\Phi)$ is in fact a lattice (see
\cite{\BRWaAB} for a uniform proof), and moreover self-dual
(this is obvious from the definition). Clearly, the minimal element
in $NC(\Phi)$ is the identity element in $W$, which we denote by
$\ep$, and the maximal element
in $NC(\Phi)$ is the chosen Coxeter element $c$.
The term ``non-crossing partition lattice" is used because
$NC(A_n)$ is isomorphic to the lattice of non-crossing partitions
originally introduced by Kreweras \cite{\KrewAC}
(see also \cite{\FoReAB}), and because also
$NC(B_n)$ and $NC(D_n)$ can be realised as lattices of
non-crossing partitions (see \cite{\AtReAA, \ReivAG}).

The poset of {\it $m$-divisible non-crossing partitions} has as a
ground-set the following subset of $(NC(\Phi))^{m+1}$,
$$\multline
NC^m(\Phi)=\{(w_0;w_1,\dots,w_m):w_0w_1\cdots w_m=c\text{ and }\\
\ell_T(w_0)+\ell_T(w_1)+\dots+\ell_T(w_m)=\ell_T(c)\}.
\endmultline\tag\AA$$
The order relation is defined by
$$(u_0;u_1,\dots,u_m)\le(w_0;w_1,\dots,w_m)\quad \text{if and only
if}\quad u_i\ge_T w_i,\ 1\le i\le m.$$
We emphasize that, according to this definition, $u_0$ and $w_0$ need
not be related in any way. The poset $NC^m(\Phi)$ is graded by 
the rank function
$$\rk\big((w_0;w_1,\dots,w_m)\big)=\ell_T(w_0).$$
Thus, there is a unique maximal element, namely $(c;\ep,\dots,\ep)$,
where $\ep$ stands for the identity element in $W$, but, if $m>1$, there 
are several
different minimal elements. In particular, there is no global minimum
in $NC^m(\Phi)$ if $m>1$ and, hence, $NC^m(\Phi)$ is not a lattice
for $m>1$. (It is, however, a graded join-semilattice, see
\cite{\ArmDAA, Theorem~3.4.4}.) 

The central object in the present paper is the ``$M$-triangle" of 
$NC^m(\Phi)$, which is the polynomial defined by
$$
M^m_\Phi(x,y)=\sum _{u,w\in NC^m(\Phi)} ^{}\mu(u,w)\,x^{\rk u}y^{\rk w},$$
where $\mu(u,w)$ is the M\"obius function in $NC^m(\Phi)$.
It is called ``triangle" because the M\"obius function
$\mu(u,w)$ vanishes unless $u\le w$, and, thus, the only coefficients
in the polynomial which may be non-zero are the coefficients of $x^ky^l$
with $0\le k\le l\le n$.

An equivalent object is the {\it dual $M$-triangle}, which
is defined by
$$
(M^m_\Phi)^*(x,y)=
\sum _{u,w\in(NC^m(\Phi))^*}
^{}\mu^*(u,w)\,x^{\rk^*w}y^{\rk^*u},
$$
where $(NC^m(\Phi))^*$ denotes the poset {\rm dual} to $NC^m(\Phi)$
{\rm(}i.e., the poset which arises from $NC^m(\Phi)$ by reversing all
order relations{\rm)}, where $\mu^*$ denotes the M\"obius
function in\linebreak 
$(NC^m(\Phi))^*$, and where $\rk^*$ denotes the rank
function in $(NC^m(\Phi))^*$. It is equivalent since, obviously, we
have
$$
(M^m_\Phi)^*(x,y)=(xy)^nM^m_\Phi(1/x,1/y).
\tag\AB$$
The dual $M$-triangle of $NC^m(\Phi)$ (and, thus, its $M$-triangle as
well) was computed explicitly in
\cite{\KratCB} for all types, except for $D_n$ (a conjectural expression
appears, however, in \cite{\KratCB, Sec.~11, Prop.~D}), for $E_7$, and
for $E_8$. In Sections~6 and 7 below, we fill this gap for $E_7$ and
$E_8$. Thus, it is only the case of $D_n$ which remains open.

\subhead 3. How to compute the $M$-triangle of the generalised
non-crossing partitions for a specific root system\endsubhead
We follow the strategy outlined in \cite{\KratCB, Sec.~12}, which,
however, has to be complemented by additional ideas.
These additional ideas will be described in the next section.

Let us recall from \cite{\KratCB} that
the dual $M$-triangle (and, thus, the $M$-triangle as well) 
can be expressed in terms of 
certain {\it decomposition numbers} and {\it characteristic
polynomials} of non-crossing partitions posets, which we define now.
The decomposition number $N_\Phi(T_1,T_2,\dots,T_d)$ is 
the number of ``minimal" products $c_1c_2\cdots c_d$
less than or equal to the Coxeter
element $c$ in absolute order, 
``minimal" meaning that all the $c_i$'s are different from
$\ep$ and that
$\ell_T(c_1)+\ell_T(c_2)+\dots+\ell_T(c_d)=\ell_T(c_1c_2\cdots c_d)$, 
such that the
type of $c_i$ as a parabolic Coxeter element is $T_i$, $i=1,2,\dots,d$.
(Here, the term ``parabolic Coxeter element" means a Coxeter element
in some parabolic subgroup.
The reader must recall that it follows from \cite{\BesDAA, Lemma~1.4.3}
that any element $c_i$ is indeed a Coxeter element in a parabolic subgroup
of $W=W(\Phi)$. By definition, the type of $c_i$ is the type of this
parabolic subgroup.)
On the other hand, we denote by
$\chi^*_{NC(\Psi)}(y)$ the reciprocal polynomial
of the characteristic polynomial of 
$NC(\Psi)$, that is, using self-duality of $NC(\Psi)$,
$$\chi^*_{NC(\Psi)}(y)=\sum _{u\in NC(\Psi)}
\mu(\hat 0_{NC(\Psi)},u)\,y^{\rk \Psi-\rk u}=\sum _{u\in NC(\Psi)}
\mu(u,\hat 1_{NC(\Psi)})\,y^{\rk u},$$
where $\hat 0_{NC(\Psi)}$ stands for the minimal element 
and $\hat 1_{NC(\Psi)}$ stands for the maximal element in
$NC(\Psi)$. (The reader should recall from Section~2 that
$\hat 0_{NC(\Psi)}$ is the identity element in $W(\Psi)$
and that $\hat 1_{NC(\Psi)}$ is the chosen Coxeter element in
$W(\Psi)$.) Using this notation, a combination of
Eqs.~(12.3) and (12.4) from \cite{\KratCB} reads as follows.

\proclaim{Proposition \PA}
For any finite root system $\Phi$ of rank $n$, we have
$$\multline
(M^m_\Phi)^*(x,y)=
\sum _{d=0} ^{n}
\sum _{(T_1,\dots,T_d)} ^{}x^{\rk T_1+\dots+\rk T_d}\cdot
N_\Phi(T_1,T_2,\dots,T_d)\\
\cdot\chi^*_{NC(T_1)}(y)
\chi^*_{NC(T_2)}(y)\cdots \chi^*_{NC(T_d)}(y)\binom md,
\endmultline\tag\AC$$
where the inner sum is over all possible $d$-tuples
$(T_1,T_2,\dots,T_d)$ of types {\rm(}not necessarily irreducible
types{\rm)}.
In {\rm(\AC)}, and in the sequel, the notation $NC(T)$  
means $NC(\Psi)$, where $\Psi$ is a root system
of type $T$, and $\rk T$ denotes the rank of $\Psi$.
\endproclaim

So, what we have to do to apply Formula~(\AC) to compute the
(dual) $M$-triangle is, first, to determine all the 
decomposition numbers
$N_\Phi(T_1,T_2,\dots,T_d)$, and, taking advantage of the
multiplicativity 
$$\chi^*_{NC(\Psi_1*\Psi_2)}(y)
=\chi^*_{NC(\Psi_1)\times NC(\Psi_2)}(y)
=\chi^*_{NC(\Psi_1)}(y)\cdot\chi^*_{NC(\Psi_2)}(y)$$
of the characteristic polynomial (here, $\Psi_1*\Psi_2$ denotes the
direct sum of the root systems $\Psi_1$ and $\Psi_2$), second, 
one needs a list of the characteristic
polynomials $\chi^*_{NC(\Psi)}(y)$ for all 
{\it irreducible} root systems $\Psi$.
We describe in the following section how we compute the decomposition 
numbers, and in the subsequent section how to use this knowledge
to compute as well the characteristic polynomials that we need.

\subhead 4. How to compute the decomposition numbers\endsubhead
Our strategy to compute the decomposition numbers 
$N_\Phi(T_1,T_2,\dots,T_d)$ for a fixed root system $\Phi$ is
to find as many linear relations between them as possible, and
eventually solve this system of linear equations. Indeed,
the decomposition numbers have many relations
between themselves, some of which have already been stated and used
in \cite{\KratCB}. We recall these here in Proposition~\PB\ below.
Equation~(\AD) says that the order of the types $T_1$, $T_2$, \dots,
$T_d$ is not relevant. Equation~(\AE) reduces the computation to
the computation of the decomposition numbers $N_\Phi(T_1,T_2,\dots,T_d)$
with ``full" rank, that is, with $\rk T_1+\rk T_2+\dots+\rk
T_d=\rk\Phi$. 

In Equation~(\AF) we add another family of
relations.
They involve the decomposition numbers $N_T(T'_1,T'_2,\dots,T'_e)$
of ``smaller" root systems than $\Phi$, that is, decomposition
numbers with $\rk T<\rk\Phi$. (Similarly to Proposition~\PA, by abuse
of notation, $N_T(T'_1,T'_2,\dots,T'_e)$ stands for the decomposition
number $N_\Psi(T'_1,T'_2,\dots,T'_e)$, where $\Psi$ is a root system
of type $T$.)
If we had already computed these
beforehand, then the relations (\AF) are linear relations
between full rank decomposition numbers $N_\Phi(T_1,T_2,\dots,T_d)$.
Proposition~\PD\ allows one to reduce these beforehand computations
to irreducible root systems of smaller rank.

A further set of linear relations comes from an identity featuring
the zeta polynomials of (ordinary and generalized) non-crossing
partitions posets given in Proposition~\PE. Indeed, as we explain
in Remark~(2) after the proof of Proposition~\PE, all the zeta
polynomials appearing in (\AG) are explicitly known, so that,
by comparing coefficients of $m^iz^j$ on both sides of (\AG),
we obtain a set of linear relations between the decomposition numbers
$N_\Phi(T_1,T_2,\dots,T_d)$.

While the number of equations which result from Propositions~\PB\ and
\PE\ exceeds the number of variables
(that is, the number of decomposition numbers of full rank) by far, 
it turns out that they
are not sufficient to determine them uniquely.
To remedy this somewhat, we add Proposition~\PF\ which provides the values
for three special decomposition numbers,
and we add Proposition~\PC\ which, as we illustrate in the
Remark after the statement of the proposition, allows one to compute
all the (full rank) decomposition numbers $N_\Phi(T,A_1)$ with 
$\rk T=\rk \Phi-1$. 

As we shall see in Sections~6 and 7, the system of linear equations
resulting from Propositions~\PB, \PE, \PF\ and \PC\ still does not
yield a unique solution for the decomposition numbers for $\Phi=E_7$ and
$\Phi=E_8$. However, they allow one to come very close, so that, upon
adding some arithmetic considerations and, in type~$E_8$, a
{\sl Maple} calculation using Stembridge's {\tt coxeter} package
\cite{\StemAZ}, one eventually succeeds in finding 
all the decomposition numbers.

\proclaim{Proposition \PB}
Let $\Phi$ be a finite root system. Then,
for any permutation $\si$ of $\{1,2,\dots,d\}$, we have
$$
N_\Phi(T_{\si(1)},T_{\si(2)},\dots,T_{\si(d)})=N_\Phi(T_1,T_2,\dots,T_d).
\tag\AD$$

Furthermore,
$$
N_\Phi(T_1,T_2,\dots,T_d)=
\sum _{T} ^{}N_\Phi(T_1,T_2,\dots,T_d,T),
\tag\AE$$
where the sum is over all types $T$ of rank $\rk\Phi-\rk T_1-\rk T_2-
\dots -\rk T_d$.

If 
$\rk T_1+\rk T_2+\dots +\rk T_d+
\rk T'_1+\rk T'_2+\dots +\rk T'_e=\rk\Phi$,
then
$$
N_\Phi(T_1,T_2,\dots,T_d,T'_1,T'_2,\dots,T'_e)=\sum _{T} ^{}
N_T(T'_1,T'_2,\dots,T'_e)N_\Phi(T_1,T_2,\dots,T_d,T),
\tag\AF$$
where the sum is over all types $T$ of rank $\rk T'_1+\rk T'_2+
\dots +\rk T'_e$.

Finally, if one of $T_1,T_2,\dots,T_d$ is not the type
of a sub-diagram of the Dynkin diagram of $\Phi$, then
$N_\Phi(T_1,T_2,\dots,T_d)=0$.
\endproclaim

\example{Examples}
An example of a relation of the form (\AE) is (\ALi). 
To be precise, Equation~(\ALi) is the special case of (\AE)
where $\Phi=E_8$, $d=1$, and $T_1=D_4$.

An example of a relation of the form (\AF) is
(we make also use of (\AD))
$$\align
N_{E_7}(A_1*A_3,A_1^2,A_1)
&=N_{A_1^5}(A_1*A_3,A_1)N_{E_7}(A_1^5,A_1^2)\\
&\kern1cm+
N_{A_1^3*A_2}(A_1*A_3,A_1)N_{E_7}(A_1^3*A_2,A_1^2)\\
&\kern1cm+
N_{A_1*A_2^2}(A_1*A_3,A_1)N_{E_7}(A_1*A_2^2,A_1^2)\\
&\kern1cm+
N_{A_1^2*A_3}(A_1*A_3,A_1)N_{E_7}(A_1^2*A_3,A_1^2)\\
&\kern1cm+
N_{A_2*A_3}(A_1*A_3,A_1)N_{E_7}(A_2*A_3,A_1^2)\\
&\kern1cm+
N_{A_1*A_4}(A_1*A_3,A_1)N_{E_7}(A_1*A_4,A_1^2)\\
&\kern1cm+
N_{A_1*D_4}(A_1*A_3,A_1)N_{E_7}(A_1*D_4,A_1^2)\\
&\kern1cm+
N_{A_5}(A_1*A_3,A_1)N_{E_7}(A_5,A_1^2)\\
&\kern1cm+
N_{D_5}(A_1*A_3,A_1)N_{E_7}(D_5,A_1^2)\\
&=2N_{E_7}(A_1^2*A_3,A_1^2)+
3N_{E_7}(A_2*A_3,A_1^2)+
5N_{E_7}(A_1*A_4,A_1^2)\\
&\kern1cm+
9N_{E_7}(A_1*D_4,A_1^2)+
6N_{E_7}(A_5,A_1^2)+
4N_{E_7}(D_5,A_1^2).
\endalign$$
To be precise, the above equation is the special case of (\AF)
where $\Phi=E_7$, $d=1$, $e=2$, $T_1=A_1^2$, $T_1'=A_1*A_3$ and
$T_2'=A_1$. The decomposition numbers $N_\Psi(\dots)$ with
$\Psi=A_1^5,A_1^3*A_2,A_1*A_2^2,
A_1^2*A_3,A_2*A_3,A_1*A_4,A_1*D_4,A_5,D_5$ which are needed
here are
computed by using Proposition~\PD\ and the decomposition numbers for
the irreducible root systems $A_1,A_2,A_3,A_4,A_5,D_4,D_5$ tabulated
in the Appendix.
\endexample

\demo{Proof of Proposition \PB} 
Equation (\AD) is \cite{\KratCB, Eq.~(45)}, while
Equation~(\AE) is \cite{\KratCB, Eq.~(46)}.

In order to see (\AF), we recall that, by definition and by taking
into account the rank assumption, the number
$N_\Phi(T_1,T_2,\dots,T_d,T'_1,T'_2,\dots,T'_e)$ is the number
of decompositions 
$$c=c_1c_2\cdots c_dc'_1c'_2\cdots c'_e,\tag\AFa$$
with $c_i$ of type $T_i$ and $c'_i$ of type $T'_i$ for all $i$.
The decompositions (\AFa) can also be determined by first
decomposing
$$c=c_1c_2\cdots c_dc',\tag\AFb$$
with $c_i$ of type $T_i$ for all $i$,
and subsequently
$$c'=c'_1c'_2\cdots c'_e,\tag\AFc$$
with $c'_i$ of type $T'_i$ for all $i$.
If we fix the type of $c'$, $T$ say, then the number of decompositions
(\AFb) is $N_\Phi(T_1,T_2,\dots,T_d,T)$. For determining the number
of decompositions (\AFc), we recall from \cite{\BesDAA, Cor.~1.6.2
and Def.~1.6.3} that $c'$ is a Coxeter element in a parabolic
subgroup, denoted by $W_{c'}$ in \cite{\BesDAA}, and that the absolute
order $\le_T$ on $W=W(\Phi)$ 
restricted to $W_{c'}$ is identical with absolute order
on $W_{c'}$. In other words, the decompositions (\AFc) with
$c'_i\in W$ are identical with
the decompositions (\AFc) with $c'_i\in W_{c'}$. Hence, the number of
decompositions (\AFc) is the number 
$N_T(T'_1,T'_2,\dots,T'_e)$. Finally, in order to get the
{\it total\/} number $N_\Phi(T_1,T_2,\dots,T_d)$, we must sum the 
products
$$N_T(T'_1,T'_2,\dots,T'_e)N_\Phi(T_1,T_2,\dots,T_d,T)$$
over all possible types $T$, thus arriving at Equation~(\AF).

The last assertion follows again from the fact \cite{\BesDAA,
Lemma~1.4.3} that any element which is less than or equal to a Coxeter
element $c$ of $W$ is the Coxeter element in some parabolic subgroup of $W$.
Hence, its type must by a sub-type of $\Phi$, the type of 
$W$.\quad \quad \qed
\enddemo

More equations come from the following proposition, featuring the
zeta polynomial of the non-crossing partitions posets
(ordinary and generalized).
Recall that, given a poset $P$, its {\it zeta polynomial\/} $Z_P(z)$
is the number of multichains $x_1\le x_2\le \dots\le x_{z-1}$ in
$P$. (It can be shown that this is indeed a
polynomial in $z$. The reader should consult \cite{\StanAP, Sec.~3.11} for
more information on this topic.)

\proclaim{Proposition \PE}
For any finite root system $\Phi$ of rank $n$, we have
$$\multline
Z_{NC^m(\Phi)}(z)
=\sum _{d=0} ^{n}
\sum _{(T_1,\dots,T_d)} ^{}
N_\Phi(T_1,T_2,\dots,T_d)\\
\cdot Z_{NC(T_1)}(z-1)
Z_{NC(T_2)}(z-1)\cdots Z_{NC(T_d)}(z-1)\binom md,
\endmultline\tag\AG$$
where the inner sum is over all possible $d$-tuples
$(T_1,T_2,\dots,T_d)$ of types,
and where $N_\Phi(T_1,T_2,\dots,T_d)$ is as before.
\endproclaim
\demo{Proof} 
Let $\hat 1_{NC^m(\Phi)}$ denote the maximal element 
$(c;\ep,\dots,\ep)$ in $NC^m(\Phi)$. If, given a multichain
$x_1\le x_2\le \dots\le x_{z-1}$ in $NC^m(\Phi)$, 
we remove the minimal element
$w=x_1$, then a multichain remains which is by $1$ shorter.
Hence, by summing over all possible $w$'s, we obtain
$$Z_{NC^m(\Phi)}(z)=\sum _{w\in NC^m(\Phi)} ^{}
Z_{[w,\hat 1_{NC^m(\Phi)}]}(z-1).\tag\AGa$$
Now, by definition of $NC^m(\Phi)$, $w$ is of the form
$(w_0;w_1,\dots,w_m)$ with $w_0w_1\cdots w_m=c$ and
$\ell_T(w_0)+\ell_T(w_1)+\dots+\ell_T(w_m)=\ell_T(c)=n$. 
Moreover, as we already recalled in Section~3, it follows from
\cite{\BesDAA, Lemma~1.4.3} that any element $w_j$
is a Coxeter element in a parabolic subgroup of $W=W(\Phi)$
(or, in the language used earlier, a ``parabolic Coxeter element").
On the other hand, we have
$$[w,\hat 1_{NC^m(\Phi)}]\cong[\ep,w_{1}]\times[\ep,w_{2}]\times
\dots\times[\ep,w_{m}],\tag\AGb$$
where each interval $[\ep,w_{j}]$ is an interval in $NC(\Phi)$.
In fact, some of the $w_j$'s may be equal to the identity $\ep$.
If we denote by $w_{i_1}, w_{i_2},\dots,w_{i_d}$ those among the
$w_j$'s which are {\it not\/} equal to $\ep$, then (\AGb) reduces to
$$[w,\hat 1_{NC^m(\Phi)}]\cong[\ep,w_{i_1}]\times[\ep,w_{i_2}]\times
\dots\times[\ep,w_{i_d}].$$
More precisely, since each $w_{i_j}$ is a parabolic Coxeter element,
each interval $[\ep,w_{i_j}]$ is isomorphic to some non-crossing
partition lattice $NC(\Psi)$, where $\Psi$ is the
root system of this parabolic subgroup.
The zeta polynomial being multiplicative, this implies
$$Z_{[w,\hat 1_{NC^m(\Phi)}]}(z-1)=
Z_{[\ep,w_{i_1}]}(z-1)Z_{[\ep,w_{i_2}]}(z-1)
\cdots Z_{[\ep,w_{i_d}]}(z-1).$$
If we take into account that the number of possibilities to choose
the indices $\{i_1<i_2<\dots<i_d\}$ out of $\{1,2,\dots,m\}$ is equal
to $\binom md$ and combine this with the previous observations,
then (\AGa) turns into (\AG).\quad \quad
\qed
\enddemo

\remark{Remarks} 
(1) It is striking to note the similarities between (\AC) and (\AG).
It would be interesting to find an intrinsic explanation.
Even in lack of such an explanation, the relations (\AC) and (\AG)
underline the significance of the decomposition numbers
$N_\Phi(T_1,T_2,\dots,T_d)$.

(2) The zeta polynomial of the non-crossing partition lattice
$NC(\Phi)$, where $\Phi$ is a root system of rank $n$, 
has the elegant formula (see \cite{\ChaFAA, Prop.~9})
$$Z_{NC(\Phi)}(z)=\prod _{i=1} ^{n}\frac {(z-1)h+d_i} {d_i},\tag\AGc$$
where $h$ is the {\it Coxeter number} and the $d_i$'s are the {\it
degrees} of
$W=W(\Phi)$. (The reader should be warned that the convention chosen
in \cite{\ChaFAA} for the zeta polynomial is different from the one
here.)
This uniform formula was originally conjectured by
Chapoton on the basis of the already known formulae in types $A_n$,
$B_n$ and $D_n$, and of calculations he did for some exceptional
groups. The formula was finally confirmed for the large exceptional groups
by {\sl Mathematica} and {\sl MATLAB} calculations done by Reiner.

On the basis of (\AGc), Armstrong could also determine the zeta
polynomial for all $m$-divisible non-crossing partitions posets.
Again, there is a uniform formula, namely (see \cite{\ArmDAA,
Theorem~3.5.2}; we warn the reader that also in \cite{\ArmDAA} 
this different convention for the zeta polynomial is used)
$$Z_{NC^m(\Phi)}(z)=\prod _{i=1} ^{}\frac {(z-1)mh+d_i} {d_i}.\tag\AGd$$

If we use formulae (\AGc) and (\AGd) in (\AG), then,
by comparing coefficients of $m^iz^j$ on both sides of (\AG),
we obtain $(n+1)^2$ linear equations for the decomposition numbers
$N_\Phi(T_1,T_2,\dots,T_d)$. Although some of them turn out to be
trivial ($1=1$ or $0=0$), this is nevertheless a considerable number
of linear relations.
\endremark

\proclaim{Proposition \PF}
{\rm(1)} For any finite root system $\Phi$ of type $T$ we have $N_{\Phi}(T)=1$.

{\rm(2)} The decomposition number $N_{\Phi}(A_1)$
is equal to the total number of reflections in $W(\Phi)$ or,
equivalently, 
to the number of positive roots in $\Phi$. These numbers are
tabulated in \cite{\HumpAC, Table~2 on p.~80}.

{\rm(3)} The decomposition number $N_{\Phi}(A_1,A_1,\dots,A_1)$
(with $n=\rk\Phi$ occurrences of $A_1$) is equal to $n!\,h^n/\vert
W(\Phi)\vert$, where $h$ is the Coxeter number of $W(\Phi)$. 
\endproclaim
\demo{Proof} The claim (1) is trivial.

To see (2), we have to show that any reflection $t$ is less than or
equal to the Coxeter element $c$ in absolute order. Indeed, we have
$c=t(tc)$, with $\ell_T(t)=1$ and $\ell_T(tc)=\ell_T(c)\pm1$,
by general properties of the absolute length. Now, the maximal
possible absolute length of an element in $W(\Phi)$ is
$\rk\Phi=\ell_T(c)$
(cf\. \cite{\BesDAA, Lemma~1.2.1(ii)}). Hence, we have
$\ell_T(tc)=\ell_T(c)-1$, and, thus, $t\le_Tc$. 

Finally we prove claim (3).
The number of decompositions $c=t_1t_2\cdots t_n$ is equal to the
number of maximal chains $\ep<w_1<w_2<\dots<w_{n-1}<c$ in $NC(\Phi)$, 
because of the obvious bijection between decompositions and maximal chains 
defined by the identification $w_i=t_1t_2\dots t_i$,
$i=1,2,\dots,n-1$. It is a general fact (see 
\cite{\StanAP, Prop.~3.11.1})
that the number of maximal chains in a poset is equal to the
leading coefficient of its zeta polynomial multiplied by the
factorial of the rank of the poset. 
The claim then follows from the explicit
formula (\AGc) for the zeta polynomial of $NC(\Phi)$.\quad \quad \qed
\enddemo

Our next goal is to determine the decomposition numbers $N_\Phi(T,A_1)$
with $\rk T=\rk\Phi-1$. Proposition~\PC\ will provide the means to do
that. For the proof of the proposition, we need an
auxiliary lemma, Lemma~\LBa\ below. 
It is an extended version of Lemma~1.3.4 from
\cite{\BesDAA}. The extension makes the theoretical assertion
of \cite{\BesDAA, Lemma~1.3.4} concrete for each type. Actually, in
the present paper, we shall need this ``concretisation" only for the types
$D_6$ and $E_7$. However, as it may be useful in other situations,
we work it out here for all types.

\proclaim{Lemma \LBa}
Let $\Phi$ be an irreducible root systems of rank $n$. Furthermore, 
let 
$$S=\{s_1,s_2,\dots,s_n\}=\{s_1,\dots,s_r\}\cup\{s_{r+1},\dots,s_n\}$$ 
be a choice of simple reflections with the property that $s_i$ and $s_j$
commute for all $i,j$ with $1\le i,j\le r$, and for all $i,j$ with
$r+1\le i,j\le n$, $r$ being chosen appropriately {\rm(}cf\.
\cite{\HumpAC, Sec.~3.17}{\rm)}, and let
$c=s_1s_2\cdots s_n$ be the corresponding Coxeter element.
Then, for any reflection $t$, the cardinality of the orbit 
$\Om(t)=\{c^ktc^{-k}:k\in\Bbb Z\}$ of $t$ under conjugation by $c$ is 
\roster
\item $h$ if $\vert\Om(t)\cap S\vert=2$,
\item $h/2$ if $\vert\Om(t)\cap S\vert=1$,
\endroster
where $h$ is the Coxeter number of $W=W(\Phi)$ {\rm(}cf\.
\cite{\HumpAC, Table~2 on p.~80}{\rm)},
and there are no other possibilities.
Specifically, the cardinality of $\Om(t)$ is
$$
\cases h&\text {if $\Phi=A_n$ except if $n$ is odd and $tc$ is of
type $A_{(n-1)/2}^2$,}\\
&\text {if $\Phi=D_n$, $n$ is odd, and $tc$ is of type $A_{n-1}$,}\\
&\text {if $\Phi=E_6$ and $tc$ is of type $D_5$ or $A_1*A_4$,}\\
h/2&\text {otherwise.}\endcases$$
\endproclaim

\demo{Proof} 
The assertion on the cardinality of $\Om(t)$ 
as it depends on the
cardinality of the intersection $\Om(t)\cap S$ is \cite{\BesDAA,
Lemma~1.3.4}. The concrete cardinality assertion for the type $A_n$
can be worked out by using the well-known combinatorial realisation of 
$W(A_n)$ as the symmetric group on $n+1$ elements, while for the
types $B_n$ and $D_n$ this can be worked out by using the
combinatorial realisations of the corresponding reflection groups
as subgroups of the symmetric group on $2n$ elements (see e.g\.
\cite{\BjBrAB, Sections~8.1 and 8.2}.) Since this does not contain
any surprises, we leave the details to the reader.

\midinsert
\vbox{
$$
\Einheit1cm
\Pfad(0,0),112611\endPfad
\DickPunkt(0,0)
\DickPunkt(1,0)
\DickPunkt(2,0)
\DickPunkt(2,1)
\DickPunkt(3,0)
\DickPunkt(4,0)
\Label\u{1}(0,0)
\Label\o{2}(2,1)
\Label\u{3}(1,0)
\Label\u{4}(2,0)
\Label\u{5}(3,0)
\Label\u{6}(4,0)
\hskip4cm
$$
\centerline{\eightpoint Dynkin diagram of $E_6$}
\vskip7pt
\centerline{\eightpoint Figure 1}
}
\endinsert

For the type $E_6$ we argue as follows. 
Let $s_1,s_2,s_3,s_4,s_5,s_6$ denote the simple
reflections of $E_6$, with commutation relations coded by the Dynkin
diagram of $E_6$ as given by Figure~1.
(The labelling of the nodes is the one which
Stembridge's {\tt coxeter} package \cite{\StemAZ} chooses.)
To prepare the following arguments, we note that
if $T$ is the type of the sub-diagram of the Dynkin diagram of $E_6$
obtained by deleting the node corresponding to the simple reflection
$s_i$, say, then, because of 
$c=s_1s_2s_3s_4s_5s_6=t(s_1\cdots s_{i-1}s_{i+1}\cdots s_6)$
with $t=(s_1\cdots s_{i-1})s_i(s_1\cdots s_{i-1})^{-1}$,
there is (at least) one
orbit $\{c^ktc^{-k}:k\in\Bbb Z\}$ such that the type of $tc$
is $T$. The types which
occur as types of sub-diagrams of the Dynkin diagram of $E_6$
obtained by deleting one node from it are
$D_5$, $A_1*A_4$, $A_5$ and $A_1*A_2^2$. 

{\sl Maple} computations
using Stembridge's package yield
the following facts. There are four orbits of reflections
under conjugation by $c$,
two of size $12=h$, and two of size $6=h/2$. Hence, 
for each of the types $D_5$, $A_1*A_4$, $A_5$, and $A_1*A_2^2$,
there corresponds exactly one orbit in the sense described above,
that is, if $T$ is one of these four types then there is exactly one
orbit such that $tc$ is of type $T$ for each reflection $t$ in the orbit.

The orbit of $s_1$ has 12
elements, and it contains also $s_3,s_4,s_5,s_6$. From the fact that
it contains $s_1$, we can conclude that $tc$ is of type $D_5$ for all
reflections $t$ in this orbit. The second orbit with $12$ elements is
the one of the reflection
$s_2s_3s_4s_2s_3$. Since the order of $s_2s_3s_4s_2s_3c$ is 10, we 
must necessarily have that its type is $A_1*A_4$. (The other
remaining sub-diagrams of the Dynkin diagram of $E_6$ have types 
$A_5$ and $A_1*A_2^2$. Coxeter elements of these types have the order
$6$.) The orbits with $6$ elements are the ones of $s_2$ and
$s_4s_5s_6s_5s_4$, respectively. If $t$ denotes any reflection in one
of these two orbits, then the type of $tc$ is necessarily either $A_5$
or $A_1*A_2^2$, one type applying for {\it all\/} reflections in one
orbit, the other type applying for {\it all\/} reflections in the
other orbit. For our purpose it is not important which type is associated
to which orbit, since both possibilities are in agreement with our claim.

\smallskip
In the case that $\Phi=I_2(a)$, there is just one orbit.
For all other exceptional types, one can again do calculations using 
Stembridge's {\tt coxeter} package. The result is that
in type $H_3$ one obtains $3$ orbits of $5$ elements each,
in type $H_4$ one obtains $4$ orbits of $15$ elements each,
in type $F_4$ one obtains $4$ orbits of $6$ elements each,
in type $E_7$ one obtains $7$ orbits of $9$ elements each, and
in type $E_8$ one obtains $8$ orbits of $15$ elements each.
All this is
in accordance with our claim.\quad \quad \qed
\enddemo

\proclaim{Proposition \PC}
Let $\Phi$ be a finite irreducible root system, and let $T$
be a type with $\rk T=\rk \Phi-1$. Then, if $T$ is the type of a 
sub-diagram of the Dynkin diagram of $\Phi$
obtained by deleting one node from it, we have
$$\frac {N_\Phi(T,A_1)}
{N(T\subseteq\Phi)}=\frac {h} {2},
\tag\AR$$
where $N(T\subseteq\Phi)$ denotes the number of times $T$
arises as type of a sub-diagram of the Dynkin diagram of $\Phi$,
and where $h$ is the Coxeter number of $W=W(\Phi)$.
Otherwise, we have $N_\Phi(T,A_1)=0$.
\endproclaim

\demo{Proof} 
We start with the following observation.
By definition, the number $N_\Phi(T,A_1)$ counts the number of product
decompositions $c=wt$, where $w$ is a parabolic Coxeter element of 
type $T$ and $t$ is a reflection. Given $w$ and $t$,
we obtain another product decomposition by conjugation,
$$c=ccc^{-1}=cwtc^{-1}=(cwc^{-1})(ctc^{-1}),$$
where $ctc^{-1}$ is also a reflection, and where the type of $cwc^{-1}$
is still $T$. Hence, all the elements of the complete orbit
$\{(c^kwc^{-k},c^ktc^{-k}):k\in\Bbb Z\}$
of $(w,t)$ under conjugation by $c$ provide product decompositions of
$c$ in a parabolic Coxeter element of type $T$ and a reflection.
Since the reflection $t'$ already determines its ``companion"
$w'$ in the decomposition $c=w't'$ uniquely via $w'=ct'$, 
we may as well concentrate on the orbit
$\{c^ktc^{-k}:k\in\Bbb Z\}$
of $t$ under conjugation by $c$.

In what follows, we assume the setup of Lemma~\LBa, that is, 
we assume that
$$S=\{s_1,s_2,\dots,s_n\}=\{s_1,\dots,s_r\}\cup\{s_{r+1},\dots,s_n\}$$ 
is a choice of simple reflections with the property that $s_i$ and $s_j$
commute for all $i,j$ with $1\le i,j\le r$, and for all $i,j$ with
$r+1\le i,j\le n$, $r$ being chosen appropriately, and we let
$c=s_1s_2\cdots s_n$ be the corresponding Coxeter element.

Now, let $T$ be the type of the
sub-diagram of the Dynkin diagram of $\Phi$
obtained by deleting the node corresponding to $s_i$ from it.
If $1\le i\le r$, we have
$$\align 
c&=\big(s_1\cdots s_{i-1}s_{i+1}\cdots s_r(s_is_{r+1}s_i)\cdots
(s_is_ns_i)\big)\cdot s_i\\
&=\big((s_is_1s_i)\cdots (s_is_{i-1}s_i)(s_is_{i+1}s_i)\cdots
(s_is_ns_i)\big)\cdot s_i,
\tag\SIa
\endalign$$
and, if $r+1\le i\le n$, we have
$$c=\big(s_1\cdots s_{i-1}s_{i+1}\cdots s_n\big)\cdot s_i.$$
In both cases, the type of $cs_i$ is $T$, in the former case since 
the system
$$\{(s_is_1s_i),\dots ,(s_is_{i-1}s_i)(s_is_{i+1}s_i),\dots,
(s_is_ns_i)\}
\tag\SIb$$
is conjugate to $\{s_1,\dots,s_{i-1},s_{i+1},\dots,s_n\}$, the latter
being a system of simple reflections of type $T$.

Next, let $\Om_1,\Om_2,\dots,\Om_a$ and
$\bar\Om_1,\bar\Om_2,\dots,\bar\Om_b$ be the orbits of reflections $t$
under conjugation by $c$ such that $ct$ is of type 
$T$,\footnote{A case-by-case analysis
shows that there can be at most $2$ such orbits, and $2$ orbits only if $\Phi$
is of type $D_n$, $n$ is even, and $T=A_{n-1}$. However, this is of no
importance here.} 
the former those of cardinality $h$, the
latter those of cardinality $h/2$. Lemma~\LBa\ says that
there can be no others, while the argument involving (\SIa) and (\SIb)
implies that there is at least one such orbit, namely the orbit of
$s_i$. The last statement implies in particular that at least
one of $a$ and $b$ is non-zero. We have
$$N(T,A_1)=ah+b\frac {h} {2}.\tag\SIc$$
On the other hand, from Lemma~\LBa~(1),(2) and the argument involving
(\SIa) and (\SIb) we infer that
$$N(T\subseteq \Phi)=2a+b.\tag\SId$$
Dividing right-hand and left-hand sides, respectively, of (\SIc) and
(\SId) we obtain (\AR).

If $T$ is a type which is {\it not\/} the type of a
sub-diagram of the Dynkin diagram of $\Phi$ and $N(T,A_1)\ne0$, 
then we obtain a contradiction: let $t$ be a reflection such that 
$c=wt$, where $w$ is of type $T$. 
By Lemma~\LBa, the orbit $\Om(t)$ of
$t$ under conjugation by $c$ contains a simple reflection, $s_i$ say.
Now the argument involving (\SIa) and (\SIb) yields that $T$ is the
type of a sub-diagram of the Dynkin diagram of $\Phi$, which is
absurd.\quad \quad \qed
\enddemo

The final proposition in this section reduces the computation of
the decomposition numbers $N_{\Phi}(T_1,T_2,\dots,T_d)$ to
irreducible root systems $\Phi$. On the right-hand side of
(\AH) below, we use an extended definition of decomposition numbers,
where we allow some of the types $T_i$ in the argument to be empty,
which we denote by $T_i=\emptyset$:
we set $N_{\Phi}(\emptyset,T_2,\dots,T_d)=N_{\Phi}(T_2,\dots,T_d)$, 
with analogous conventions if one or more of the other $T_i$'s 
should be empty.

\proclaim{Proposition \PD}
If 
$$\rk\Phi_1+\rk\Phi_2=\rk T_1+\rk T_2+\dots+\rk T_d,\tag\AHa$$ 
then
$$N_{\Phi_1*\Phi_2}(T_1,T_2,\dots,T_d)=\sum _{T'_1*T''_1=T_1,\dots,
T'_d*T''_d=T_d} ^{}
N_{\Phi_1}(T'_1,T'_2,\dots,T'_d)\cdot
N_{\Phi_2}(T''_1,T''_2,\dots,T''_d),
\tag\AH$$
where in the sum on the right-hand side any of $T'_1,T'_2,\dots,T'_d,
T''_1,T''_2,\dots,T''_d$ could also be empty.
\endproclaim

\remark{Remark} The reader should note that,
in order to have a non-vanishing summand in the sum
on the right-hand side of (\AH), we must necessarily have
$\rk\Phi_1=\rk T'_1+\rk T'_2+\dots+\rk T'_d$ and
$\rk\Phi_2=\rk T''_1+\rk T''_2+\dots+\rk T''_d$, because otherwise
one of $N_{\Phi_1}(T'_1,T'_2,\dots,T'_d)$ or
$N_{\Phi_2}(T''_1,T''_2,\dots,T''_d)$ vanishes.
\endremark
\demo{Proof of Proposition~\PD} 
By definition of the decomposition numbers
and by (\AHa), the number
$N_{\Phi_1*\Phi_2}(T_1,T_2,\dots,T_d)$ is equal to the number of
products 
$$c_1c_2\cdots c_d=c,\tag\AHb$$ 
where $c$ denotes the fixed Coxeter
element of type $\Phi_1*\Phi_2$, such that
$\ell_T(c_1)+\ell_T(c_2)+\dots+\ell_T(c_d)=\ell_T(c)$, and
such that the
type of $c_i$ as a parabolic Coxeter element is $T_i$, $i=1,2,\dots,d$.
Each of $c,c_1,c_2,\dots,c_d$ decomposes uniquely into a product of an element
of $W(\Phi_1)$ and an element of $W(\Phi_2)$, say 
$c=c'c''$,
$c_1=c_1'c_1''$,
$c_2=c_2'c_2''$, \dots,
$c_d=c_d'c_d''$ where $c',c_1',c_2',\dots,c_d'\in W(\Phi_1)$ and
$c'',c_1'',c_2'',\dots,c_d''\in W(\Phi_2)$. Since elements
of $W(\Phi_1)$ commute with elements of $W(\Phi_2)$ in
$W(\Phi_1*\Phi_2)=W(\Phi_1)\times W(\Phi_2)$, the relation (\AHb)
is equivalent to the two relations
$$\align c'_1c'_2\cdots c'_d&=c',\tag\AHc\\
c''_1c''_2\cdots c''_d&=c''.\tag\AHd
\endalign$$
Hence, we may compute $N_{\Phi_1*\Phi_2}(T_1,T_2,\dots,T_d)$
also by first determining the number 
of all possible products (\AHc) such that
the type of $c'_i$ is $T'_i$ (where we cover the case that $c'_i$ is
the identity element in $W(\Phi_1)$ by declaring that in that case
$T'_i=\emptyset$) and the number of all possible products (\AHd) such that
the type of $c''_i$ is $T''_i$ (with the same convention concerning 
identity elements), taking their product, and summing the results
over all possible $T'_1,\dots,T'_d,T''_1,\dots,T''_d$ with
$T'_1*T''_1=T_1$, \dots, $T'_d*T''_d=T_d$. This leads exactly to the
right-hand side of (\AH).\quad \quad \qed
\enddemo

\subhead 5. How to compute the characteristic polynomials\endsubhead
Aside from the decomposition numbers, the
second ingredient which we need for applying (\AC)
to compute the (dual) $M$-triangle of the $m$-divisible partitions
posets $NC^m(\Phi)$ is a list of the characteristic
polynomials $\chi^*_{NC(\Psi)}(y)$ for all 
{\it irreducible} root systems $\Psi$ of rank at most the rank of
$\Phi$. 
(By the multiplicativity of the characteristic polynomial,
this then gives also formulae for the characteristic polynomials
of all the reducible types.)
In fact, the numbers $N_\Psi(T_1,T_2,\dots,T_d)$ carry all the
information which is necessary to do this recursively. 
Namely, by the definition of $NC(\Psi)$ and of
the decomposition numbers
$N_\Psi(T_1,T_2,\dots,T_d)$, we have
$$\chi^*_{NC(\Psi)}(y)=
\sum _{T_1,T_2} ^{}N_{\Psi}(T_1,T_2)\,
\mu_{NC(T_2)}\!\left(\hat 0_{NC(T_2)},\hat 1_{NC(T_2)}\right)
y^{\rk T_1},\tag\AI$$
where $\mu_{NC(T_2)}(.,.)$ denotes the M\"obius function in
$NC(T_2)$,
and where $\hat 0_{NC(T_2)}$ and $\hat 1_{NC(T_2)}$ are, 
respectively, the minimal and the maximal element in $NC(T_2)$.
Indeed, inductively, the M\"obius functions 
$\mu_{NC(T_2)}\!\left(\hat 0_{NC(T_2)},\hat 1_{NC(T_2)}\right)$
are already known for all $T_2$ of lower rank than the rank
of $\Psi$. Hence, the only unknown in (\AI) is
$\mu_{NC(\Psi)}\!\left(\hat 0_{NC(\Psi)},\hat 1_{NC(\Psi)}\right)$.
However, the latter can be computed by setting $y=1$ in (\AI) and
using the fact that $\chi^*_{NC(\Psi)}(1)=0$ for all root
systems $\Psi$ of rank at least $1$. (This fact is
equivalent to the statement that 
$\sum _{u\in NC(\Psi)} ^{}\mu_{NC(\Psi)}\!
\left(u,\hat 1_{NC(\Psi)}\right)=0$, 
which is nothing but a part of the
definition of the M\"obius function.
Alternatively, one may use the uniform formula (\AGc) for the zeta 
polynomial of the non-crossing partition lattices, in which one
specializes the variable to $-1$, cf\. \cite{\StanAP, Sec.~3.11}.)

Below we list the values of the characteristic polynomials of the irreducible
root systems that we need in Sections~6 and 7 (and also for the
computations which are behind the numbers in the Appendix).
$$\allowdisplaybreaks\align \chi^*_{A_1}(y)&=y-1,\\
\chi^*_{A_2}(y)&=y^2-3y+2,\\
\chi^*_{A_3}(y)&=y^3-6y^2+10y-5,\\
\chi^*_{A_4}(y)&=y^4-10y^3+30y^2-35y+14,\\
\chi^*_{A_5}(y)&=y^5-15y^4+70y^3-140y^2+126y-42,\\
\chi^*_{A_6}(y)&=y^6-21y^5+140y^4-420y^3+630y^2-462y+132,\\
\chi^*_{A_7}(y)&=y^7-28y^6+252y^5-1050y^4+2310y^3-2772y^2+1716y-429,\\
\chi^*_{D_4}(y)&=y^4-12y^3+39y^2-48y+20,\\
\chi^*_{D_5}(y)&=y^5-20y^4+106y^3-230y^2+220y-77,\\
\chi^*_{D_6}(y)&=y^6-30y^5+235y^4-780y^3+1260y^2-980y+294,\\
\chi^*_{D_7}(y)&=y^7-42y^6+456y^5-2135y^4+5110y^3-6552y^2+4284y-1122,\\
\chi^*_{E_6}(y)&=y^6-36y^5+300y^4-1035y^3+1720y^2-1368y+418,\\
\chi^*_{E_7}(y)&=y^7-63y^6+777y^5-3927y^4+9933y^3-13299y^2+9009y-2431,\\
\chi^*_{E_8}(y)&= y^8-120y^7+2135y^6-15120y^5 \\
&\hskip2cm          +54327y^4-108360y^3+121555y^2-71760y+17342.
\tag\AJ\endalign$$

\subhead 6. The $M$-triangle of generalised non-crossing partitions 
of type $E_7$\endsubhead
We now implement the programme outlined in Sections~4 and 5
to compute the decomposition numbers for $E_7$ and, thus, 
via (\AC) and (\AJ), the $M$-triangle of the $m$-divisible
non-crossing partitions of type $E_7$.

To begin with, we have to compute the decomposition numbers for
types of smaller ranks, that is, of ranks $\le 6$. The decomposition
numbers for the irreducible types of rank $\le 6$ that we need, 
for $A_1$, $A_2$, $A_3$, $A_4$,
$A_5$, $A_6$, $D_4$, $D_5$, $D_6$, and $E_6$, are given
in the Appendix. Those for the reducible types of rank $\le 6$ then 
can be computed by using Proposition~\PD.\footnote{Although
theoretically totally sound, this is not the most convenient 
way to do it: the
actual way we computed these was to take the setup for an irreducible
root system of the same rank following the
programme outlined in Section~4, to change the values of the special
decomposition numbers given by Propositions~\PF\ and \PC, 
and to solve the new system of linear equations.}
Then we use (\AD) and (\AE) to express all the decomposition numbers
in terms of full rank decomposition numbers
$N_\Phi(T_1,T_2,\dots,T_d)$ (that is, those with $\rk\Phi=\rk T_1+\rk
T_2+\dots+\rk T_d$), in which the types $T_i$ are ordered
(for example, according to lexicographic order). 
Subsequently we produce the equations which one obtains from
(\AF) (using the earlier computed decomposition numbers of smaller
rank) and from comparing
coefficients of $m^iz^j$, $i,j\in\{0,1,\dots,7\}$, on both sides of
(\AG). Finally, we use Proposition~\PF\ to determine $N_{E_7}(E_7)$,
$N_{E_7}(A_1)$ and $N_{E_7}(A_1,A_1,A_1,A_1,A_1,A_1,A_1)$,
Proposition~\PC\ to compute the numbers
$N_{E_7}(T,A_1)$ with $\rk T=6$, and we use the last assertion
in Proposition~\PB\ for the type $A_1^5$ (which does not occur as the
type of a sub-diagram of the Dynkin diagram of $E_7$). To wit,
the latter leads to the following
special values of decomposition numbers:
$$\allowdisplaybreaks\gather
N_{E_7}(E_7)=1,\\
N_{E_7}(A_1)=63,\\
N_{E_7}(A_1,A_1,A_1,A_1,A_1,A_1,A_1)=1062882,\\
\matrix
N_{E_7}(E_6,A_1)=
N_{E_7}(D_6,A_1)=
N_{E_7}(A_6,A_1)=
N_{E_7}(A_1*A_5,A_1)\hfill\\
\kern2cm=
N_{E_7}(A_1*D_5,A_1)=
N_{E_7}(A_2*A_4,A_1)=
N_{E_7}(A_1*A_2*A_3,A_1)=9,\hfill
\endmatrix\\
\matrix
N_{E_7}(A_2*D_4,A_1)=
N_{E_7}(A_3^2,A_1)=
N_{E_7}(A_1^2*D_4,A_1)=
N_{E_7}(A_1^2*A_4,A_1)\hfill\\
\kern1cm=
N_{E_7}(A_1^3*A_3,A_1)=
N_{E_7}(A_2^3,A_1)=
N_{E_7}(A_1^2*A_2^2,A_1)=
N_{E_7}(A_1^4*A_2,A_1)\hfill\\
\kern1cm=
N_{E_7}(A_1^6,A_1)=
N_{E_7}(A_1^5,A_2)=
N_{E_7}(A_1^5,A_1^2)=
N_{E_7}(A_1^5,A_1,A_1)=0.\hfill
\endmatrix
\endgather$$

We now let {\sl Mathematica~5.2} solve the described system of
linear equations for the full rank decomposition
numbers.\footnote{\NoBlackBoxes 
The {\sl Mathematica} input is available at
{\tt http://www.mat.univie.ac.at/\~{}kratt/artikel/cluster2.html}.} 
Although this is a system of more than 200 equations
with 115 variables, the solution space is two-dimensional.
{\sl Mathematica} expresses all the variables in terms of
$X=N_{E_7}(A_1^4,A_1^3)$ and $Y=N_{E_7}(A_1^2*A_2,A_1^3)$.
For the subsequent considerations, we work with the following
selection of relations output by {\sl Mathematica}:
\def\BOX#1{\hbox to 20pt{\hss$\dsize#1$\hss}}
$$\align   N_{E_7}(A_5, A_2) &= \hphantom{-}\BOX{\frac {1272} {25}} 
  - \frac {58} {75} X- 
    \frac {58} {225} Y,\tag\AKa\\
N_{E_7}(A_4, A_3) &= \hphantom{-}\BOX{\frac {594} {25}} + \frac {18} {25} X+ 
    \hphantom{1}\frac {11} {25} Y,\tag\AKb\\
N_{E_7}(D_4, A_3) &= \hphantom{-}\BOX{\frac {126} {5}} 
  - \hphantom{1}\frac {2} {5} X- 
    \hphantom{1}\frac {7} {30} Y,\tag\AKc\\
  N_{E_7}(A_1^3*A_2, A_1^2) &= \hphantom{-}\BOX{\frac {423} {5}} - \frac {14} {5} X- 
    \hphantom{1}\frac {14} {15} Y,\tag\AKe\\
  N_{E_7}(A_1^3*A_2, A_2) &= -\BOX{\frac {192} {5}} + \frac {28} {15} X+ 
    \hphantom{1}\frac {28} {45} Y.\tag\AKf
\endalign$$
Relation (\AKa) implies
$$\align X-8Y&\equiv 2\pmod{25},\tag\AKg\\
\endalign$$
Relation (\AKb) implies
$$18X+11Y\equiv 6\pmod{25},\tag\AKi$$
and Relation (\AKc) implies
$$Y\equiv 0\pmod6.\tag\AKj$$
Solving (\AKg), (\AKi) and (\AKj) for $X$ and $Y$ yields
$$\align X&=25x+15y+4,\\
Y&=30y-6,
\endalign$$
for some integers $x$ and $y$ with $y>0$. 
If we substitute this in (\AKe) and
(\AKf), then we obtain
$$\align N_{E_7}(A_1^3*A_2, A_1^2) &=79-70(x+y),\\
N_{E_7}(A_1^3*A_2, A_2) &=\frac {4} {3}(35(x+y)-26).
\endalign$$
Since the decomposition numbers must be non-negative integers, we
conclude that $x+y=1$. Thus, we have $X=29-10y$ and $Y=30y-6$,
with $y=1$ or $y=2$. To decide which of the two values is the true
value, we appeal to Lemma~\LBa. Given a decomposition
$c=c_1c_2$, where $c_1$ is of type $A_1^4$ and $c_2$ is of type
$A_1^3$, the entire orbit 
$\Cal O=\{(c^kc_1c^{-k},c^kc_2c^{-k}):k\in\Bbb Z\}$
consists of pairs $(x,y)$ with $c=xy$, 
where $x$ is of type $A_1^4$ and $y$ is of type $A_1^3$.
Since Lemma~\LBa\ in the case of $E_7$ says that the orbits of
reflections have all size $9$, the size of $\Cal O$ must be a divisor
of $9$. We claim that it cannot be $1$.
For, in that case $c_1$ and $c_2$ commute with the Coxeter element
$c$. By a result of Carter \cite{\CartAA, Prop.~30},
the centralizer of $c$ is known to be the cyclic group
(of order $h=18$) generated
by $c$. Hence, $c_1$ and $c_2$ would have to be powers of
$c$. Since $c_1$ is of type $A_1^4$ and $c_2$ is of type $A_1^3$,
both $c_1$ and $c_2$ are of order $2$. However, the only element of
order $2$ in the cyclic group generated by $c$ is $c^9$, whence
$c^9=c_1=c_2$, which is absurd. The conclusion is that $3$ divides
the orbit $\Cal O$. Thus, $3$ must divide $X=N_{E_7}(A_1^4,A_1^3)$.
Out of the two possible values for $y$, we must therefore choose
$y=2$, and so $X=9$ and $Y=N_{E_7}(A_1^2*A_2,A_1^3)=54$. 

Now we substitute these values in the expressions in terms of $X$
and $Y$ found by {\sl Mathematica} for the other full rank decomposition
numbers. The result is that
$N_{E_7}(E_7)=1$, 
$N_{E_7}(E_6, A_1) = 9$,
$N_{E_7}(D_6, A_1) = 9$,
$N_{E_7}(A_6, A_1) = 9$,
$N_{E_7}(A_1*D_5, A_1) = 9$,
$N_{E_7}(A_1*A_5, A_1) = 9$,
$N_{E_7}(A_2*D_4, A_1) = 0$,
$N_{E_7}(A_2*A_4, A_1) = 9$,
$N_{E_7}(A_1^2*D_4, A_1) = 0$,
$N_{E_7}(A_1^2*A_4, A_1) = 0$,
$N_{E_7}(A_3^2, A_1) = 0$,
$N_{E_7}(A_1*A_2*A_3, A_1) = 9$,
$N_{E_7}(A_1^3*A_3, A_1) = 0$, 
$N_{E_7}(A_2^3, A_1) = 0$,
$N_{E_7}(A_1^2*A_2^2, A_1) = 0$,
$N_{E_7}(A_1^4*A_2, A_1) = 0$, 
$N_{E_7}(A_1^6, A_1) = 0$, 
$N_{E_7}(D_5, A_2) = 18$,
$N_{E_7}(A_5, A_2) = 30$,
$N_{E_7}(A_1*A_4, A_2) = 54$,
$N_{E_7}(A_1*D_4, A_2) = 9$,
$N_{E_7}(A_2*A_3, A_2) = 36$,
$N_{E_7}(A_1^2*A_3, A_2) = 36$,
$N_{E_7}(A_1*A_2^2, A_2) = 36$,
$N_{E_7}(A_1^3*A_2, A_2) = 12$,
$N_{E_7}(A_1^5, A_2) = 0$,
$N_{E_7}(D_5, A_1^2) = 54$,
$N_{E_7}(A_5, A_1^2) = 63$,
$N_{E_7}(A_1*D_4, A_1^2) = 27$,
$N_{E_7}(A_1*A_4, A_1^2) = 81$,
$N_{E_7}(A_2*A_3, A_1^2) = 27$,
$N_{E_7}(A_1^2*A_3, A_1^2) = 27$,
$N_{E_7}(A_1*A_2^2, A_1^2) = 27$, 
$N_{E_7}(A_1^3*A_2, A_1^2) = 9$,
$N_{E_7}(A_1^5, A_1^2) = 0$, 
$N_{E_7}(D_5, A_1, A_1) = 162$,
$N_{E_7}(A_5, A_1, A_1) = 216$,
$N_{E_7}(A_1*D_4, A_1, A_1) = 81$,
$N_{E_7}(A_1*A_4, A_1, A_1) = 324$,
$N_{E_7}(A_2*A_3, A_1, A_1) = 162$,
$N_{E_7}(A_1^2*A_3, A_1, A_1) = 162$,
$N_{E_7}(A_1*A_2^2, A_1, A_1) = 162$,
$N_{E_7}(A_1^3*A_2, A_1, A_1) = 54$,
$N_{E_7}(A_1^5, A_1, A_1) = 0$, 
$N_{E_7}(D_4,A_3)=9$, 
$N_{E_7}(A_4,A_3)=54$,
$N_{E_7}(A_1*A_3, A_3) = 135$,
$N_{E_7}(A_2^2, A_3) = 54$,
$N_{E_7}(A_1^2*A_2, A_3) = 162$,
$N_{E_7}(A_1^4, A_3) = 27$,
$N_{E_7}(D_4, A_1*A_2) = 45$,
$N_{E_7}(A_4, A_1*A_2) = 162$,
$N_{E_7}(A_1*A_3, A_1*A_2) = 243$,
$N_{E_7}(A_2^2, A_1*A_2) = 54$,
$N_{E_7}(A_1^2*A_2, A_1*A_2) = 162$,
$N_{E_7}(A_1^4, A_1*A_2) = 27$,
$N_{E_7}(D_4, A_1^3) = 30$,
$N_{E_7}(A_4, A_1^3) = 99$,
$N_{E_7}(A_1*A_3, A_1^3) = 126$,
$N_{E_7}(A_2^2, A_1^3) = 18$,
$N_{E_7}(A_1^2*A_2, A_1^3) = 54$,
$N_{E_7}(A_1^4, A_1^3) = 9$, 
$N_{E_7}(D_4, A_2, A_1) = 81$,
$N_{E_7}(A_4, A_2, A_1) = 378$,
$N_{E_7}(A_1*A_3, A_2, A_1) = 783$,
$N_{E_7}(A_2^2, A_2, \mathbreak A_1) = 270$,
$N_{E_7}(A_1^2*A_2, A_2, A_1) = 810$,
$N_{E_7}(A_1^4, A_2, A_1) = 135$,
$N_{E_7}(D_4, A_1^2, A_1) = 243$,
$N_{E_7}(A_4, A_1^2, A_1) = 891$,
$N_{E_7}(A_1*A_3, A_1^2, A_1) = 1377$,
$N_{E_7}(A_2^2, A_1^2, A_1) = 324$,\linebreak
$N_{E_7}(A_1^2*A_2, A_1^2, A_1) = 972$,
$N_{E_7}(A_1^4, A_1^2, A_1) = 162$,
$N_{E_7}(D_4, A_1, A_1, A_1) = 729$,
$N_{E_7}(A_4, \mathbreak A_1, A_1, A_1) = 2916$,
$N_{E_7}(A_1*A_3, A_1, A_1, A_1) = 5103$,
$N_{E_7}(A_2^2, A_1, A_1, A_1) = 1458$,\linebreak
$N_{E_7}(A_1^2*A_2, A_1, A_1, A_1) = 4374$,
$N_{E_7}(A_1^4, A_1, A_1, A_1) = 729$,
$N_{E_7}(A_3, A_3, A_1) = 486$,
$N_{E_7}(A_3, A_1*A_2, A_1) = 1458$,
$N_{E_7}(A_3, A_1^3, A_1) = 891$,
$N_{E_7}(A_1*A_2, A_1*A_2, A_1) = 2430$,
$N_{E_7}(A_1*A_2, A_1^3, A_1) = 1215$,
$N_{E_7}(A_1^3, A_1^3, A_1) = 540$,
$N_{E_7}(A_3, A_2, A_2) = 432$,\linebreak
$N_{E_7}(A_1*A_2, A_2, A_2) = 1188$,
$N_{E_7}(A_1^3, A_2, A_2) = 711$,
$N_{E_7}(A_3, A_2, A_1^2) = 1053$,\linebreak
$N_{E_7}(A_1*A_2, A_2, A_1^2) = 2349$,
$N_{E_7}(A_1^3, A_2, A_1^2) = 1323$,
$N_{E_7}(A_3, A_1^2, A_1^2) = 2430$,\linebreak
$N_{E_7}(A_1*A_2, A_1^2, A_1^2) = 3402$,
$N_{E_7}(A_1^3, A_1^2, A_1^2) = 1539$, 
$N_{E_7}(A_3, A_2, A_1, A_1) = 3402$,
$N_{E_7}(A_1*A_2, A_2, A_1, A_1) = 8262$,
$N_{E_7}(A_1^3, A_2, A_1, A_1) = 4779$,
$N_{E_7}(A_3, A_1^2, A_1, A_1) = 8019$,
$N_{E_7}(A_1*A_2, A_1^2, A_1, A_1) = 13851$,
$N_{E_7}(A_1^3, A_1^2, A_1, A_1) = 7047$,
$N_{E_7}(A_3, A_1, A_1, A_1, \mathbreak A_1) = 26244$,
$N_{E_7}(A_1*A_2, A_1, A_1, A_1, A_1) = 52488$,
$N_{E_7}(A_1^3, A_1, A_1, A_1, A_1) = 28431$,
$N_{E_7}(A_2, A_2, A_2, A_1) = 2916$,
$N_{E_7}(A_2, A_2, A_1^2, A_1) = 6561$,
$N_{E_7}(A_2, A_1^2, A_1^2, A_1) = 13122$,
$N_{E_7}(A_1^2, A_1^2, A_1^2, A_1) = 19683$, 
$N_{E_7}(A_2, A_2, A_1, A_1, A_1) = 21870$,
$N_{E_7}(A_2, A_1^2, A_1, A_1, \mathbreak A_1) = 45927$,
$N_{E_7}(A_1^2, A_1^2, A_1, A_1, A_1) = 78732$,
$N_{E_7}(A_2, A_1, A_1, A_1, A_1, A_1) = 157464$,
$N_{E_7}(A_1^2, A_1, A_1, A_1, A_1, A_1) = 295245$,
$N_{E_7}(A_1,A_1,A_1,A_1,A_1,A_1,A_1)=1062882$,
plus the assignments implied by
(\AD) and (\AE),
all other numbers $N_{E_7}(T_1,\dots,T_d)$ being zero. 

Finally, we substitute the values found for the decomposition numbers
and the formulae for the characteristic polynomials from (\AJ)
in (\AC), and we obtain 

{\eightpoint
$$\allowdisplaybreaks\multline 
(M^m_{E_7})^*\(x,y\)= 
\frac {1} {280}{{m( 9 m -2)  ( 9 m -4)  ( 9 m -5)  ( 9 m -8)  
      ( 3 m -1)  ( 3 m -2)    {x^7}{y^7}} }\\ -
\frac {9} {40}{{ {m^2}( 9 m -2)  ( 9 m -5)  ( 3 m -1)  ( 3 m -2)   
      ( 9 m -4)   {x^7}{y^6}}} \\+ 
\frac {3} {40}{{ {m^2}( 9 m -2)  ( 9 m -4)  ( 3 m -1)   
      ( 243 {m^2} -81m-14)   {x^7}{y^5}}} \\-
\frac {3} {8}{{ {m^2}( 3 m -1)  ( 9 m -2)  ( 9 m +2)  
      ( 81 {m^2} -9m-4)   {x^7}{y^4}}} \\+ 
\frac {3} {8}{{ {m^2}( 9 m -2)  ( 3 m +1) ( 9 m +2)  ( 81 {m^2} +9m-4)   {x^7}{y^3}}}
\\- 
\frac {3} {40}{{ {m^2}( 3 m +1)  ( 9 m +2)  ( 9 m +4)  ( 243 {m^2} +81m-14)  
       {x^7}{y^2}}} \\+ 
\frac {9} {40}{{ {m^2} ( 3 m +1)  ( 3 m +2)  ( 9 m +2)  ( 9 m +4)  
      ( 9 m +5)  {x^7}y}} \\-
\frac {1} {280}{{m ( 3 m +1)  ( 3 m +2)  ( 9 m +2)  ( 9 m +4)  
      ( 9 m +5)  ( 9 m +8)  {x^7}}} \\+
\frac {9} {40}{{ m( 9 m -2)  ( 9 m -5)  ( 3 m -1)  ( 3 m -2)  
      ( 9 m -4)   {x^6}{y^6}}} \\-
\frac {27} {20}{{ {m^2}( 27 m -13)  ( 9 m -2)  
  ( 9 m -4)  ( 3 m -1)  {x^6}{y^5}} } \\+ 
\frac {27} {8}{{ {m^2}( 3 m -1)  ( 9 m -2)  
      ( 243 {m^2} -45m-10)   {x^6}{y^4}}} -
\frac {27} {2}{{ {m^2}( 9 m -2) ( 9 m +2)  
      ( 27 {m^2} -1)   {x^6}{y^3}}} \\+ 
\frac {27} {8}{{ {m^2}( 3 m +1)   ( 9 m +2)  
      ( 243 {m^2} +45m-10)   {x^6}{y^2}}} -
\frac {27} {20}{{ {m^2} ( 3 m +1)  ( 9 m +2)  
  ( 9 m +4)  ( 27 m +13)  {x^6}y}} \\+ 
\frac {9} {40}{{ m ( 3 m +1) ( 3 m +2)  ( 9 m +2)  
      ( 9 m +4)  ( 9 m +5)  {x^6}}} + 
\frac {3} {40}{{ m( 207 m -103)  ( 9 m -2)  
  ( 9 m -4)  ( 3 m -1)    {x^5}{y^5}} } \\-
\frac {27} {8}{{ {m^2}( 207 m -71)  ( 3 m -1)   
      ( 9 m -2)   {x^5}{y^4}}} + 
\frac {9} {4}{{ {m^2}( 9 m -2)   
      ( 1863 {m^2} -144m-55)   {x^5}{y^3}}} \\-
\frac {9} {4}{{ {m^2}( 9 m +2)   
      ( 1863 {m^2} +144m-55)   {x^5}{y^2}}} + 
\frac {27} {8}{{ {m^2} ( 3 m +1)  ( 9 m +2)  
      ( 207 m +71) {x^5}y}} \\-
\frac {3} {40}{{ m ( 3 m +1)  ( 9 m +2)  
   (  9 m +4)  ( 207 m +103)  {x^5}}} + 
\frac {21} {8}{{ m( 63 m -23)  ( 3 m -1)   ( 9 m -2)   {x^4}{y^4}}}
\\-
\frac {189} {2}{{ {m^2}( 21 m -5)  ( 9 m -2)   {x^4}{y^3}}} + 
\frac {63} {4}{{ {m^2} ( 1701 {m^2} -37)   {x^4}{y^2}}} -
\frac {189} {2}{{ {m^2} ( 9 m +2)  ( 21 m +5) {x^4}y}} \\+ 
\frac {21} {8}{{ m ( 3 m +1)  ( 9 m+2 )  ( 63 m +23)  {x^4}}} + 
\frac {21} {2}{{ m( 27 m-7 )  ( 9 m -2)    {x^3}{y^3}}} -
\frac {189} {2}{{ {m^2} ( 81 m -13)   {x^3}{y^2}}} \\+ 
\frac {189} {2}{{ {m^2} (  81 m +13) {x^3}y}} -
\frac {21} {2}{{ m ( 9 m +2)  ( 27 m +7)  {x^3}}} + 
\frac {21} {2}{{m ( 63 m -11)   {x^2}{y^2}}} \\- 
1323 {m^2}{x^2}y + 
\frac {21} {2}{{ m ( 63 m +11)  {x^2}}} + 
63 m xy -
63 m x + 
1  .
\endmultline\tag{\tenpoint\AK}$$}%

\subhead 7. The $M$-triangle of generalised non-crossing partitions 
of type $E_8$\endsubhead
In this section we implement the programme outlined in Sections~4 and 5
to compute the decomposition numbers for $E_8$ and, thus, 
via (\AC) and (\AJ), the $M$-triangle of the $m$-divisible
non-crossing partitions of type $E_8$.

Again, to begin with, we have to compute the decomposition numbers for
types of smaller ranks, that is, of ranks $\le 7$. The decomposition
numbers for the irreducible types of rank $\le 7$ that we need, 
for $A_1$, $A_2$, $A_3$, $A_4$,
$A_5$, $A_6$, $A_7$, $D_4$, $D_5$, $D_6$, $D_7$, $E_6$, and $E_7$, are given
in the Appendix, respectively in Section~6. 
Those for the reducible types of rank $\le 7$ then 
can be computed by using Proposition~\PD.\footnote{See Footnote~1.}
Then we use (\AD) and (\AE) to express all the decomposition numbers
in terms of full rank decomposition numbers
$N_\Phi(T_1,T_2,\dots,T_d)$ (that is, those with $\rk\Phi=\rk T_1+\rk
T_1+\dots+\rk T_d$), in which the types $T_i$ are ordered
(for example, according to lexicographic order). 
Subsequently we produce the equations which one obtains from
(\AF) (using the earlier computed decomposition numbers of smaller
rank) and from comparing
coefficients of $m^iz^j$, $i,j\in\{0,1,\dots,8\}$, on both sides of
(\AG). Finally, we use Proposition~\PF\ to determine $N_{E_8}(E_8)$,
$N_{E_8}(A_1)$ and $N_{E_8}(A_1,A_1,A_1,A_1,A_1,A_1,A_1,A_1)$,
Proposition~\PC\ to compute the numbers
$N_{E_8}(T,A_1)$ with $\rk T=7$, and we use the last assertion
in Proposition~\PB\ for the types $A_1^6$, $A_1^5$, 
$A_1^3*A_3$, $A_2^3$ and $A_1^2*D_4$ (which do not occur as the
type of a sub-diagram of the Dynkin diagram of $E_8$). To wit,
the latter leads to the following
special values of decomposition numbers:
$$\allowdisplaybreaks\gather
N_{E_8}(E_8)=1,\\
N_{E_8}(A_1)=120,\\
N_{E_8}(A_1,A_1,A_1,A_1,A_1,A_1,A_1,A_1)=37968750,\\
\matrix
N_{E_8}(E_7,A_1)=
N_{E_8}(D_7,A_1)=
N_{E_8}(A_7,A_1)=
N_{E_8}(A_1*E_6,A_1)\hfill\\
\kern2cm=
N_{E_8}(A_1*A_6,A_1)=
N_{E_8}(A_2*D_5,A_1)\hfill\\
\kern2cm=
N_{E_8}(A_3*A_4,A_1)=
N_{E_8}(A_1*A_2*A_4,A_1)=15,\hfill
\endmatrix\\
\matrix
N_{E_8}(A_1*D_6,A_1)=
N_{E_8}(A_2*A_5,A_1)=
N_{E_8}(A_1^2*D_5,A_1)=
N_{E_8}(A_1^2*A_5,A_1)\hfill\\
\kern2cm=
N_{E_8}(A_3*D_4,A_1)=
N_{E_8}(A_1*A_2*D_4,A_1)=
N_{E_8}(A_1^3*D_4,A_1)\hfill\\
\kern2cm=
N_{E_8}(A_1^3*A_4,A_1)=
N_{E_8}(A_1*A_3^2,A_1)=
N_{E_8}(A_2^2*A_3,A_1)\hfill\\
\kern2cm=
N_{E_8}(A_1^2*A_2*A_3,A_1)=
N_{E_8}(A_1^4*A_3,A_1)=
N_{E_8}(A_1*A_2^3,A_1)\hfill\\
\kern2cm=
N_{E_8}(A_1^3*A_2^2,A_1)=
N_{E_8}(A_1^5*A_2,A_1)=
N_{E_8}(A_1^7,A_1)\hfill\\
\kern2cm=
N_{E_8}(A_1^6,A_2)=
N_{E_8}(A_1^6,A_1^2)=
N_{E_8}(A_1^6,A_1,A_1)=
N_{E_8}(A_1^5,A_3)\hfill\\
\kern2cm=
N_{E_8}(A_1^5,A_1*A_2)=
N_{E_8}(A_1^5,A_1^3)=
N_{E_8}(A_1^5,A_2,A_1)\hfill\\
\kern2cm=
N_{E_8}(A_1^5,A_1^2,A_1)=
N_{E_8}(A_1^5,A_1,A_1,A_1)=
N_{E_8}(A_1^3*A_3,A_2)\hfill\\
\kern2cm=
N_{E_8}(A_1^3*A_3,A_1^2)=
N_{E_8}(A_1^3*A_3,A_1,A_1)=
N_{E_8}(A_2^3,A_2)\hfill\\
\kern2cm=
N_{E_8}(A_2^3,A_1^2)=
N_{E_8}(A_1^4*A_2,A_2)=
N_{E_8}(A_1^4*A_2,A_1^2)\hfill\\
\kern2cm=
N_{E_8}(A_1^2*D_4,A_2)=
N_{E_8}(A_1^2*D_4,A_1^2)=0.\hfill
\endmatrix
\endgather$$

Again, we let {\sl Mathematica~5.2} solve the described system of
linear equations for the full rank decomposition
numbers.\footnote{\NoBlackBoxes 
The {\sl Mathematica} input is available at
{\tt http://www.mat.univie.ac.at/\~{}kratt/artikel/cluster2.html}.} 
Although this is now a system of more than 600 equations
with only about 250 variables, the situation here is even worse as 
the solution space is four-dimensional.
{\sl Mathematica} expresses all the variables in terms of
$X=N_{E_8}(A_5,\mathbreak A_1*A_2)$, $Y=N_{E_8}(D_5,A_1*A_2)$,
$N_{E_8}(A_4,A_1*A_3)$, and $N_{E_8}(D_4,A_4)$.
For the subsequent considerations, we work with the following
selection of relations output by {\sl Mathematica}:
$$\align   
N_{E_8}(A_5, A_3) &= \frac {1} {4}(750-X),\tag\ALa\\
N_{E_8}(D_5, A_3) &= \frac {1} {4}(375-Y),\tag\ALb\\
N_{E_8}(A_5, A_1^3) &= \frac {5} {6}(750-X),\tag\ALc\\
N_{E_8}(D_5, A_1^3) &= \frac {5} {6}(375-Y),\tag\ALd\\
N_{E_8}(A_2*A_3, A_1*A_2) &= \frac {3} {25}(4125-8X+16Y),\tag\ALe\\
N_{E_8}(A_2*A_3, A_1^3) &= \frac {1} {5}(1125+4X-8Y),\tag\ALf\\
N_{E_8}(A_1*D_4, A_3) &= Y-150,\tag\ALg\\
N_{E_8}(A_1^2*A_3, A_1^3) &= \frac {1} {5}(49875 - 56X - 138Y),\tag\ALh\\
N_{E_8}(A_1*A_4, A_1*A_2) &= 2 (2295 - 3 X - 4 Y ),\tag\ALgg\\
N_{E_8}(A_1^3*A_2, A_1^3) &= \frac {1} {15}(112 X + 226 Y - 86625).\tag\ALhh
\endalign
$$
There is one further relation which we shall make use of,
$$\align N_{E_8}(D_4)&=
N_{E_8}(D_4,D_4)+N_{E_8}(D_4,A_4)+N_{E_8}(D_4,A_3*A_1)\\
&\kern2cm+
N_{E_8}(D_4,A_2^2)+N_{E_8}(D_4,A_2*A_1^2)+N_{E_8}(D_4,A_1^4)\tag\ALi\\
&=\frac {27263} {168} - \frac {1} {40}N_{E_8}(A_4, A_1*A_3) 
+ \frac {1} {40} N_{E_8}(D_4, A_4) + \frac {283} {1500}X 
+   \frac {7957} {15750} Y,\tag\ALj
\endalign$$
where the first line is an instance of (\AE), and where the second
line follows from the first upon substituting the relations for
the full rank decomposition numbers output by {\sl Mathematica}.

\BlackBoxes
Relations (\ALa) and (\ALc) imply 
$$X\equiv 6\pmod{12},\tag\ALk$$
and Relations (\ALb) and (\ALd) imply 
$$Y\equiv 3\pmod{12}.\tag\ALl$$
Next, Relation (\ALe) implies 
$$X\equiv 2Y\pmod{25}.\tag\ALm$$
Solving (\ALk), (\ALl) and (\ALm) for $X$ and $Y$ yields
$$\align X&=300x+24y+6,\tag\ALmm\\
Y&=12y+3,\tag\ALmmm
\endalign$$
for some integers $x$ and $y$ with $y\ge0$. If we substitute this in
(\ALa), (\ALf)--(\ALhh), then we obtain
$$\align 
N_{E_8}(A_5, A_3) &=3(62-25x-2y),\tag\ALn\\
N_{E_8}(A_2*A_3, A_1^3) &=225+240x,\tag\ALo\\
N_{E_8}(A_1*D_4, A_3) &=12y-147,\tag\ALp\\
N_{E_8}(A_1^2*A_3, A_1^3) &=15(655-224x-40y),\tag\ALq\\
N_{E_8}(A_1*A_4, A_1*A_2) &= 30(151 - 60x - 8y),\tag\ALpp\\
N_{E_8}(A_1^3*A_2, A_1^3) &= 5(448x + 72y-1137).\tag\ALr
\endalign$$
From (\ALn) we infer $x\le 2$, while (\ALo) implies $x\ge0$.
Furthermore, from (\ALp) we infer $y\ge 13$, while (\ALq) implies
$y\le16$. If $x\ge1$, then (\ALpp) implies $y\le 91/8<13$, a
contradiction. Hence, $x=0$. In this case, Equation~(\ALr) implies
$y\ge 1137/72> 15$, and so $y=16$. Substituting this in
(\ALmm) and (\ALmmm), we obtain
$X=N_{E_8}(A_5,A_1*A_2)=390$ and $Y=N_{E_8}(D_5,A_1*A_2)=195$.

Unfortunately, similar arithmetic and positivity considerations do
not suffice to determine the remaining two ``free variables," the
decomposition numbers $N_{E_8}(A_4,A_1*A_3)$ and $N_{E_8}(D_4,A_4)$.
Instead, by a 12~days {\sl Maple} computation
using Stembridge's {\tt coxeter} package, we found that
$N_{E_8}(D_4)=325$ and $N_{E_8}(D_4,A_4)=15$. If we use these
values, together with the already determined values for $X$ and $Y$,
in (\ALj), we eventually find that $N_{E_8}(A_4,A_1*A_3)=390$.

We now substitute the above values for
$N_{E_8}(A_5,A_1*A_2)$, $N_{E_8}(D_5,A_1*A_2)$,
$N_{E_8}(A_4,\mathbreak A_1*A_3)$, and $N_{E_8}(D_4,A_4)$
in the expressions 
found by {\sl Mathematica} for the other full rank decomposition
numbers. The result is that
$N_{E_8}(E_8)= 1$, 
$N_{E_8}(E_7, A_1)= 15$,
$N_{E_8}(D_7, A_1)= 15$,
$N_{E_8}(A_7, A_1)= 15$,
$N_{E_8}(A_1*E_6, A_1)= 15$,
$N_{E_8}(A_1*D_6, A_1)= 0$,
$N_{E_8}(A_1*A_6, A_1)= 15$,
$N_{E_8}(A_2*D_5, A_1)= 15$,
$N_{E_8}(A_2*A_5, A_1)= 0$,
$N_{E_8}(A_1^2*D_5, A_1)= 0$,
$N_{E_8}(A_1^2*A_5, A_1)= 0$,
$N_{E_8}(A_3*D_4, A_1)= 0$,
$N_{E_8}(A_3*A_4, A_1)= 15$,
$N_{E_8}(A_1*A_2*D_4, A_1)= 0$,
$N_{E_8}(A_1*A_2*A_4, A_1)= 15$,
$N_{E_8}(A_1^3*D_4, A_1)= 0$,
$N_{E_8}(A_1^3*A_4, A_1)= 0$,
$N_{E_8}(A_1*A_3^2, A_1)= 0$,
$N_{E_8}(A_2^2*A_3, A_1)= 0$,
$N_{E_8}(A_1^2*A_2*A_3, A_1)= 0$,
$N_{E_8}(A_1^4*A_3, A_1)= 0$,
$N_{E_8}(A_1*A_2^3, A_1)= 0$,
$N_{E_8}(A_1^3*A_2^2, A_1)= 0$,
$N_{E_8}(A_1^5*A_2, A_1)= 0$,
$N_{E_8}(A_1^7, A_1)= 0$,
$N_{E_8}(E_6, A_2)= 20$,
$N_{E_8}(D_6, A_2)= 15$,
$N_{E_8}(A_6, A_2)= 60$,
$N_{E_8}(A_1*D_5, A_2)= 60$,
$N_{E_8}(A_1*A_5, A_2)= 60$,
$N_{E_8}(A_2*D_4, A_2)= 20$,
$N_{E_8}(A_2*A_4, A_2)= 90$,
$N_{E_8}(A_3^2, A_2)= 45$,
$N_{E_8}(A_1^2*D_4, A_2)= 0$,
$N_{E_8}(A_1^2*A_4, A_2)= 90$,
$N_{E_8}(A_1*A_2*A_3, A_2)= 90$,
$N_{E_8}(A_1^3*A_3, A_2)= 0$,
$N_{E_8}(A_2^3, A_2)= 0$,
$N_{E_8}(A_1^2*A_2^2, A_2)= 45$,
$N_{E_8}(A_1^4*A_2, A_2)= 0$,
$N_{E_8}(A_1^6, A_2)= 0$,
$N_{E_8}(E_6, A_1^2)= 45$,
$N_{E_8}(D_6, A_1^2)= 90$,
$N_{E_8}(A_6, A_1^2)= 135$,
$N_{E_8}(A_1*D_5, A_1^2)= 135$,
$N_{E_8}(A_1*A_5, A_1^2)= 135$,
$N_{E_8}(A_2*D_4, A_1^2)= 45$,
$N_{E_8}(A_2*A_4, A_1^2)= 90$,
$N_{E_8}(A_3^2, A_1^2)= 45$,
$N_{E_8}(A_1^2*D_4, A_1^2)= 0$,
$N_{E_8}(A_1^2*A_4, A_1^2)= 90$,
$N_{E_8}(A_1*A_2*A_3, A_1^2)= 90$,
$N_{E_8}(A_1^3*A_3, A_1^2)= 0$,
$N_{E_8}(A_2^3, A_1^2)= 0$,
$N_{E_8}(A_1^2*A_2^2, A_1^2)= 45$,
$N_{E_8}(A_1^4*A_2, A_1^2)= 0$,
$N_{E_8}(A_1^6, A_1^2)= 0$,
$N_{E_8}(E_6, A_1, A_1)= 150$,
$N_{E_8}(D_6, A_1, A_1)= 225$,
$N_{E_8}(A_6, A_1, A_1)= 450$,
$N_{E_8}(A_1*D_5, A_1, A_1)= 450$,
$N_{E_8}(A_1*A_5, A_1, A_1)= 450$,
$N_{E_8}(A_2*D_4, A_1, A_1)= 150$,
$N_{E_8}(A_2*A_4, A_1, A_1)= 450$,
$N_{E_8}(A_3^2, A_1, A_1)= 225$,
$N_{E_8}(A_1^2*D_4, A_1, A_1)= 0$,
$N_{E_8}(A_1^2*A_4, A_1, A_1)= 450$,
$N_{E_8}(A_1*A_2*A_3, A_1, A_1)= 450$,
$N_{E_8}(A_1^3*A_3, A_1, A_1)= 0$,
$N_{E_8}(A_2^3, A_1, A_1)= 0$,
$N_{E_8}(A_1^2*A_2^2, A_1, A_1)= 225$,
$N_{E_8}(A_1^4*A_2, A_1, A_1)= 0$,
$N_{E_8}(A_1^6, A_1, A_1)= 0$,
$N_{E_8}(D_5, A_3)= 45$,
$N_{E_8}(A_5, A_3)= 90$,
$N_{E_8}(A_1*A_4, A_3)= 315$,
$N_{E_8}(A_1*D_4, A_3)= 45$,
$N_{E_8}(A_2*A_3, A_3)= 270$,
$N_{E_8}(A_1^2*A_3, A_3)= 270$,
$N_{E_8}(A_1*A_2^2, A_3)= 225$,
$N_{E_8}(A_1^3*A_2, A_3)= 225$,
$N_{E_8}(A_1^5, A_3)= 0$,
$N_{E_8}(D_5, A_1*A_2)= 195$,
$N_{E_8}(A_5, A_1*A_2)= 390$,
$N_{E_8}(A_1*A_4, A_1*A_2)= 690$,
$N_{E_8}(A_1*D_4, A_1*A_2)= 195$,
$N_{E_8}(A_2*A_3, A_1*A_2)= 495$,
$N_{E_8}(A_1^2*A_3, A_1*A_2)= 495$,
$N_{E_8}(A_1*A_2^2, A_1*A_2)= 300$,
$N_{E_8}(A_1^3*A_2, A_1*A_2)= 300$,
$N_{E_8}(A_1^5, A_1*A_2)= 0$,
$N_{E_8}(D_5, A_1^3)= 150$,
$N_{E_8}(A_5, A_1^3)= 300$,
$N_{E_8}(A_1*A_4, A_1^3)= 375$,
$N_{E_8}(A_1*D_4, A_1^3)= 150$,
$N_{E_8}(A_2*A_3, A_1^3)= 225$,
$N_{E_8}(A_1^2*A_3, A_1^3)= 225$,
$N_{E_8}(A_1*A_2^2, A_1^3)= 75$,
$N_{E_8}(A_1^3*A_2, A_1^3)= 75$,
$N_{E_8}(A_1^5, A_1^3)= 0$,
$N_{E_8}(D_5, A_2, A_1)= 375$,
$N_{E_8}(A_5, A_2, A_1)= 750$,
$N_{E_8}(A_1*A_4, A_2, A_1)= 1950$,
$N_{E_8}(A_1*D_4, A_2, A_1)= 375$,
$N_{E_8}(A_2*A_3, A_2, A_1)= 1575$,
$N_{E_8}(A_1^2*A_3, A_2, A_1)= 1575$,
$N_{E_8}(A_1*A_2^2, A_2, A_1)= 1200$,
$N_{E_8}(A_1^3*A_2, A_2, A_1)= 1200$,
$N_{E_8}(A_1^5, A_2, A_1)= 0$,
$N_{E_8}(D_5, A_1^2, A_1)= 1125$,
$N_{E_8}(A_5, A_1^2, A_1)= 2250$,
$N_{E_8}(A_1*A_4, A_1^2, A_1)= 3825$,
$N_{E_8}(A_1*D_4, A_1^2, A_1)= 1125$,
$N_{E_8}(A_2*A_3, A_1^2, A_1)= 2700$,
$N_{E_8}(A_1^2*A_3, A_1^2, A_1)= 2700$,
$N_{E_8}(A_1*A_2^2, A_1^2, A_1)= 1575$,
$N_{E_8}(A_1^3*A_2, A_1^2, A_1)= 1575$,
$N_{E_8}(A_1^5, A_1^2, A_1)= 0$,
$N_{E_8}(D_5, A_1, A_1, A_1)= 3375$,
$N_{E_8}(A_5, A_1, A_1, A_1)= 6750$,
$N_{E_8}(A_1*A_4, A_1, A_1, A_1)= 13500$,
$N_{E_8}(A_1*D_4, A_1, A_1, A_1)= 3375$,
$N_{E_8}(A_2*A_3, A_1, A_1, A_1)= 10125$,
$N_{E_8}(A_1^2*A_3, A_1, A_1, A_1)= 10125$,
$N_{E_8}(A_1*A_2^2, A_1, A_1, A_1)= 6750$,
$N_{E_8}(A_1^3*A_2, A_1, A_1, A_1)= 6750$,
$N_{E_8}(A_1^5, A_1, A_1, A_1)= 0$,\linebreak
$N_{E_8}(D_4, D_4)= 5$,
$N_{E_8}(D_4, A_4)= 15$,
$N_{E_8}(A_4, A_4)= 138$,
$N_{E_8}(D_4, A_1*A_3)= 105$,
$N_{E_8}(A_4, A_1*A_3)= 390$,
$N_{E_8}(A_1*A_3, A_1*A_3)= 1155$,
$N_{E_8}(D_4, A_2^2)= 35$,
$N_{E_8}(A_4, A_2^2)= 180$,
$N_{E_8}(A_1*A_3, A_2^2)= 360$,
$N_{E_8}(A_2^2, A_2^2)= 95$,
$N_{E_8}(D_4, A_1^2*A_2)= 135$,
$N_{E_8}(A_4, A_1^2*A_2)= 630$,
$N_{E_8}(A_1*A_3, A_1^2*A_2)= 1035$,
$N_{E_8}(A_2^2, A_1^2*A_2)= 270$,
$N_{E_8}(A_1^2*A_2, A_1^2*A_2)= 495$,
$N_{E_8}(D_4, A_1^4)= 30$,
$N_{E_8}(A_4, A_1^4)= 165$,
$N_{E_8}(A_1*A_3, A_1^4)= 255$,
$N_{E_8}(A_2^2, A_1^4)= 60$,
$N_{E_8}(A_1^2*A_2, A_1^4)= 135$,
$N_{E_8}(A_1^4, A_1^4)= 30$,
$N_{E_8}(D_4, A_3, A_1)= 225$,
$N_{E_8}(A_4, A_3, A_1)= 1215$,
$N_{E_8}(A_1*A_3, A_3, A_1)= 4050$,
$N_{E_8}(A_2^2, A_3, A_1)= 1575$,
$N_{E_8}(A_1^2*A_2, A_3, A_1)= 5400$,
$N_{E_8}(A_1^4, A_3, A_1)= 1350$,
$N_{E_8}(D_4, A_1*A_2, A_1)= 975$,
$N_{E_8}(A_4, A_1*A_2, A_1)= 4590$,
$N_{E_8}(A_1*A_3, A_1*A_2, A_1)= 10800$,
$N_{E_8}(A_2^2, A_1*A_2, A_1)= 3450$,
$N_{E_8}(A_1^2*A_2, A_1*A_2, A_1)= 9900$,
$N_{E_8}(A_1^4, A_1*A_2, A_1)= 2475$,
$N_{E_8}(D_4, A_1^3, A_1)= 750$,
$N_{E_8}(A_4, A_1^3, A_1)= 3375$,
$N_{E_8}(A_1*A_3, A_1^3, A_1)= 6750$,
$N_{E_8}(A_2^2, A_1^3, A_1)= 1875$,
$N_{E_8}(A_1^2*A_2, A_1^3, A_1)= 4500$,
$N_{E_8}(A_1^4, A_1^3, A_1)= 1125$,
$N_{E_8}(D_4, A_2, A_2)= 175$,
$N_{E_8}(A_4, A_2, A_2)= 1140$,
$N_{E_8}(A_1*A_3, \mathbreak A_2, A_2)= 3300$,
$N_{E_8}(A_2^2, A_2, A_2)= 1300$,
$N_{E_8}(A_1^2*A_2, A_2, A_2)= 4500$,
$N_{E_8}(A_1^4, A_2, A_2)= 1125$,
$N_{E_8}(D_4, A_2, A_1^2)= 675$,
$N_{E_8}(A_4, A_2, A_1^2)= 3015$,
$N_{E_8}(A_1*A_3, A_2, A_1^2)= 8550$,
$N_{E_8}(A_2^2, A_2, A_1^2)= 2925$,
$N_{E_8}(A_1^2*A_2, A_2, A_1^2)= 9000$,
$N_{E_8}(A_1^4, A_2, A_1^2)= 2250$,
$N_{E_8}(D_4, \mathbreak A_1^2, A_1^2)= 1800$,
$N_{E_8}(A_4, A_1^2, A_1^2)= 8640$,
$N_{E_8}(A_1*A_3, A_1^2, A_1^2)= 17550$,
$N_{E_8}(A_2^2, A_1^2, \mathbreak A_1^2)= 5175$,
$N_{E_8}(A_1^2*A_2, A_1^2, A_1^2)= 13500$,
$N_{E_8}(A_1^4, A_1^2, A_1^2)= 3375$,
$N_{E_8}(D_4, A_2, A_1, \mathbreak A_1)= 1875$,
$N_{E_8}(A_4, A_2, A_1, A_1)= 9450$,
$N_{E_8}(A_1*A_3, A_2, A_1, A_1)= 27000$,
$N_{E_8}(A_2^2, A_2, \mathbreak A_1, A_1)= 9750$,
$N_{E_8}(A_1^2*A_2, A_2, A_1, A_1)= 31500$,
$N_{E_8}(A_1^4, A_2, A_1, A_1)= 7875$,
$N_{E_8}(D_4, \mathbreak A_1^2, A_1, A_1)= 5625$,
$N_{E_8}(A_4, A_1^2, A_1, A_1)= 26325$,
$N_{E_8}(A_1*A_3, A_1^2, A_1, A_1)= 60750$,
$N_{E_8}(A_2^2, A_1^2, A_1, A_1)= 19125$,
$N_{E_8}(A_1^2*A_2, A_1^2, A_1, A_1)= 54000$,
$N_{E_8}(A_1^4, A_1^2, A_1, A_1)= 13500$,
$N_{E_8}(D_4, A_1, A_1, A_1, A_1)= 16875$,
$N_{E_8}(A_4, A_1, A_1, A_1, A_1)= 81000$,
$N_{E_8}(A_1*A_3, \mathbreak A_1, A_1, A_1, A_1)= 202500$,
$N_{E_8}(A_2^2, A_1, A_1, A_1, A_1)= 67500$,
$N_{E_8}(A_1^2*A_2, A_1, A_1, A_1, \mathbreak A_1)= 202500$,
$N_{E_8}(A_1^4, A_1, A_1, A_1, A_1)= 50625$,
$N_{E_8}(A_3, A_3, A_2)= 1350$,
$N_{E_8}(A_3, \mathbreak A_1*A_2, A_2)= 5175$,
$N_{E_8}(A_3, A_1^3, A_2)= 3825$,
$N_{E_8}(A_1*A_2, A_1*A_2, A_2)= 15000$,\linebreak
$N_{E_8}(A_1*A_2, A_1^3, A_2)= 9825$,
$N_{E_8}(A_1^3, A_1^3, A_2)= 6000$,
$N_{E_8}(A_3, A_3, A_1^2)= 4050$,
$N_{E_8}(A_3, \mathbreak A_1*A_2, A_1^2)= 13500$,
$N_{E_8}(A_3, A_1^3, A_1^2)= 9450$,
$N_{E_8}(A_1*A_2, A_1*A_2, A_1^2)= 30825$,\linebreak
$N_{E_8}(A_1*A_2, A_1^3, A_1^2)= 17325$,
$N_{E_8}(A_1^3, A_1^3, A_1^2)= 7875$,
$N_{E_8}(A_3, A_3, A_1, A_1)= 12150$,
$N_{E_8}(A_3, A_1*A_2, A_1, A_1)= 42525$,
$N_{E_8}(A_3, A_1^3, A_1, A_1)= 30375$,
$N_{E_8}(A_1*A_2, A_1*A_2, A_1, \mathbreak A_1)= 106650$,
$N_{E_8}(A_1*A_2, A_1^3, A_1, A_1)= 64125$,
$N_{E_8}(A_1^3, A_1^3, A_1, A_1)= 33750$,
$N_{E_8}(A_3, \mathbreak A_2, A_2, A_1)= 10575$,
$N_{E_8}(A_3, A_2, A_1^2, A_1)= 29700$,
$N_{E_8}(A_3, A_1^2, A_1^2, A_1)= 76950$,\linebreak
$N_{E_8}(A_1*A_2, A_2, A_2, A_1)= 35700$,
$N_{E_8}(A_1*A_2, A_2, A_1^2, A_1)= 84825$,
$N_{E_8}(A_1*A_2, A_1^2, A_1^2, \mathbreak A_1)= 171450$,
$N_{E_8}(A_1^3, A_2, A_2, A_1)= 25125$,
$N_{E_8}(A_1^3, A_2, A_1^2, A_1)= 55125$,
$N_{E_8}(A_1^3, A_1^2, \mathbreak A_1^2, A_1)= 94500$,
$N_{E_8}(A_3, A_2, A_1, A_1, A_1)= 91125$,
$N_{E_8}(A_3, A_1^2, A_1, A_1, A_1)= 243000$,
$N_{E_8}(A_1*A_2, A_2, A_1, A_1, A_1)= 276750$,
$N_{E_8}(A_1*A_2, A_1^2, A_1, A_1, A_1)= 597375$,
$N_{E_8}(A_1^3, \mathbreak A_2, A_1, A_1, A_1)= 185625$,
$N_{E_8}(A_1^3, A_1^2, A_1, A_1, A_1)= 354375$,
$N_{E_8}(A_3, A_1, A_1, A_1, A_1, \mathbreak A_1)= 759375$,
$N_{E_8}(A_1*A_2, A_1, A_1, A_1, A_1, A_1)= 2025000$,
$N_{E_8}(A_1^3, A_1, A_1, A_1, A_1, A_1)= 1265625$,
$N_{E_8}(A_2, A_2, A_2, A_2)= 9350$,
$N_{E_8}(A_2, A_2, A_2, A_1^2)= 24975$,
$N_{E_8}(A_2, A_2, A_1^2, \mathbreak A_1^2)= 64350$,
$N_{E_8}(A_2, A_1^2, A_1^2, A_1^2)= 143100$,
$N_{E_8}(A_1^2, A_1^2, A_1^2, A_1^2)= 261225$,
$N_{E_8}(A_2, A_2, \mathbreak A_2, A_1, A_1)= 78000$,
$N_{E_8}(A_2, A_2, A_1^2, A_1, A_1)= 203625$,
$N_{E_8}(A_2, A_1^2, A_1^2, A_1, A_1)= \mathbreak 479250$,
$N_{E_8}(A_1^2, A_1^2, A_1^2, A_1, A_1)= 951750$,
$N_{E_8}(A_2, A_2, A_1, A_1, A_1, A_1)= 641250$,\linebreak
$N_{E_8}(A_2, A_1^2, A_1, A_1, A_1, A_1)= 1569375$,
$N_{E_8}(A_1^2, A_1^2, A_1, A_1, A_1, A_1)= 3341250$,
$N_{E_8}(A_2, \mathbreak A_1, A_1, A_1, A_1, A_1, A_1)= 5062500$,
$N_{E_8}(A_1^2, A_1, A_1, A_1, A_1, A_1, A_1)= 11390625$,
$N_{E_8}(A_1, \mathbreak A_1, A_1, A_1, A_1, A_1, A_1, A_1)= 37968750$,
plus the assignments implied by
(\AD) and (\AE),
all other numbers $N_{E_8}(T_1,\dots,T_d)$ being zero. 

Finally, we substitute the values found for the decomposition numbers 
and the formulae for the characteristic polynomials from (\AJ)
in (\AC), and we obtain 

{\eightpoint
$$\allowdisplaybreaks\multline 
(M^m_{E_8})^*\(x,y\)= 
\frac {1} {1344}{{m( 15 m -8)  ( 15 m -11)  ( 15 m -14)  ( 5 m -1)  
      ( 5 m -2)  ( 3 m -1)  ( 5 m -3)   {x^8}{y^8}}} \\-
\frac {5} {56}{{ {m^2}( 15 m -8)  ( 15 m -11)  ( 5 m -1)  ( 5 m -2)  
      ( 5 m -3)  ( 3 m -1)  {x^8}{y^7}} } \\+ 
\frac {5} {48}{{ {m^2}( 15 m -8)  ( 5 m -1)  ( 3 m -1) ( 5 m -2)  
      ( 225 {m^2} -90m-13)   {x^8}{y^6}}} \\-
\frac {15} {8}{{ {m^2}( 5 m -1)  ( 5 m -2)  ( 3 m -1)  ( 5 m +1)  
      ( 75 {m^2} -15m-4)   {x^8}{y^5}}} \\+ 
\frac {1} {32}{{{m^2} ( 5 m +1)  ( 5 m -1)   
      ( 84375 {m^4} - 12375 {m^2} +436)   {x^8}{y^4}}} \\- 
\frac {15} {8}{{ {m^2}( 5 m -1)  ( 3 m +1)  ( 5 m +1)  ( 5 m +2)  
      ( 75 {m^2}+15m-4 )   {x^8}{y^3}}} \\+ 
\frac {5} {48}{{ {m^2}( 3 m +1) ( 5 m +1)  ( 5 m +2)  ( 15 m +8)  
      ( 225 {m^2} +90m-13)   {x^8}{y^2}}} \\- 
\frac {5} {56}{{ {m^2} ( 3 m +1)  ( 5 m +1)  ( 5 m +2)  ( 5 m +3)  
      ( 15 m +8)  ( 15 m +11) {x^8}y}} \\+ 
\frac {1} {1344}{{m ( 3 m +1)  ( 5 m +1)  ( 5 m +2)  ( 5 m +3)  
      ( 15 m +8)  ( 15 m +11)  ( 15 m +14)  {x^8}}} \\+ 
\frac {5} {56}{{ m( 15 m -8)  ( 15 m -11)  ( 5 m -1)  
      ( 5 m -2)  ( 5 m -3)  ( 3 m -1)  {x^7}{y^7}}} \\-
\frac {25} {8}{{ {m^2}( 5 m -1)  ( 3 m -1)   
      ( 5 m -2)  {{( 15 m -8) }^2}  {x^7}{y^6}}} \\+ 
\frac {375} {8}{{ {m^2}( 5 m -1)  ( 5 m -2)  ( 3 m -1)   
      ( 45 {m^2} -12m-2)   {x^7}{y^5}}} \\-
\frac {15} {8}{{ {m^2}( 5 m +1)  ( 5 m -1)  
       ( 5625 {m^3} - 2250 {m^2} -75m+74)    {x^7}{y^4}}} \\+ 
\frac {15} {8}{{ {m^2}( 5 m -1)  ( 5 m +1)  
      ( 5625 {m^3} +2250m^2-75m-74)   {x^7}{y^3}} } \\-
\frac {375} {8}{{ {m^2}( 3 m +1)   ( 5 m +1)  ( 5 m +2)  
      ( 45 {m^2} +12m-2)   {x^7}{y^2}}} \\+ 
\frac {25} {8}{{ {m^2} ( 3 m +1)  ( 5 m +1)  
      ( 5 m +2)  {{( 15 m +8) }^2}{x^7}y} } \\-
\frac {5} {56}{{ m ( 3 m +1)  ( 5 m +1)  ( 5 m +2)  ( 5 m +3)  
      ( 15 m +8)  ( 15 m +11)  {x^7}}} \\+ 
\frac {5} {48}{{ m( 15 m -8)  ( 5 m -1)  ( 5 m -2)  ( 3 m -1)  
      ( 195 m -107)   {x^6}{y^6}}} \\-
\frac {375} {8}{{{m^2} ( 39 m -16)  ( 5 m -1)  ( 5 m -2)  
      ( 3 m -1)   {x^6}{y^5}}} \\+ 
\frac {5} {16}{{ {m^2} ( 5 m -1)  
  ( 219375 {m^3} - 103500 {m^2} + 3675 m +2342 )   {x^6}{y^4}}} \\-
\frac {75} {4}{{ {m^2}(5 m -1)  ( 5 m +1)  ( 975 {m^2} -32)   {x^6}{y^3}}
     } \\+ 
\frac {5} {16}{{ {m^2}( 5 m +1)  
  (  219375 {m^3} + 103500 {m^2} +3675 m -2342)  {x^6}{y^2}}} \\-
\frac {375} {8}{{ {m^2} ( 3 m +1)  ( 5 m +1)  
      ( 5 m +2)  ( 39 m +16) {x^6}y}} \\+ 
\frac {5} {48}{{ m ( 3 m +1)  ( 5 m +1)  ( 5 m +2)  ( 15 m +8)  
      ( 195 m +107)  {x^6}}} + 
15 m( 30 m -13)  ( 5 m -1)  ( 5 m -2)  ( 3 m -1)    {x^5}{y^5} \\-
15 {m^2} ( 5 m -1)  ( 2250 {m^2} - 1395 m+218 )   {x^5}{y^4} + 
75 {m^2}( 5 m -1)  ( 900 {m^2} -66m-23)   {x^5}{y^3} \\-
75 {m^2}( 5 m +1)  ( 900 {m^2} +66m-23)   {x^5}{y^2} + 
15 {m^2} ( 5 m +1)  (  2250 {m^2} + 1395 m+218) {x^5}y \\-
15 m ( 3 m +1)  ( 5 m +1)  ( 5 m +2)  ( 30 m +13)  {x^5} + 
\frac {1} {2}{{m ( 5 m -1)  (  10350 {m^2} - 6675 m+1084)   {x^4}{y^4}}}
\\-
15{m^2} ( 1380 m -307)  ( 5 m -1)  {x^4}{y^3} + 
5 {m^2} ( 31050 {m^2} -577)   {x^4}{y^2} -
15 {m^2} ( 5 m +1)  ( 1380 m +307) {x^4}y \\+ 
\frac {1} {2}{{m ( 5 m +1)  ( 10350 {m^2} + 6675 m+1084 )  {x^4}}} +
45 m( 45 m -11)  ( 5 m -1)  {x^3}{y^3} -
1125 {m^2} ( 27 m -4)   {x^3}{y^2} \\+ 
1125 {m^2} ( 27 m +4) {x^3}y -
45 m ( 5 m +1)  ( 45 m +11)  {x^3} + 
\frac {35} {2}{{ m ( 105 m -17)   {x^2}{y^2}}} -
3675 {m^2}{x^2}y \\+
\frac {35} {2}{{ m ( 105 m +17)  {x^2}}} + 
120 m xy - 
120 m x + 
1  .
\endmultline\tag{\tenpoint\AL}$$}%

\subhead 8. The $F=M$ Conjecture\endsubhead
In this section we describe briefly the $F=M$ Conjecture of
Armstrong \cite{\ArmDAA} predicting a surprising
relation between the $M$-triangle of the $m$-divisible partitions poset 
and the $F$-triangle of the generalised cluster complex of Fomin and
Reading \cite{\FoReAA}, we review the progress from \cite{\KratCB},
and we explain that our results from Sections~6 and 7 provide proofs
of the $F=M$ Conjecture for $E_7$ and $E_8$. If we put this together with
the results from \cite{\KratCB}, then it remains only the $D_n$ case of the
conjecture which is not completely proven.

For a non-negative integer $m$, the generalised cluster complex 
$\De^m(\Phi)$ is a certain simplicial complex on a certain set of
``coloured" roots, the roots being from $\Phi$. The precise definition
will not be important here, we refer the reader to \cite{\FoReAA, Sec.~2}.
The only fact which is important here is that some of the coloured
roots can be positive, others negative. Let $f_{k,l}(\Phi,m)$ denote
the number of faces of $\De^m(\Phi)$ which contain exactly $k$
positive and $l$ negative coloured roots. Define the $F$-triangle of
$\De^m(\Phi)$, denoted by $F^m_\Phi(x,y)$, as the two-variable polynomial
$$F^m_\Phi(x,y)=
\sum _{k,l\ge0} ^{}f_{k,l}(\Phi,m)\,x^ky^l.
\tag\AM$$
It is called ``triangle" because all faces have cardinality at most
$n$ and, thus, in the summation in (\AM) we can restrict the summation
indices to the triangle $k+l\le n$, $k,l\ge0$.

The generalised version of Chapoton's (ex-)conjecture \cite{\ChaFAA,
Conjecture~1}, due to Armstrong \cite{\ArmDAA, Conjecture~5.3.2}, 
is the following.

\proclaim{Conjecture FM}
For any finite root system $\Phi$ of rank $n$, we have
$$F^m_\Phi(x,y)=y^n\,M^m_\Phi\(\frac {1+y} {y-x},\frac {y-x} {y}\).
\tag\AMa$$
Equivalently, 
$$
(1-xy)^n F^m_\Phi\left(\frac {x(1+y)} {1-xy},\frac {xy}
{1-xy}\right)=(M^m_\Phi)^*(-x,-y).\tag\AN
$$ 
\endproclaim

It is easy to see (cf\. \cite{\KratCB, Sec.~8}) that 
it is enough to prove the conjecture for the irreducible
root systems. In \cite{\KratCB}, this has been done for the root
systems of type $A_n$, $B_n$, $I_2(a)$, $H_3$, $H_4$, $F_4$,
and $E_6$. The paper contains also
a partial proof for the root system of type $D_n$.
Given that we computed the $M$-triangle of the $m$-divisible
non-crossing partitions poset for $E_7$ and $E_8$ and that the
$F$-triangle of the generalised cluster complex has been computed in
\cite{\KratCB} for {\it all\/} irreducible root systems (thus, in
particular, for $E_7$ and $E_8$), the verification of (\AN) 
for $E_7$ and $E_8$ is pure routine on a computer algebra system.

If we combine these observations with the proof of the $m=1$ case of
the $F=M$ Conjecture by Athanasiadis \cite{\AthaAI} (an alternative
case-by-case proof is provided by the results in \cite{\KratCB} and
in the present paper), then we obtain the following theorem.

\proclaim{Theorem FM}
Conjecture FM is true with the possible exception of
root systems containing a copy of the root system $D_k$
for some $k$.
If $m=1$, Conjecture FM is true unconditionally.
\endproclaim

\subhead 9. A reciprocity phenomenon\endsubhead
By staring at the explicit expressions for the $M$-triangles 
of the $m$-divisible non-crossing partitions posets for the
irreducible root systems which we computed in \cite{\KratCB}
and in the present paper (in the $D_n$ case, this expression is only
conjectural), one observes a curious reciprocity relation between
the original $M$-triangle and the $M$-triangle with $m$ replaced by
$-m$. From the outset, the $M$-triangle with $-m$ in place of $m$ has
no combinatorial meaning since there is no meaning for
``$(-m)$-divisible non-crossing partitions." 
Nevertheless, it would be interesting to find an intrinsic explanation
of the phenomenon given in the theorem below. 
In any case, examples of situations
where combinatorial meaning was given to
non-combinatorial parameters are not so uncommon;
for example,
there exist reciprocity theorems for {\it $P$-partitions}
and the {\it Ehrhart quasi-polynomial} of polytopes
(cf.\ \cite{\StanAP, Theorems~4.5.7, 4.6.26, Cor.~4.5.15}),
which have their explanation in Stanley's reciprocity theorem
for linear homogeneous diophantine equations \cite{\StanAP,
Theorem~4.6.14}, and there exist reciprocity theorems for
monomer-dimer coverings of rectangles 
(cf.\ \cite{\AnBBAA, \PropAI, \StanBZ}).

\proclaim{Theorem \TA}
For any finite root system of rank $n$, 
except possibly for those which contain a copy of $D_k$
for some $k$, we have
$$
y^nM^{-m}_\Phi(xy,1/y)=M^m_\Phi(x,y).
\tag\AO$$
If the (conjectural) expression for $M^m_{D_n}(x,y)$ in 
\cite{\KratCB, Sec.~11, Prop.~D} is correct, then {\rm(\AO)} is 
true without any exceptions. 
Phrased differently, if the $F=M$ Conjecture is true, then
{\rm(\AO)} is true unconditionally.
\endproclaim

Formula (\AO) contains a reciprocity between numbers of 
maximal simplices in the generalised cluster complex
$\De^m(\Phi)$,
observed by Fomin and Reading \cite{\FoReAA, Eq.~(11.2) and
Prop.~11.4}, and revealed as an instance of Ehrhart
reciprocity by Athanasiadis and Tza\-na\-ki \cite{\AtTzAA, Sec.~7},
as a special case. (Again, for root systems
containing a copy of $D_k$ for some $k$, this is a conditional
statement.)
To see this, we first translate (\AO)
into a reciprocity relation for the $F$-triangle $F^m_\Phi(x,y)$
via (\AMa):
$$F^m_\Phi(x,y)=(1+x)^n\,
F^{-m}_\Phi\(-\frac {x} {1+x},\frac {y-x} {1+x}\).
\tag\AOa$$
If we compare coefficients of $x^ny^0$ on both sides of this
equation, then we obtain the relation
$$f_{n,0}(\Phi,m)=
\sum _{k=0} ^{n}(-1)^kf_k(\Phi,-m),\tag\AOb$$
where $f_k(\Phi,m)=
\sum _{l=0} ^{k}f_{l,k-l}(\Phi,m)$ is the {\it total\/} number
of faces in $\De^m(\Phi)$ containing exactly $k$ roots
(positive or negative). In the notation of \cite{\AtTzAA, \FoReAA},
we have $f_{n,0}(\Phi,m)=N^+(\Phi,m)$, and by \cite{\FoReAA, 
Eq.~(10.2)} the sum on the right-hand side of (\AOb) is
equal to $(-1)^nf_n(\Phi,-m-1)$ 
($(-1)^nN(\Phi,-m-1)$ in the notation of \cite{\FoReAA}).
Thus, we arrive at the reciprocity relation
$$N^+(\Phi,m)=(-1)^{n}N(\Phi,-m-1),$$
which is exactly \cite{\FoReAA, Eq.~(11.2)}.

The combinatorial significance of (\AOa) in general is less
clear. Comparison of coefficients of $x^ky^l$ on both sides
leads to the identity
$$f_{k,l}(\Phi,m)=
\sum _{r,s\ge0} ^{}(-1)^{r+s+l}\binom {n-r-s}{k+l-r-s}\binom sl
f_{r,s}(\Phi,-m).\tag\AOc$$

\subhead 10. A formula for the $A_n$ decomposition
numbers\endsubhead
Since the decomposition numbers $N_\Phi(T_1,T_2,\dots,T_d)$
carry so much information on the enumerative structure of
(ordinary and generalised) non-crossing partitions (see,
for example, Propositions~\PA\ and \PE), it would be 
of intrinsic interest to compute these numbers also for the
classical root systems. As a matter of fact, the type $A_n$ decomposition
numbers of full rank are known due to a result of
Goulden and Jackson \cite{\GoJaAS, Theorem~3.2} on the minimal factorization
of a long cycle. (The condition on the sum 
$l(\al_1)+l(\al_2)+\dots+l(\al_m)$ is misstated throughout the latter
paper. It should be replaced by
$l(\al_1)+l(\al_2)+\dots+l(\al_m)=(m-1)n+1$.) 
In the 
following theorem, we state this result in our language.

\proclaim{Theorem 9} 
Let $\Phi=A_n$, and let $T_1,T_2,\dots,T_d$ be types
with $\rk T_1+\rk T_2+\dots+\rk T_d=n$,
where 
$$T_i=A_1^{m_1^{(i)}}*A_2^{m_2^{(i)}}*\dots*A_n^{m_n^{(i)}},\quad 
i=1,2,\dots,d.$$ 
Then
$$
N_{A_n}(T_1,T_2,\dots,T_d)=(n+1)^{d-1}\prod _{i=1} ^{d}
\frac {1} {n-\rk T_i+1}\binom {n-\rk T_i+1}{m_1^{(i)},m_2^{(i)},\dots,
m_n^{(i)}},
\tag\AOe$$
where the multinomial coefficient
is defined by
$$\binom {M}{m_1,m_2,\dots,m_n}=\frac {M!} {m_1!\,m_2!\cdots
m_n!\,(M-m_1-m_2-\dots-m_n)!}.$$
\endproclaim

This result allows one to derive a compact formula
for {\it all\/} type $A_n$ decomposition numbers.

\proclaim{Theorem 10}
Let $\Phi=A_n$, and let the types $T_1,T_2,\dots,T_d$ be given,
where 
$$T_i=A_1^{m_1^{(i)}}*A_2^{m_2^{(i)}}*\dots*A_n^{m_n^{(i)}},\quad 
i=1,2,\dots,d.$$ 
Then
$$\multline 
N_{A_n}(T_1,T_2,\dots,T_d)=(n+1)^{d-1}
\binom {n+1}{\rk T_1+\rk T_2+\dots+\rk T_d+1}\\
\times
\prod _{i=1} ^{d}
\frac {1} {n-\rk T_i+1}\binom {n-\rk T_i+1}{m_1^{(i)},m_2^{(i)},\dots,
m_n^{(i)}}.
\endmultline\tag\AOd$$
\endproclaim

\demo{Proof} 
If we write $r$ for $n-\rk T_1-\rk T_2-\dots-\rk T_d$,
then the relation (\AE) becomes
$$N_{A_n}(T_1,T_2,\dots,T_d)=
{\sum _{T:\rk T=r}}
^{}N_{A_n}(T_1,T_2,\dots,T_d,T).
$$
Upon letting $T=A_1^{m_1}*A_2^{m_2}*\dots*A_n^{m_n}$,
substitution of (\AOe) in the above equation yields
$$\align
N_{A_n}(T_1,T_2,\dots,T_d)&=\sum _{m_1+2m_2+\dots+nm_n=r} ^{}
(n+1)^{d}\frac {1} {n-r+1}\binom {n-r+1}{m_1,m_2,\dots,m_n}\\
&\kern2cm\cdot
\prod _{i=1} ^{d}
\frac {1} {n-\rk T_i+1}\binom {n-\rk T_i+1}{m_1^{(i)},m_2^{(i)},\dots,
m_n^{(i)}}.
\endalign$$
Now, by comparison of coefficients of $z^r$ on both sides of
$$(1+z+z^2+z^3+\cdots)^M=(1-z)^{-M},$$
we infer
$$\sum _{m_1+2m_2+\dots+rm_r=r} ^{}\binom M{m_1,m_2,\dots,m_r}=
\binom {M+r-1}r.
$$
If we use this identity with $M=n-r+1$, we obtain our claim after
little simplification.\quad \quad \qed
\enddemo

The reader should note that formula (\AOd) generalises Kreweras'
formula \cite{\KrewAC, Theorem~4} 
for the number of non-crossing partitions of $n+1$
with given block sizes, to which it reduces for $d=1$.

\subhead Appendix: The decomposition numbers for types $A_n$ and
$D_n$ for small $n$, and for $E_6$\endsubhead
In this appendix we list the full rank decomposition numbers,
that is, the decomposition numbers
$N_{\Phi}(T_1,T_2,\dots,T_d)$, where
$\rk T_1+\rk T_2+\dots+\rk T_d=\rk \Phi$, 
for $\Phi=A_1,A_2,A_3,A_4,A_5,A_6,A_7,D_4,D_5,D_6,D_7,E_6$.
These numbers are required for setting up the system of equations
(\AF) in Sections~6 and 7. We do not list decomposition numbers
which are zero. It was explained in Section~4 how to compute
these numbers by setting up a system of linear equations
for them. Whenever, aside from the listing of the numbers, 
there is no further comment, then the computation of the decomposition
numbers is either
trivial (for $A_1$ and $A_2$), or the equations in Propositions~\PB\
and \PE, together with the assignments from Propositions~\PF\ and
\PC\ determine them already uniquely. If not, then we mention the
additional assignments, respectively considerations, which are
required for the computation.\footnote{The {\sl Mathematica} inputs
for the computations are available at\newline
{\tt http://www.mat.univie.ac.at/\~{}kratt/artikel/cluster2.html}.} 
Clearly, the decomposition numbers for $A_n$, $n=1,2,\dots,$ can be
computed directly from Theorem~9.

\NoBlackBoxes
\subsubhead The decomposition numbers for $A_1$\endsubsubhead
$N_{A_1}(A_1)=1$.

\subsubhead The decomposition numbers for $A_2$\endsubsubhead
$N_{A_2}(A_2)=1$, 
$N_{A_2}(A_1,A_1)=3$.

\subsubhead The decomposition numbers for $A_3$\endsubsubhead
$N_{A_3}(A_3)=1$, 
$N_{A_3}(A_2,A_1)=4$,
$N_{A_3}(A_1^2,A_1)=2$,
$N_{A_3}(A_1,A_1,A_1)=16$.

\subsubhead The decomposition numbers for $A_4$\endsubsubhead
$N_{A_4}(A_4)=1$, 
$N_{A_4}(A_3,A_1)=5$,
$N_{A_4}(A_1*A_2,A_1)=5$,
$N_{A_4}(A_2,A_2)=5$,
$N_{A_4}(A_2,A_1^2)=5$,
$N_{A_4}(A_1^2,A_1^2)=5$,
$N_{A_4}(A_2,A_1,A_1)=25$,
$N_{A_4}(A_1^2,\mathbreak A_1,A_1)=25$,
$N_{A_4}(A_1,A_1,A_1,A_1)=125$.

\subsubhead The decomposition numbers for $A_5$\endsubsubhead
$N_{A_5}(A_5)=1$, 
$N_{A_5}(A_4,A_1)=6$,
$N_{A_5}(A_1*A_3,A_1)=6$,
$N_{A_5}(A_2^2,A_1)=3$,
$N_{A_5}(A_3,A_2)=6$,
$N_{A_5}(A_1*A_2,A_2)=12$,
$N_{A_5}(A_1^3,A_2)=2$,
$N_{A_5}(A_3,\mathbreak A_1^2)=9$,
$N_{A_5}(A_1*A_2,A_1^2)=18$,
$N_{A_5}(A_1^3,A_1^2)=3$,
$N_{A_5}(A_3,A_1,A_1)=36$,
$N_{A_5}(A_1*A_2,\mathbreak A_1,A_1)=72$,
$N_{A_5}(A_1^3,A_1,A_1)=12$,
$N_{A_5}(A_2,A_2,A_1)=36$,
$N_{A_5}(A_2,A_1^2,A_1)=54$,
$N_{A_5}(A_1^2,A_1^2,A_1)=81$,
$N_{A_5}(A_2,A_1,A_1,A_1)=216$,
$N_{A_5}(A_1^2,A_1,A_1,A_1)=324$,
$N_{A_5}(A_1,\mathbreak A_1,A_1,A_1,A_1)=1296$.

\subsubhead The decomposition numbers for $A_6$\endsubsubhead
$N_{A_6}(A_6)=1$, 
$N_{A_6}(A_5, A_1) = 7$,  
$N_{A_6}(A_1*A_4, A_1) = 7$, 
$N_{A_6}(A_2*A_3, A_1) = 7$, 
$N_{A_6}(A_4, A_2) = 7$, 
$N_{A_6}(A_1*A_3, A_2) = 14$, 
$N_{A_6}(A_2^2, A_2) = 7$, 
$N_{A_6}(A_1^2*A_2, A_2) = 7$, 
$N_{A_6}(A_4, A_1^2) = 14$, 
$N_{A_6}(A_1*A_3, A_1^2) = 28$, 
$N_{A_6}(A_2^2, A_1^2) = 14$, 
$N_{A_6}(A_1^2*A_2, A_1^2) = 14$, 
$N_{A_6}(A_3, A_3) = 7$, 
$N_{A_6}(A_3, A_1*A_2) = 21$, 
$N_{A_6}(A_3, A_1^3) = 7$, 
$N_{A_6}(A_1*A_2, A_1*A_2) = 63$, 
$N_{A_6}(A_1*A_2, A_1^3) = 21$,
$N_{A_6}(A_1^3, A_1^3) = 7$, 
$N_{A_6}(A_4, A_1, A_1) = 49$, 
$N_{A_6}(A_1*A_3, A_1, A_1) = 98$, 
$N_{A_6}(A_2^2, A_1, A_1) = 49$, 
$N_{A_6}(A_1^2*A_2, A_1, A_1) = 49$, 
$N_{A_6}(A_3, A_2, A_1) = 49$, 
$N_{A_6}(A_3, A_1^2, A_1) = 98$, 
$N_{A_6}(A_1*A_2, A_2, A_1) = 147$, 
$N_{A_6}(A_1*A_2, A_1^2, A_1) = 294$,
$N_{A_6}(A_1^3, A_2, A_1) = 49$, 
$N_{A_6}(A_1^3, A_1^2, A_1) = 98$, 
$N_{A_6}(A_3, A_1, A_1, A_1) = 343$, 
$N_{A_6}(A_1*A_2, A_1, A_1, A_1) = 1029$, 
$N_{A_6}(A_1^3, A_1, A_1, A_1) = 343$, 
$N_{A_6}(A_2, A_2, A_2) = 49$,
$N_{A_6}(A_2, A_2, A_1^2) = 98$, 
$N_{A_6}(A_2, A_1^2, A_1^2) = 196$, 
$N_{A_6}(A_1^2, A_1^2, A_1^2) = 392$, 
$N_{A_6}(A_2, A_2, A_1,\mathbreak A_1) = 343$,
$N_{A_6}(A_2, A_1^2, A_1, A_1) = 686$, 
$N_{A_6}(A_1^2, A_1^2, A_1, A_1) = 1372$, 
$N_{A_6}(A_2, A_1, A_1, A_1,\mathbreak A_1) = 2401$,
$N_{A_6}(A_1^2, A_1, A_1, A_1, A_1) = 4802$, 
$N_{A_6}(A_1, A_1, A_1, A_1, A_1, A_1) = 16807$.

\subsubhead The decomposition numbers for $A_7$\endsubsubhead
$N_{A_7}(A_7)=1$, 
$N_{A_7}(A_6, A_1) = 8$,
$N_{A_7}(A_1*A_5, A_1) = 8$,
$N_{A_7}(A_2*A_4, A_1) = 8$,
$N_{A_7}(A_3^2, A_1) = 4$,
$N_{A_7}(A_5, A_2) = 8$,
$N_{A_7}(A_1*A_4, A_2) = 16$,
$N_{A_7}(A_2*A_3, A_2) = 16$,
$N_{A_7}(A_1^2*A_3, A_2) = 8$,
$N_{A_7}(A_1*A_2^2, A_2) = 8$,
$N_{A_7}(A_5, A_1^2) = 20$,
$N_{A_7}(A_1*A_4, A_1^2) = 40$,
$N_{A_7}(A_2*A_3, A_1^2) = 40$,
$N_{A_7}(A_1^2*A_3, A_1^2) = 20$,
$N_{A_7}(A_1*A_2^2, A_1^2) = 20$,
$N_{A_7}(A_5, A_1, A_1) = 64$,
$N_{A_7}(A_1*A_4, A_1, A_1) = 128$,
$N_{A_7}(A_2*A_3, A_1, A_1) = 128$,
$N_{A_7}(A_1^2*A_3, A_1, A_1) = 64$,
$N_{A_7}(A_1*A_2^2, A_1, A_1) = 64$,
$N_{A_7}(A_4, A_3) = 8$,
$N_{A_7}(A_1*A_3, A_3) = 24$,
$N_{A_7}(A_2^2, A_3) = 12$,
$N_{A_7}(A_1^2*A_2, A_3) = 24$,
$N_{A_7}(A_1^4, A_3) = 2$,
$N_{A_7}(A_4, A_1*A_2) = 32$,
$N_{A_7}(A_1*A_3, A_1*A_2) = 96$,
$N_{A_7}(A_2^2, A_1*A_2) = 48$,
$N_{A_7}(A_1^2*A_2, A_1*A_2) = 96$,
$N_{A_7}(A_1^4, A_1*A_2) = 8$,
$N_{A_7}(A_4, A_1^3) = 16$,
$N_{A_7}(A_1*A_3, A_1^3) = 48$,
$N_{A_7}(A_2^2, A_1^3) = 24$,
$N_{A_7}(A_1^2*A_2, A_1^3) = 48$,
$N_{A_7}(A_1^4, A_1^3) = 4$,
$N_{A_7}(A_4, A_2, A_1) = 64$,
$N_{A_7}(A_1*A_3, A_2, A_1) = 192$,
$N_{A_7}(A_2^2, A_2, A_1) = 96$,
$N_{A_7}(A_1^2*A_2, A_2, A_1) = 192$,
$N_{A_7}(A_1^4, A_2, A_1) = 16$,
$N_{A_7}(A_4,\mathbreak A_1^2, A_1) = 160$,
$N_{A_7}(A_1*A_3, A_1^2, A_1) = 480$,
$N_{A_7}(A_2^2, A_1^2, A_1) = 240$,
$N_{A_7}(A_1^2*A_2, A_1^2, A_1) = 480$,
$N_{A_7}(A_1^4, A_1^2, A_1) = 40$,
$N_{A_7}(A_4, A_1, A_1, A_1) = 512$,
$N_{A_7}(A_1*A_3, A_1, A_1, A_1) = 1536$,
$N_{A_7}(A_2^2, A_1, A_1, A_1) = 768$,
$N_{A_7}(A_1^2*A_2, A_1, A_1, A_1) = 1536$,
$N_{A_7}(A_1^4, A_1, A_1, A_1) = 128$,
$N_{A_7}(A_3, A_3, A_1) = 64$,
$N_{A_7}(A_3, A_1*A_2, A_1) = 256$,
$N_{A_7}(A_3, A_1^3, A_1) = 128$,
$N_{A_7}(A_1*A_2, A_1*A_2, A_1) = 1024$,
$N_{A_7}(A_1*A_2, A_1^3, A_1) = 512$,
$N_{A_7}(A_1^3, A_1^3, A_1) = 256$,
$N_{A_7}(A_3, A_2,\mathbreak A_2) = 64$,
$N_{A_7}(A_3, A_2, A_1^2) = 160$,
$N_{A_7}(A_3, A_1^2, A_1^2) = 400$,
$N_{A_7}(A_1*A_2, A_2, A_2) = 256$,
$N_{A_7}(A_1*A_2, A_2, A_1^2) = 640$,
$N_{A_7}(A_1*A_2, A_1^2, A_1^2) = 1600$,
$N_{A_7}(A_1^3, A_2, A_2) = 128$,
$N_{A_7}(A_1^3, A_2, A_1^2) = 320$,
$N_{A_7}(A_1^3, A_1^2, A_1^2) = 800$,
$N_{A_7}(A_3, A_2, A_1, A_1) = 512$,
$N_{A_7}(A_3, A_1^2,\mathbreak A_1, A_1) = 1280$,
$N_{A_7}(A_1*A_2, A_2, A_1, A_1) = 2048$,
$N_{A_7}(A_1*A_2, A_1^2, A_1, A_1) = 5120$,
$N_{A_7}(A_1^3, A_2, A_1, A_1) = 1024$,
$N_{A_7}(A_1^3, A_1^2, A_1, A_1) = 2560$,
$N_{A_7}(A_3, A_1, A_1, A_1, A_1) =\mathbreak 4096$,
$N_{A_7}(A_1*A_2, A_1, A_1, A_1, A_1) = 16384$,
$N_{A_7}(A_1^3, A_1, A_1, A_1, A_1) = 8192$,
$N_{A_7}(A_2, A_2,\mathbreak A_2, A_1) = 512$,
$N_{A_7}(A_2, A_2, A_1^2, A_1) = 1280$,
$N_{A_7}(A_2, A_1^2, A_1^2, A_1) = 3200$,
$N_{A_7}(A_1^2, A_1^2, A_1^2,\mathbreak A_1) = 8000$,
$N_{A_7}(A_2, A_2, A_1, A_1, A_1) = 4096$,
$N_{A_7}(A_2, A_1^2, A_1, A_1, A_1) = 10240$,
$N_{A_7}(A_1^2,\mathbreak A_1^2, A_1, A_1, A_1) = 25600$,
$N_{A_7}(A_2, A_1, A_1, A_1, A_1, A_1) = 32768$,
$N_{A_7}(A_1^2, A_1, A_1, A_1, A_1,\mathbreak A_1) = 81920$,
$N_{A_7}(A_1, A_1, A_1, A_1, A_1, A_1, A_1) = 262144$.

\subsubhead The decomposition numbers for $D_4$\endsubsubhead
$N_{D_4}(D_4)=1$, 
$N_{D_4}(A_3,A_1)=9$,
$N_{D_4}(A_1^3,A_1)=3$,
$N_{D_4}(A_2,A_2)=6$,
$N_{D_4}(A_2,A_1^2)=9$,
$N_{D_4}(A_2,A_1,A_1)=36$,
$N_{D_4}(A_1^2,A_1,A_1)=27$,
$N_{D_4}(A_1,A_1,A_1,A_1)=162$.

\subsubhead The decomposition numbers for $D_5$\endsubsubhead
$N_{D_5}(D_5)=1$, 
$N_{D_5}(D_4,A_1)=4$,
$N_{D_5}(A_4,A_1)=8$,
$N_{D_5}(A_1*A_3,A_1)=4$,
$N_{D_5}(A_1^2*A_2,A_1)=4$,
$N_{D_5}(A_3,A_2)=12$,
$N_{D_5}(A_1*A_2,A_2)=16$,
$N_{D_5}(A_1^3,A_2)=8$,
$N_{D_5}(A_3,A_1^2)=22$,
$N_{D_5}(A_1*A_2,A_1^2)=8$,
$N_{D_5}(A_1^3,A_1^2)=4$,
$N_{D_5}(A_3,A_1,A_1)=80$,
$N_{D_5}(A_1*A_2,A_1,A_1)=64$,
$N_{D_5}(A_1^3,A_1,A_1)=32$,
$N_{D_5}(A_2,A_2,\mathbreak A_1)=64$,
$N_{D_5}(A_2,A_1^2,A_1)=96$,
$N_{D_5}(A_1^2,A_1^2,A_1)=80$,
$N_{D_5}(A_2,A_1,A_1,A_1)=384$,
$N_{D_5}(A_1^2,A_1,A_1,A_1)=448$,
$N_{D_5}(A_1,A_1,A_1,A_1,A_1)=2048$.

\subsubhead The decomposition numbers for $D_6$\endsubsubhead
Unfortunately, the solution space of the system of linear equations 
does not consist of a unique solution but is two-dimensional.\linebreak
{\sl Mathematica~5.2} expresses all the decomposition numbers in terms of
$X=N_{D_6}(A_1^3,A_1^3)$ and
$Y=N_{D_6}(A_1^4,A_1^2)$. In the sequel we make use of the two
relations
$$\align   
  N_{D_6}(A_1^2*A_2, A_1^2) &= 25 - 
\hphantom{1}\frac {27}{40}X - 3Y, \\
N_{D_6}(A_3, A_3) &= 20 + \frac9{100}X.
\endalign$$
They imply $X\equiv0$ (mod $200$), $X\le 1000/27<200$,
$Y\equiv2$ (mod $3$), and $Y\le 25/3<9$.
Hence, we have $X=N_{D_6}(A_1^3,A_1^3)=0$, and 
the three possibilities $Y=2,5,8$.
To decide which of the three values is the true
value, we argue as in Section~6 when we had to decide
which of two possible values for $N_{E_7}(A_1^2*A_2,A_1^3)$
was the correct one (see the paragraph after (\AKj)).
Since Lemma~\LBa\ in the case of $D_6$ says that the orbits of
reflections have all size $5$, here the conclusion is that
$5$ must divide $Y=N_{D_6}(A_1^4,A_1^2)$. Thus,
$N_{D_6}(A_1^4,A_1^2)=5$.
If we substitute this in the relations for the decomposition numbers
found by {\sl Mathematica}, then we obtain
$N_{D_6}(D_6)=1$, 
$N_{D_6}(D_5, A_1) = 5$,  
$N_{D_6}(A_5, A_1) = 10$,  
$N_{D_6}(A_1*D_4, A_1) = 5$, 
$N_{D_6}(A_2*A_3, A_1) = 5$, 
$N_{D_6}(A_1^2*A_3, A_1) = 5$, 
$N_{D_6}(D_4, A_2) = 5$, 
$N_{D_6}(A_4, A_2) = 10$, 
$N_{D_6}(A_1*A_3, A_2) = 30$, 
$N_{D_6}(A_2^2, A_2) = 10$, 
$N_{D_6}(A_1^2*A_2, A_2) = 10$, 
$N_{D_6}(A_1^4, A_2) = 5$, 
$N_{D_6}(D_4, A_1^2) = 5$, 
$N_{D_6}(A_4, A_1^2) = 35$, 
$N_{D_6}(A_1*A_3, A_1^2) = 30$, 
$N_{D_6}(A_2^2, A_1^2) = 10$, 
$N_{D_6}(A_1^2*A_2, A_1^2) = 10$, 
$N_{D_6}(A_1^4, A_1^2) = 5$, 
$N_{D_6}(A_3, A_3) = 20$, 
$N_{D_6}(A_3,\mathbreak A_1*A_2) = 45$, 
$N_{D_6}(A_3, A_1^3) = 25$, 
$N_{D_6}(A_1*A_2, A_1*A_2) = 70$, 
$N_{D_6}(A_1*A_2, A_1^3) = 25$,
$N_{D_6}(D_4, A_1, A_1) = 25$, 
$N_{D_6}(A_4, A_1, A_1) = 100$, 
$N_{D_6}(A_1*A_3, A_1, A_1) = 150$, 
$N_{D_6}(A_2^2, A_1,\mathbreak A_1) = 50$, 
$N_{D_6}(A_1^2*A_2, A_1, A_1) = 50$, 
$N_{D_6}(A_1^4, A_1, A_1) = 25$, 
$N_{D_6}(A_3, A_2, A_1) = 125$, 
$N_{D_6}(A_3, A_1^2, A_1) = 250$, 
$N_{D_6}(A_1*A_2, A_2, A_1) = 250$, 
$N_{D_6}(A_1*A_2, A_1^2, A_1) = 375$,
$N_{D_6}(A_1^3,\mathbreak A_2, A_1) = 125$, 
$N_{D_6}(A_1^3, A_1^2, A_1) = 125$, 
$N_{D_6}(A_3, A_1, A_1, A_1) = 875$, 
$N_{D_6}(A_1*A_2, A_1, A_1,\mathbreak A_1) = 1500$, 
$N_{D_6}(A_1^3, A_1, A_1, A_1) = 625$, 
$N_{D_6}(A_2, A_2, A_2) = 100$,
$N_{D_6}(A_2, A_2, A_1^2) = 225$, 
$N_{D_6}(A_2, A_1^2, A_1^2) = 350$, 
$N_{D_6}(A_1^2, A_1^2, A_1^2) = 475$, 
$N_{D_6}(A_2, A_2, A_1, A_1) = 750$,
$N_{D_6}(A_2, A_1^2,\mathbreak A_1, A_1) = 1375$, 
$N_{D_6}(A_1^2, A_1^2, A_1, A_1) = 2000$, 
$N_{D_6}(A_2, A_1, A_1, A_1, A_1) = 5000$,
$N_{D_6}(A_1^2,\mathbreak A_1, A_1, A_1, A_1) = 8125$, 
$N_{D_6}(A_1, A_1, A_1, A_1, A_1, A_1) = 31250$.

\subsubhead The decomposition numbers for $D_7$\endsubsubhead
In addition to the described equations and assignments, we use also
the last assertion in Proposition~\PB\ for the type $A_1^5$, that is,
$N_{D_7}(A_1^5,\mathbreak A_2)=N_{D_7}(A_1^5,A_1^2)=0$.
Unfortunately, the solution space of the system of linear equations 
does not consist of a unique solution but is two-dimensional.
{\sl Mathematica~5.2} expresses all the decomposition numbers in terms of
$X=N_{D_7}(A_1^4,A_1^3)$ and
$Y=N_{D_7}(A_1^2*A_2,A_1^3)$. In the sequel we make use of the two
relations
$$\align   
  N_{D_7}(A_1^4, A_1*A_2) &= \hphantom{5}\frac65(36 - X),\tag\AT \\
  N_{D_7}(A_1^2*A_2, A_3) &= \frac3{10}(Y+162),\tag\AU\\ 
  N_{D_7}(A_1*A_3, A_1^3) &= 84 - \frac {1} {3}Y, \tag\AV\\
N_{D_7}(A_1*A_2^2, A_2) &= \frac{14}{15}(36 - 3X - Y),\tag\AW\\
N_{D_7}(A_1*A_2^2, A_1^2) &= \hphantom{7}\frac75(3X + Y-36).\tag\AX
\endalign$$
Relation (\AT) implies $X\equiv1$ (mod $5$), while
Relations~(\AU) and (\AV) imply
$Y\equiv18$ (mod~$30$). Since decomposition numbers must be
non-negative, Relations~(\AW) and (\AX) force $3X+Y$ to be equal to
$36$. Hence, we must have $X=N_{D_7}(A_1^4,A_1^3)=6$ and 
$Y=N_{D_7}(A_1^2*A_2,A_1^3)=18$.
If we substitute this in the relations for the decomposition numbers
found by {\sl Mathematica}, then we obtain
$N_{D_7}(D_7)=1$, 
$N_{D_7}(D_6, A_1) = 6$,
$N_{D_7}(A_6, A_1) = 12$,
$N_{D_7}(A_1*D_5, A_1) = 6$,
$N_{D_7}(A_2*D_4, A_1) = 6$,
$N_{D_7}(A_1^2*A_4, A_1) = 6$,
$N_{D_7}(A_3^2, A_1) = 6$,
$N_{D_7}(D_5, A_2) = 6$,
$N_{D_7}(A_5, A_2) = 12$,
$N_{D_7}(A_1*D_4, A_2) = 12$,
$N_{D_7}(A_1*A_4, A_2) = 24$,
$N_{D_7}(A_2*A_3, A_2) = 36$,
$N_{D_7}(A_1^2*A_3, A_2) = 18$,
$N_{D_7}(A_1^3*A_2, A_2) = 12$,
$N_{D_7}(D_5, A_1^2) = 9$,
$N_{D_7}(A_5, A_1^2) = 54$,
$N_{D_7}(A_1*D_4, A_1^2) = 18$,
$N_{D_7}(A_1*A_4, A_1^2) = 36$,
$N_{D_7}(A_2*A_3, A_1^2) = 54$,
$N_{D_7}(A_1^2*A_3, A_1^2) = 27$,
$N_{D_7}(A_1^3*A_2, A_1^2) = 18$,
$N_{D_7}(D_5, A_1, A_1) = 36$,
$N_{D_7}(A_5, A_1, A_1) = 144$,
$N_{D_7}(A_1*D_4, A_1, A_1) = 72$,
$N_{D_7}(A_1*A_4, A_1, A_1) = 144$,
$N_{D_7}(A_2*A_3, A_1, A_1) = 216$,
$N_{D_7}(A_1^2*A_3, A_1, A_1) = 108$,
$N_{D_7}(A_1^3*A_2, A_1, A_1) = 72$,
$N_{D_7}(D_4, A_3) = 6$,
$N_{D_7}(A_4, A_3) = 18$,
$N_{D_7}(A_1*A_3, A_3) = 72$,
$N_{D_7}(A_2^2, A_3) = 27$,
$N_{D_7}(A_1^2*A_2, A_3) = 54$,
$N_{D_7}(A_1^4, A_3) = 18$,
$N_{D_7}(D_4, A_1*A_2) = 12$,
$N_{D_7}(A_4, A_1*A_2) = 72$,
$N_{D_7}(A_1*A_3, A_1*A_2) = 180$,
$N_{D_7}(A_2^2, A_1*A_2) = 72$,
$N_{D_7}(A_1^2*A_2, A_1*A_2) = 108$,
$N_{D_7}(A_1^4, A_1*A_2) = 36$,
$N_{D_7}(D_4, A_1^3) = 2$,
$N_{D_7}(A_4, A_1^3) = 60$,
$N_{D_7}(A_1*A_3, A_1^3) = 78$,
$N_{D_7}(A_2^2, A_1^3) = 36$,
$N_{D_7}(A_1^2*A_2, A_1^3) = 18$,
$N_{D_7}(A_1^4, A_1^3) = 6$,
$N_{D_7}(D_4, A_2, A_1) = 36$,
$N_{D_7}(D_4, A_1^2, A_1) = 54$,
$N_{D_7}(A_4, A_2, A_1) = 144$,
$N_{D_7}(A_4, A_1^2, A_1) = 432$,
$N_{D_7}(A_1*A_3, A_2, A_1) = 468$,
$N_{D_7}(A_1*A_3, A_1^2, A_1) = 918$,
$N_{D_7}(A_2^2, A_2, A_1) = 180$,
$N_{D_7}(A_2^2, A_1^2, A_1) = 378$,
$N_{D_7}(A_1^2*A_2, A_2, A_1) = 324$,
$N_{D_7}(A_1^2*A_2, A_1^2, A_1) = 486$,
$N_{D_7}(A_1^4, A_2, A_1) = 108$,
$N_{D_7}(A_1^4, A_1^2, A_1) = 162$,
$N_{D_7}(D_4, A_1, A_1, A_1) = 216$,
$N_{D_7}(A_4, A_1, A_1, A_1) = 1296$,
$N_{D_7}(A_1*A_3, A_1, A_1, A_1) = 3240$,
$N_{D_7}(A_2^2, A_1, A_1, A_1) = 1296$,
$N_{D_7}(A_1^2*A_2, A_1, A_1, A_1) = 1944$,
$N_{D_7}(A_1^4, A_1, A_1, A_1) = 648$,
$N_{D_7}(A_3, A_3, A_1) = 216$,
$N_{D_7}(A_3, A_1*A_2, A_1) = 648$,
$N_{D_7}(A_3, A_1^3, A_1) = 396$,
$N_{D_7}(A_1*A_2, A_1*A_2, A_1) = 1728$,
$N_{D_7}(A_1*A_2, A_1^3, A_1) = 864$,
$N_{D_7}(A_1^3, A_1^3, A_1) = 240$,
$N_{D_7}(A_3, A_2, A_2) = 180$,
$N_{D_7}(A_1*A_2, A_2, A_2) = 576$,
$N_{D_7}(A_1^3, A_2, A_2) = 384$,
$N_{D_7}(A_3, A_2, A_1^2) = 486$,
$N_{D_7}(A_1*A_2, A_2, A_1^2) = 1296$,
$N_{D_7}(A_1^3, A_2, A_1^2) = 648$,
$N_{D_7}(A_3, A_1^2, A_1^2) = 1053$,
$N_{D_7}(A_1*A_2, A_1^2,\mathbreak A_1^2) = 2592$,
$N_{D_7}(A_1^3, A_1^2, A_1^2) = 1080$,
$N_{D_7}(A_3, A_2, A_1, A_1) = 1512$,
$N_{D_7}(A_3, A_1^2, A_1, A_1) = 3564$,
$N_{D_7}(A_1*A_2, A_2, A_1, A_1) = 4320$,
$N_{D_7}(A_1*A_2, A_1^2, A_1, A_1) = 9072$,
$N_{D_7}(A_1^3, A_2, A_1,\mathbreak A_1) = 2448$,
$N_{D_7}(A_1^3, A_1^2, A_1, A_1) = 4104$,
$N_{D_7}(A_3, A_1, A_1, A_1, A_1) = 11664$,
$N_{D_7}(A_1*A_2, A_1, A_1, A_1, A_1) = 31104$,
$N_{D_7}(A_1^3, A_1, A_1, A_1, A_1) = 15552$,
$N_{D_7}(A_2, A_2, A_2, A_1) = 1296$,
$N_{D_7}(A_2, A_2, A_1^2, A_1) = 3240$,
$N_{D_7}(A_2, A_1^2, A_1^2, A_1) = 6804$,
$N_{D_7}(A_1^2, A_1^2, A_1^2, A_1) = 13122$,
$N_{D_7}(A_2, A_2, A_1, A_1, A_1) = 10368$,
$N_{D_7}(A_2, A_1^2, A_1, A_1, A_1) = 23328$,
$N_{D_7}(A_1^2, A_1^2,\mathbreak A_1, A_1, A_1) = 46656$,
$N_{D_7}(A_2, A_1, A_1, A_1, A_1, A_1) = 77760$,
$N_{D_7}(A_1^2, A_1, A_1, A_1, A_1, A_1) = 163296$,
$N_{D_7}(A_1, A_1, A_1, A_1, A_1, A_1, A_1) = 559872$.

\subsubhead The decomposition numbers for $E_6$\endsubsubhead
These numbers have already been computed in \cite{\KratCB, Sec.~17} using 
Stembridge's {\tt coxeter} package \cite{\StemAZ}. An alternative,
independent way to do this is by following the method described in Section~4.
Here, in addition to the equations and assignments described
at the beginning of the Appendix, we use also
the last assertion in Proposition~\PB\ for the type $A_1^4$, that is,
$N_{E_6}(A_1^4,A_2)=N_{E_6}(A_1^4,A_1^2)=0$.
Unfortunately, the solution space of the system of linear equations 
does not consist of a unique solution but is one-dimensional.
{\sl Mathematica~5.2} expresses all the decomposition numbers in terms of
$X=N_{E_6}(A_1^3,A_1^3)$. In the sequel we make use of the two
relations
$$\align N_{E_6}(A_3,A_3)&=\frac {1} {100}(2592+9X),\tag\AP\\
N_{E_6}(A_1^3,A_1^3)&=\hphantom{11}\frac {1} {5}(192-6X).\tag\AS
\endalign$$
From (\AP) we infer that $X\equiv12\ (\text {mod 100})$, while (\AS)
implies that $X\le 32$. Hence, we have $X=N_{E_6}(A_1^3,A_1^3)=12$.
If we substitute this in the relations for the decomposition numbers
found by {\sl Mathematica}, then we obtain
$N_{E_6}(E_6)=1$, 
$N_{E_6}(D_5, A_1) = 12$, 
$N_{E_6}(A_5, A_1) = 6$,  
$N_{E_6}(A_1*A_4, A_1) = 12$, 
$N_{E_6}(A_1*A_2^2, A_1) = 6$,
$N_{E_6}(A_1^2*A_2, A_2) = 36$, 
$N_{E_6}(D_4, A_2) = 4$, 
$N_{E_6}(A_4, A_2) = 24$, 
$N_{E_6}(A_1*A_3, A_2) = 24$, 
$N_{E_6}(A_2^2, A_2) = 8$, 
$N_{E_6}(D_4, A_1^2) = 18$, 
$N_{E_6}(A_4, A_1^2) = 36$, 
$N_{E_6}(A_1*A_3, A_1^2) = 36$, 
$N_{E_6}(A_1^2*A_2, A_1^2) = 18$, 
$N_{E_6}(A_3, A_3) = 27$, 
$N_{E_6}(A_3, A_1*A_2) = 72$, 
$N_{E_6}(A_3, A_1^3) = 36$, 
$N_{E_6}(A_1*A_2, A_1*A_2) = 48$, 
$N_{E_6}(A_1*A_2, A_1^3) = 24$,
$N_{E_6}(A_1^3, A_1^3) = 12$, 
$N_{E_6}(D_4, A_1, A_1) = 48$, 
$N_{E_6}(A_4, A_1, A_1) = 144$, 
$N_{E_6}(A_1*A_3, A_1, A_1) = 144$, 
$N_{E_6}(A_2^2, A_1, A_1) = 24$, 
$N_{E_6}(A_1^2*A_2, A_1, A_1) = 144$, 
$N_{E_6}(A_3, A_2, A_1) = 180$, 
$N_{E_6}(A_3, A_1^2, A_1) = 378$, 
$N_{E_6}(A_1*A_2, A_2, A_1) = 336$, 
$N_{E_6}(A_1*A_2, A_1^2, A_1) = 360$,
$N_{E_6}(A_1^3, A_2, A_1) = 168$, 
$N_{E_6}(A_1^3, A_1^2, A_1) = 180$, 
$N_{E_6}(A_2, A_2, A_2) = 160$,
$N_{E_6}(A_2, A_2, A_1^2) = 288$, 
$N_{E_6}(A_2, A_1^2, A_1^2) = 504$, 
$N_{E_6}(A_1^2, A_1^2, A_1^2) = 432$, 
$N_{E_6}(A_2, A_2,\mathbreak A_1, A_1) = 1056$,
$N_{E_6}(A_2, A_1^2, A_1, A_1) = 1872$, 
$N_{E_6}(A_1^2, A_1^2, A_1, A_1) = 2376$, 
$N_{E_6}(A_3, A_1,\mathbreak A_1, A_1) = 1296$, 
$N_{E_6}(A_1*A_2, A_1, A_1, A_1) = 1728$, 
$N_{E_6}(A_1^3, A_1, A_1, A_1) = 864$, 
$N_{E_6}(A_2, A_1,\mathbreak A_1, A_1, A_1) = 6912$,
$N_{E_6}(A_1^2, A_1, A_1, A_1, A_1) = 10368$, 
$N_{E_6}(A_1, A_1, A_1, A_1, A_1, A_1) = \mathbreak 41472$.

\bigskip
\remark{Acknowledgment} 
I am indebted to Sergey Fomin for suggesting to me to embark 
on this project and,
thus, to the Institut Mittag--Leffler, Anders Bj\"orner and Richard
Stanley for hosting both of
us at the same time at the Institut during
the ``Algebraic Combinatorics" programme in Spring 2005.
Moreover, I wish to thank Christos Athanasiadis and Vic Reiner
for some helpful exchanges that I had with him on the subject matter. Last, but not least, 
I thank the referee for an extremely careful reading of the
manuscript, and for pointing out that Lemma~1.3.4 of \cite{\BesDAA}
simplifies the proof of Proposition~\PC.
\endremark
{
\eightpoint
\subhead Notes\endsubhead
After the first version of this paper was distributed,
Eleni Tzanaki found a uniform proof of Armstrong's $F=M$ Conjecture
(presented here in Conjecture~FM)
in ``Faces of generalized cluster complexes and noncrossing
partitions," {\tt ar$\chi$iv:math.CO/0605785}. Thus, our Theorem~\TA\ becomes
an unconditional theorem, that is to say, the reciprocity relation
(\AO) holds for any finite root system $\Phi$. Furthermore, 
the explicit form of the $M$-triangle in type $D_n$ is given by
\cite{\KratCB, Prop.~D} up to a simple substitution of variables.

The problem of computing the decomposition numbers for the type $B_n$
is solved implicitly in ``Enumeration of $m$-ary cacti" (Adv. Appl. Math. 
{\bf 24} (2000), 22--56) by Mikl\'os B\'ona, Michel Bousquet, Gilbert Labelle
and Pierre Leroux. This is explained in detail in the article
``Decomposition numbers for 
finite Coxeter groups and generalised non-crossing partitions" by
Thomas M\"uller and the author, which also solves the problem
for the remaining type $D_n$. In particular, the results in that
paper, together with the results in the present paper and in \cite{\KratCB},
constitute an independent --- case-by-case --- proof of the
$F=M$ Conjecture.

}

\enddocument